\documentclass{article}

\usepackage{amssymb}

\newtheorem{lem}{Lemma}[section]
\newtheorem{defi}[lem]{Definition}
\newtheorem{theo}[lem]{Theorem}
\newtheorem{prop}[lem]{Proposition}

\newtheorem{hyp}[lem]{Hypothesis}
\newtheorem{rem}[lem]{Remark}
\newtheorem{notation}[lem]{Notation}
\newtheorem{fait}[lem]{Fact}
\newtheorem{cor}[lem]{Corollary}

\newtheorem{config}[lem]{Configuration}

\newcommand{\proof}{\noindent{\bf Proof.}~}
\newcommand{\qed}{\ \hfill$\square$\bigskip}

\newcommand{\rk}{\hbox{\rm rk}\,}

\newcommand{\SL}{\hbox{\rm SL}\,}
\newcommand{\PSL}{\hbox{\rm PSL}\,}

\newcommand{\Tor}{\hbox{\rm Tor}\,}
\newcommand{\<}{\langle}

\renewcommand{\>}{\rangle}
\renewcommand{\o}{^\circ}

\newcommand{\Z}{\mathbb{Z}}
\newcommand{\N}{\mathbb{N}}

\title{Groups of finite Morley rank with solvable local subgroups}
\author{Adrien Deloro\footnote{Rutgers University} 
and Eric Jaligot\footnote{Corresponding author --- Universit\'e de Lyon, CNRS and 
Universit\'e Lyon 1}}
\date{December 25, 2007}

\begin{document}
\maketitle

\begin{abstract}
We lay down the fundations of the theory of groups of finite Morley rank 
in which local subgroups are solvable and we proceed to the local analysis of these 
groups. We prove a main Uniqueness Theorem, 
analogous to the Bender method in finite group theory, and derive its corollaries. 
We also consider homogeneous cases and study torsion. 
\end{abstract}



\newpage 

\renewcommand{\contentsname}{Contents}
\tableofcontents

\newpage
 
\section{Introduction}

In the Classification of the Finite Simple Groups \cite{GorLyonSol94}, the study of 
minimal simple groups has been a fundamental minimal case for the whole process. 
The local analysis of these finite simple groups, in which each proper subgroup 
is solvable, has been delineated by J. Thompson, originally for the Odd Order Theorem 
\cite{FeitThompson63, benderglauberman1994}. 
This work has later been used to get a classification of minimal simple 
groups in presence of elements of order $2$, and this classification has then 
been slightly generalized to the case of finite ``locally solvable" groups, 
that is finite groups in which each normalizer of a nontrivial solvable 
subgroup is also solvable. The simplicity assumption was replaced by 
a mere nonsolvability assumption. 
This full classification, with only very few 
extra groups in addition to the minimal simple ones, has been published 
in the series of papers \cite{Thompson68, Thompson70, Thompson71, Thompson73}. 

The present paper is the first of a series containing the same transfer of 
arguments from the minimal simple case to the locally solvable case in the 
context of groups of finite Morley rank. 
Indeed, a large body of work has been accomplished in the last years about 
{\em minimal connected simple} groups of finite Morley rank, that is connected 
simple groups of finite Morley rank in which every proper definable connected 
subgroup is solvable, and we propose here to transfer this work to the more 
general class of {\em locally solvable} groups of finite Morley rank, that is groups 
of finite Morley rank in which $N(A)$ is solvable for each nontrivial definable 
abelian subgroup $A$. 

As we prefer most of the time with groups of finite Morley rank 
to work in the connected category, we will actually weaken this definition of local 
solvability in the following three possible ways, by assuming solvability of the 
{\em connected components} only of normalizers of nontrivial definable abelian 
groups $A$, in which case we will use the terminology {\em locally$\o$}, 
and/or by considering nontrivial definable {\em connected} abelian 
subgroups $A$ only, in which case we will use the terminology {\em solvable$\o$}. 
In particular, we will most of the time work with the weakest definition 
of {\em local$\o$ solvability$\o$}, i.e. assuming only that 
$N\o(A)$ is solvable for each nontrivial definable connected abelian subgroup $A$ 
of the ambient group. In {\em local$\o$ solvability} we consider all 
nontrivial definable abelian (not necessarily connected) subgroups $A$. 

The only known infinite simple groups of finite Morley rank 
are algebraic groups over algebraically closed fields, and a long-standing conjecture 
postulates that there are no other such groups. Local solvability is a 
``smallness" condition and as in Thompson's final classification, 
the simplicity assumption is replaced here generally by a 
mere {\em nonsolvability} assumption. In particular, the only known 
nonsolvable connected locally$\o$ solvable groups of finite Morley rank 
are of the form $\PSL_{2}$ over some algebraically closed field $K$, and of the 
form $\SL_{2}$ in the slightly more general locally$\o$ solvable$\o$ case. 
For example, if we consider in $\SL_{3}(K)$ the definable connected abelian subgroup 
$$A={\left\{
 {\left(
      \begin{array}{llll} t & 0 & 0 \\
                                    0 & t & 0\\
                                    0 & 0 & t^{-2} 
      \end{array}
 \right)}
\mbox{ : } t \in K^{\times}
\right\}},$$
then $N\o(A)$ is a central product $A\cdot E$ where $E$ is a definable connected subgroup 
isomorphic to $\SL_{2}(K)$, so that $N\o(A)$ is not solvable. More precisely, for 
connected locally$\o$ solvable and locally$\o$ solvable$\o$ groups of finite Morley rank 
there are in the classical algebraic case no other groups than $\PSL_{2}$ and $\SL_{2}$, 
and in particular no groups of Lie rank $2$ and more. 

All the classes of locally solvable groups of finite Morley rank defined here contain 
of course all solvable groups of finite Morley rank, groups of the form 
$\PSL_{2}$ or $\SL_{2}$, but also many hypothetic configurations of 
semisimple so-called bad groups of finite Morley rank which appear as potential 
counterexamples to the main conjecture on simple groups. 

Hence, all the results of the present papers will lead to a kind of trichotomy 
(in a very large sense)  for locally$\o$ solvable$\o$ groups as follows. 
\begin{itemize}
\item
Solvable groups. 
\item
$\PSL_{2}$ or $\SL_{2}$. 
\item
Semisimple bad groups. 
\end{itemize}
In particular, the present work encapsulates the existing theory of solvable 
groups of finite Morley rank on the one hand, and of minimal connected simple 
groups on the other. 

In this first paper we are going to recast all the theory of solvable and 
minimal connected simple groups in this general context. 
In our second paper \cite{DeloroJaligotII} we are going to concentrate on the case 
of groups with involutions. Contrarily to the finite case, we cannot jump directly 
as in \cite{Thompson68} in the finite case to the case of groups with involutions, 
as no analog of the Feit-Thompson theorem is available in the context of groups 
of finite Morley rank. This is mostly due to the possible existence 
of bad groups, and we refer to \cite{Jaligot01} for the connection between 
the two problems. Our results towards algebraicity will only be partial, even 
in presence of involutions, but with a very severe limitation of nonalgebraic 
configurations. We refer to the introduction of \cite{DeloroJaligotII} for a more 
precise description of the case with involutions. 

The present paper contains a collection of results concerning the 
local analysis of locally$\o$ solvable$\o$ and locally$\o$ solvable groups 
of finite Morley rank which will be fully exploited in 
\cite{DeloroJaligotII}. That's why it also contains no theorem 
easily stated in the present introduction. The whole theory is naturally recasted in terms 
of generous Carter subgroups with the appeal of \cite{Jaligot06}. 

We will not consider 
the locally solvable/solvable$\o$ cases, which boil down rather to finite 
group theory and hence to Thompson's classification \cite{Thompson68}. 
We will however insist on the differences between 
local$\o$ solvability, which in general offers no new substantial phenomena 
compared to the minimal connected simple case, and the weaker 
local$\o$ solvabilty$\o$, where new phenomena can occur. This is 
at least explained by the alternative $\SL_{2}$ to $\PSL_{2}$. 

\medskip
\noindent
{\em Terminology.} 
A word should be said about the terminolgy adopted, as it might 
be confusing with the more classical notion of local solvability. In general 
group theory this refers usually to groups in which all finitely generated 
subgroups are solvable. In finite group theory, a subgroup is called {\em local} 
if it is the normalizer of a nontrivial $p$-subgroup for some prime $p$. 
This terminology goes back to Alperin. In \cite{Thompson68} a group 
in which each local subgroup is solvable was called an {\em $N$-group}, and 
Thompson's classification was stated for nonsolvable $N$-groups. We borrow 
the term ``local" to speak of subgroups normalizing subgroups similar to 
$p$-groups, and hence we hope 
that ``locally solvable" is clear enough in this context. We also 
note that a group of finite Morley rank in which every finitely generated 
subgroup is solvable --- the usual group theoretic notion --- must be solvable, 
and hence is locally solvable in our sense. 

\medskip
\noindent
{\em Historical remarks.}
A few historical remarks are necessary. Solvable groups of finite Morley have been 
highly investigated, notably by Nesin and Fr\'econ. As mentionned already, 
this theory becomes incorporated to the present one. 

With the ongoing work on simple groups of finite Morley rank with involutions, 
it became clear as corollaries of \cite{Jaligot99} and \cite{Jaligot01-a} that there was 
no ``small" simple groups of finite Morley rank of mixed type, and that the only 
specimen in even type was $\PSL_{2}(K)$, with $K$ an algebraically closed field 
of characteristic $2$. 

Then it was time to start the study of ``small" simple groups of odd type, even though there 
was almost nothing to start with. The fundations, notably the notion of minimal connected 
simple group, were laid down in the preprint \cite{jal3prep} 
which remained unpublished. 
It contained the first recognition of $\PSL_{2}$ in characteristic different from $2$ in 
this context, though under strong assumptions at that time. 
It also contained the embryo of local analysis of minimal connected simple groups of 
finite Morley rank. The original lemma, which turned out later to be an analog of 
the Bender method in finite group theory, was there given in any characteristic. 
It has unfortunately been disseminated between different characteristics later, and 
we will give here global forms and the general Uniqueness Theorem 
in Section \ref{SectionJaligLemma}. 

Because of the absence of a unipotence theory in characteristic zero at that time, 
and in order to reduce the size of an overambitious project to 
manageable size, the second author adopted the so-called ``tameness" assumption 
for the recognition of $\PSL_{2}$ with the weakest expectable hypothesis 
in this context. The nonalgebraic configurations were also studied in this tame 
context, and the full analysis algebraic/nonalgebraic appeared in 
\cite{CherlinJaligot2004}. 

In the meantime Cherlin suggested to develop a robust unipotence theory in characteristic 
$0$ for attacking certain problems concerning large groups of odd type without 
the tameness assumption. This became the main tool in Burdges' thesis 
\cite{JeffThesis} and this application corresponds to \cite{Burdges03}. 
This new abstract unipotence theory then allowed one to develop the local 
analysis of minimal connected simple groups where the above mentionned 
uniqueness theorem fails \cite{BurdgesJalLemma}. 
It was also Cherlin's idea to use this in presence of involutions to study other 
nonalgebraic configurations without tameness 
\cite[Case II]{BurdgesCherlinJaligot07}. 

With a nice unipotence theory then available in any characteristic, the 
recognition of $\PSL_{2}$ started again in the context of minimal connected 
simple groups of odd type without tameness, 
in the thesis of the first author \cite{DeloroThesis}. 
The recognition of $\PSL_{2}$ has then 
been obtained as in the tame case under the weakest expectable assumptions and 
appeared in \cite{Deloro05}. Using this new experience for the algebraic case, 
the nonalgebraic configurations were studied in \cite{Deloro07}, reaching 
essentially all the conclusions of \cite{CherlinJaligot2004} in the general case. 
The paper \cite{DeloroJaligotII} will at the same time improve and linearize 
the sequence of arguments represented by 
\cite{BurdgesCherlinJaligot07, Deloro05, Deloro07}, and also greatly simplify 
those in \cite{Deloro07}. 

The final generalization from minimal connected simple groups 
to locally solvable groups has been suggested by Borovik by analogy with finite 
group theory. 

\medskip
\noindent
{\em Organization of the paper.} 
Section \ref{SectionBackground} will contain background material, with notably 
an emphasis on the abstract unipotence theory in groups of finite Morley 
rank in Section \ref{SectionAbstractUnipTh} in continuation 
of \cite{Burdges0SylowPaper} and \cite{FreconJaligot07}. We shall 
formalize the notion of {\em soapy} subgroups, the finest approximation 
of unipotent subgroups where all the finest computations will be done in 
\cite{DeloroJaligotII}. 

Section \ref{SectLocSolvGps} will lay down the fundations concerning 
locally solvable groups of finite Morley rank. In 
Sections \ref{SectionSpeedLimits} and \ref{SectionL0RvsL0R0} we will 
focus on the new phenomena which can occur in the locally$\o$ solvable$\o$ 
case in comparison to the locally$\o$ solvable one. 

Section \ref{SectLocAnal} will concern the local analysis 
of locally solvable groups of finite Morley rank, with in Section \ref{SectionJaligLemma} 
the main Uniqueness Theorem (usually called ``Jaligot's Lemma") corresponding to 
the Bender method in finite group theory. The analysis of a maximal pair of Borel 
subgroups from \cite{BurdgesJalLemma}, a parallel technic, will follow in 
Section \ref{SectMaxPaire}. We will also derive consequences of the 
Uniqueness Theorem on generosity as in \cite{CherlinJaligot2004}. 

Section \ref{SectHomCasesTorsion} eventually concludes with several particular aspects 
concerning homogeneous cases as well as torsion. 
 
\medskip
\noindent
{\em Notations and background.}
For the basic background on groups of finite Morley rank we generally refer to 
\cite{BorovikNesin(Book)94}. The more recent 
\cite{AltBorCher(Book)} is also a very complete source. 
We will try to refer as much as possible to these when needed, but we 
assume the reader familiar with certain background facts such as 
Zilber's generation lemma and its corollaries \cite{Zilber77} 
\cite[\S 5.4]{BorovikNesin(Book)94}, notably 
the definability of subgroups generated by definable connected subgroups 
and corollaries on commutator subgroups. 

\begin{fait}\label{Fait[H,X]DefCon}
{\bf \cite[Corollary 5.29]{BorovikNesin(Book)94}}
Let $G$ be a group of finite Morley rank, $H$ a definable connected subgroup, and 
$X$ an arbitrary subset of $G$. Then $[H,X]$ is a definable connected subgroup of $G$. 
\end{fait}

We will also assume the reader familiar with the descending chain condition 
on definable subgroups, the existence of connected components, 
the uniqueness of generic types in connected 
groups \cite{Cherlin79}, and its immediate corollary concerning actions on finite sets. 

\begin{fait}\label{ActionOnFiniteSets}
A connected group acting definably on a finite set fixes it pointwise. 
\end{fait}

If $X$ is a subset, or a single element, of a group of finite Morley rank, we denote 
by $H(X)$ the {\em definable hull} of $X$, that is the smallest definable subgroup 
containing $X$. In the litterature it is the notation ``$d(\cdot)$" which is 
commonly used, but we prefer to keep the latter for certain integer-valued 
unipotence ``d"egrees, and instead use ``$H(\cdot)$" for ``H"ulls which are 
definable subgroups.

If $x$ and $y$ are elements of a group, we write $x^{y}$ for $y^{-1}xy$, and if 
$X$ and $Y$ are two subsets we denote by $X^{Y}$ the set of elements 
$x^{y}$. (This notation might be floppy, as we may for example use $x^{G}$ for 
the conjugacy class of $x$ in $G$.) We denote by $N(X)$ the set of elements $g$ such 
that $X^{g}=X$ (with an index if one wants to specify in which particular subset 
elements $g$ are taken).

\section{Background}\label{SectionBackground}

\subsection{Unipotence theory}\label{SectionAbstractUnipTh}

For the following abstract unipotence theory in groups of finite Morley rank 
\cite{JeffThesis, Burdges03, Burdges0SylowPaper}, 
we follow essentially the general exposition of \cite{FreconJaligot07}. 
We denote by $\cal P$ the set of all prime numbers. 

A {\em decent torus} is a 
divisible abelian group of finite Morley rank which coincides with the definable hull 
of its (divisible abelian) torsion subgroup. The latter is known to be in the finite 
Morley rank context a direct product, with $p$ varying in $\cal P$, 
of finite products of the Pr\"{u}fer $p$-group $\Z_{p^{\infty}}$ 
\cite{BorovikPoizat90}, and by divisibility decent tori are connected. 

If $p$ is a prime, a {\em $p$-unipotent} group of finite Morley rank 
is a definable connected nilpotent $p$-group of bounded exponent. 

A {\em unipotence parameter} is a pair 
$$\tilde p= (\mbox{\em characteristic~}p,\mbox{\em unipotence degree~}r)
\in{(\{ \infty \} \cup \cal{P})\times(\N \cup \{\infty\})}$$
satisfying $p<\infty$ if and only if $r=\infty$.  
A group of finite Morley rank is a {\em $\tilde p$-group} if it is nilpotent and 
of the following form, depending on the value of $\tilde p$. 
\begin{itemize}
\item
if $\tilde p=(\infty, 0)$, a decent torus. 
\item
if $\tilde p=(\infty, r)$, with $0<r<\infty$, a group generated by its definable 
indecomposable subgroups $A$ such that $A/\Phi(A)$ is torsion-free and of rank $r$. 
Here a group of finite Morley rank is {\em indecomposable} if it is abelian and 
not the sum of two proper definable subgroups. An indecomposable group $A$ 
must be connected \cite[Lemma 1.2]{Burdges0SylowPaper}, and 
$\Phi(A)$ denotes its maximal proper definable conected subgroup. 
\item
if $\tilde p=(p,\infty)$, with $p$ prime, a $p$-unipotent subgroup. 
\end{itemize}

We note that nilpotence of $\tilde p$-groups is imposed by definition, 
and that these groups are in any case generated by definable connected subgroups, 
and hence always connected by Zilber's generation lemma 
\cite[Corollary 5.28]{BorovikNesin(Book)94}. A {\em Sylow} $\tilde{p}$-subgroup 
of a group of finite Morley rank is a maximal definable (connected) $\tilde{p}$-subgroup. 

The term ``characteristic" for $p$ in a unipotence parameter $(p,r)$ 
clearly refers to the characteristic of the ground field for $p$-unipotent 
groups in algebraic groups when $p$ is finite. When $p$ is infinite and 
$0<r<\infty$, it refers to nontrivial torsion-free groups, which are 
potentially additive groups of fields of characteristic $0$. When $p$ is infinite 
and $r=0$, i.e. for decent tori, it conveys no special meaning. The term 
``unipotence degree" (one can also speak of ``weight") is explained in 
Fact \ref{ActionOnLessUnip} below by 
the constraints on actions of such groups on others. 

A group of finite Morley rank is 
{\em $(p,r)$-homogeneous} if every definable connected nilpotent subgroup 
is a $(p,r)$-group. We say that it is {\em homogeneous} if it is $(p,r)$-homogeneous 
for some unipotence parameter $(p,r)$.  
Following \cite{Cherlin05}, a divisible abelian $(\infty,0)$-homogeneous group of 
finite Morley rank is usually called a {\em good torus}. 

\begin{fait}\label{FactCharacttildepHom}
{\bf \cite[Lemma 2.17]{FreconJaligot07}}
Depending on the value of $\tilde p$, the $\tilde p$-homogene\-ity of 
a $\tilde p$-group is equivalent to the following:
\begin{itemize}
\item[$(1)$]
if $\tilde p=(\infty , 0)$, to being a good torus. 
\item[$(2)$] 
if $\tilde p=(\infty,r)$, with $0<r<\infty$, to having only 
$\tilde p$-subgroups as definable connected abelian subgroups. 
\item[$(3)$]
if $\tilde p=(p,\infty)$, with $p$ prime, then a $\tilde p$-group is 
always $\tilde p$-homogeneous. 
\end{itemize} 
\end{fait}

\begin{fait}\label{FactHomogenization}
{\bf \cite[Theorem 2.18]{FreconJaligot07}}
Let $G$ be a connected group of finite Morley rank acting definably on a 
${\tilde p}$-group $H$. Then $[G,H]$ is a definable $\tilde p$-homogeneous 
subgroup of $H$.
\end{fait}
\proof
The main point is when the unipotence degree $r$ of $H$ satisfies 
$0<r<\infty$ and is proved in \cite[Theorem 4.11]{FreconUnipotence}. 
When the unipotence degree of $H$ is infinite, this is just 
Fact \ref{FactCharacttildepHom} $(3)$. Decent tori are centralized by any 
connected group acting on them as an easy consequence of 
Fact \ref{ActionOnFiniteSets} called {\em rigidity} of 
decent tori (see Fact \ref{ActionOnLessUnip} $(1)$ below). Hence 
$[G,H]$ is trivial when $r=0$. 
\qed

\begin{cor}\label{Centrestildep}
Let $G$ be any $\tilde p$-group. Then $G^{n}$ and $G^{(n)}$ are definable 
homogeneous $\tilde p$-subgroups for any $n\geq 1$. 
\end{cor}

If $G$ is a group of finite Morley rank and 
$\tilde \pi$ is a set of unipotence parameters, we define 
$$U_{\tilde \pi}(G)={\<\Sigma~|~\tilde p \in \pi 
\mbox{~and $\Sigma$ is a definable $\tilde p$-subgroup of $G$}\>}.$$
The latter subgroup is always definable and connected by Zilber's generation lemma. 
When $\tilde \pi$ is empty it is trivial and when $\tilde \pi$ consists of a single 
unipotence parameter $\tilde p$ we simply write $U_{\tilde p}(G)$. If 
$\tilde p=(p,\infty)$ for some prime $p$, we also write $U_{p}(G)$ for 
$U_{\tilde p}(G)$. A {\em $U_{\tilde \pi}$-group} is a group $G$ such 
that $U_{\tilde{\pi}}(G)=G$. 

\begin{fait}\label{quotindec2}
{\bf \cite[Lemma 2.13]{FreconJaligot07}}
Let $f~:~G\longrightarrow H$ be a definable homomorphism between two 
groups of finite Morley rank. Then
\begin{itemize}
\item[$(1)$]
{\rm (}Push-forward\/{\rm )} $f(U_{\tilde \pi}(G))\leq U_{\tilde \pi}(H)$ is 
a $U_{\tilde \pi}$-group.
\item[$(2)$]
{\rm (}Pull-back\/{\rm )} Assume all unipotence degrees involved in $\tilde \pi$ 
are finite, or that $G$ is solvable. If $U_{\tilde \pi}(H)\leq f(G)$, then 
$f(U_{\tilde \pi}(G))=U_{\tilde \pi}(H)$.
\end{itemize}
In particular, an extension of a solvable 
$U_{\tilde \pi}$-group by a solvable $U_{\tilde \pi}$-group is a 
$U_{\tilde \pi}$-group.
\end{fait}

\begin{fait}\label{StrucNilpGroups}
{\bf \cite[\S 3]{Burdges0SylowPaper}}
Let $G$ be a nilpotent group of finite Morley rank. 
\begin{itemize}
\item[$(1)$]
$G$ is the central product of its 
Sylow $p$-subgroups and its Sylow $(\infty,r)$-subgroups. 
\item[$(2)$]
If $G$ is connected, then $G$ is the central 
product of its Sylow $\tilde p$-subgroups. 
\end{itemize}
\end{fait}
\proof
The connected case corresponds to \cite[Theorem 2.7]{FreconJaligot07}. 
Without connectedness we refer to the decomposition of $G$ as the central product of 
a definable divisible (connected) subgroup $D$ and a definable subgroup $B$ of 
bounded exponent of \cite{Nesin91-b} \cite[Theorem 6.8]{BorovikNesin(Book)94}, 
and to the decomposition of a nilpotent group of bounded exponent 
as the central product of its (definable) Sylow $p$-subgroups.  
\qed

\begin{fait}\label{pneqq}
A $\tilde p$-group of finite Morley rank 
cannot be a $\tilde q$-group when $\tilde{q} \neq \tilde{p}$. 
\end{fait}
\proof
It suffices to use the commutation provided by Fact \ref{StrucNilpGroups} $(2)$ 
to reduce the problem to abelian groups. Then it follows easily from the definitions. 
\qed

\bigskip
The following fact is a variation on the usual {\em normalizer condition} in finite 
nilpotent groups. 

\begin{fait}\label{NormCond}
{\bf (\cite[Lemma 2.4]{Burdges0SylowPaper}, 
\cite[Proposition 2.8]{FreconJaligot07})}
Let $G$ be a $\tilde p$-group and $H<G$ a proper definable subgroup. 
If $S_{1}$ and $S_{2}$ denote the Sylow $\tilde p$-subgroups of $H$ and 
of $N_{G}(H)$ respectively, then $S_{1}<S_{2}$. 
\end{fait}

\begin{fait}\label{CommuttildpGp}
{\bf \cite[Lemma 2.9]{FreconJaligot07}}
Let $G$ be a group of finite Morley rank, $S$ a subset of $G$, and 
$H$ a definable $\tilde p$-subgroup of $G$ normalized by $S$.
Then $[H,S]$ is a $\tilde p$-subgroup of $H$.
\end{fait}

\begin{fait}\label{FaitDefGpAutpGp}
Let $\tilde p$ be a unipotence parameter and $q$ a prime number. 
Let $H$ be a $\tilde p$-group of finite Morley rank without elements of 
order $q$, and assume $K$ is a definable solvable $q$-group of automorphisms of $H$ 
of bounded exponent. Then $C_{H}(K)$ is a definable $\tilde p$-subgroup of $H$.
\end{fait}
\proof
By descending chain condition on centralizers, $C_{H}(K)$ is the centralizer of a 
finitely generated subgroup of $K$, and by local finiteness of the latter we may assume 
$K$ finite. In particular $C_H(K)$ is connected by \cite[Fact 3.4]{Burdges03}. 

When $\tilde p=(\infty,0)$, $H$ is a good torus, and in particular 
$(\infty,0)$-homogeneous, and the connected subgroup $C_H(K)$ is 
also a good torus. Otherwise, $C_H(K)$ is also a $\tilde p$-group, by 
\cite[Lemma 3.18]{JeffThesis} \cite[Lemma 3.6]{Burdges03} when the 
unipotence parameter is finite or Fact \ref{FactCharacttildepHom} $(3)$ when the 
characteristic is finite.
\qed

\begin{defi}
Let $G$ be a group of finite Morley rank. 
\begin{itemize}
\item[$(1)$]
We say that $G$ {\em admits} the unipotence parameter $\tilde p$ if 
$U_{\tilde p}(G)\neq 1$. 
\item[$(2)$]
We denote by $d_{\infty}(G)$ the {\em maximal unipotence degree in 
characteristic $\infty$}, i.e. the maximal integer $r\in \N$ such that $G$ admits 
the unipotence parameter $(\infty , r)$, and $-1$ if $G$ admits none such. 
\item[$(3)$]
If $p$ is a prime, we denote by $d_{p}(G)$ the {\em maximal unipotence 
degree in characteristic $p$}, i.e. the $\infty$ symbol if $G$ admits the 
unipotence parameter $(p,\infty)$, and $-1$ otherwise.
\item[$(4)$]
A unipotence parameter $\tilde{p}=(p,r)$ is {\em maximal in its characteristic} 
for $G$ if $d_{p}(G)=r$ (notice here that the characteristic $p$ can be $\infty$ 
or prime). This is equivalent to saying that $r$ is the maximal unipotence parameter 
in characteristic $p$. 
\item[$(5)$]
Finally, we define the {\em absolute unipotence degree} $d(G)$ of $G$ as 
the maximum of $d_{\infty}(G)$ and $\max_{p\in \cal{P}} \{ d_{p}(G)\}$. 
\end{itemize}
\end{defi}

We say that a unipotence parameter $(p,r)$ is {\em absolutely maximal} for 
$G$ if $d(G)=d_{p}(G)=r$, i.e. if $G$ contains nontrivial 
$p$-unipotent subgroups if $p<\infty$ and otherwise admits $(\infty,r)$ and 
contains no nontrivial definable connected nilpotent subgroup of bounded exponent and 
no nontrivial definable $(\infty,r')$-subgroup with $r'>r$. 

We say that a unipotence parameter $(p,r)$ is {\em maximal} for $G$ if 
$d(G)=0$ whenever $r=0$, or $d_{p}(G)=r$ otherwise. This has essentially 
the effect of not considering good tori of $\PSL_{2}$ over a pure field of positive 
charateristic as having maximal unipotence degree. We will often mention this special 
example separately. 

The following lemma makes known facts more transparent in our notation. 

\begin{lem}\label{LemGenSmallMaxUnipParam}
Let $G$ be a group of finite Morley rank. 
\begin{itemize}
\item[$(1)$]
$G$ is finite if and only if $d(G)=-1$. 
\item[$(2)$]
$G$ is a good torus if and only if $G$ is connected solvable and $d(G)\leq 0$. 
\end{itemize}
\end{lem}
\proof
If $d(G)\geq 0$, then $G$ has a nontrivial definable connected 
nilpotent subgroup, and hence 
it cannot be finite. Conversely, if $G$ is infinite, then its minimal infinite definable 
subgroups are abelian by Reineke's Theorem \cite[Theorem 6.4]{BorovikNesin(Book)94}. 
As such subgroups are also connected, they contain a nontrivial Sylow $\tilde p$-subgroup 
for some unipotence parameter $\tilde p=(p,r)$ by Fact \ref{StrucNilpGroups} $(2)$, 
and hence $d(G)\geq r\geq 0 > -1$. 

If $G$ is a good torus, then it is abelian and connected, and any definable 
connected subgroup is a good torus, in particular a decent torus, 
and by Fact \ref{pneqq} $d(G)\leq 0$. Conversely, 
if $G$ is a connected solvable group which admits no unipotence parameter 
$\tilde p=(p,r)$ with $r\geq 1$, then $G$ is a good torus by 
\cite[Theorem 2.15]{Burdges03}. 
\qed

\bigskip
For any group $G$ of finite Morley rank we define, similarly to $U_{p}(G)$, 
the {\em unipotent radical in characteristic $\infty$} as 
$$U_{\infty}(G)=U_{(\infty,d_{\infty}(G))}(G).$$
One can also define the {\em absolute unipotent radical} $U(G)$ as 
$$U(G)=\<U_{p}(G)~|~p\mbox{~prime~}\>\mbox{~if it is nontrivial and~}
U_{\infty}(G)\mbox{~otherwise.}$$
Finally, a unipotent radical $U_{(p,r)}(G)$ is {\em maximal} for $G$ if 
$(p,r)$ is maximal for $G$. 

\subsection{Carter and soapy subgroups}

The preceding abstract unipotence theory in groups of finite Morley rank gives important 
approximations of semisimple and unipotent subgroups of algebraic groups. 
On the one hand it gives a good approximation of maximal tori in any group of 
finite Morley rank via the notion of Carter subgroup. 
On the other hand it detects, and it is a more difficult task, 
approximations of unipotent subgroups in locally solvable groups via the notion 
of soapy subgroups. 

All this is due to a good understanding of possible actions of $\tilde p$-subgroups onto 
each other in groups of finite Morley rank. These constraints can be 
summarized as follows. The first item is often called {\em rigidity} of decent tori. 

\begin{fait}\label{ActionOnLessUnip}
Let $G$ be a group of finite Morley rank, 
$\tilde \pi_{1}$ and $\tilde \pi_{2}$ two sets of unipotence parameters, 
and $r\in {\N \cup \{ \infty\}}$. 
\begin{itemize}
\item[$(1)$]
Assume $G=TH$ where $T$ is a definable decent torus of $G$ and $H$ 
is a definable connected 
subgroup normalizing $T$. Then $T\leq Z(G)$. In particular, if $T$ is a definable 
decent torus in a group of finite Morley rank, then $C\o(T)=N\o(T)$. 
\item[$(2)$]
Assume $G=U_{1}U_{2}$ where each $U_{i}=U_{\tilde \pi_{i}}(U_{i})$ 
is a definable nilpotent subgroup and $U_{1}$ is normal. Assume that all unipotence degrees 
involved in $\tilde \pi_{1}$ are $\leq r$ and that all unipotence degrees involved in 
$\tilde \pi_{2}$ are $\geq r$. Then $U_{1}U_{2}$ is nilpotent.
\item[$(3)$]
Assume $G=H_{1}H_{2}$ where each $H_{i}=U_{\tilde \pi_{i}}(H_{i})$ is definable 
and $H_{1}$ is normal and nilpotent. Assume that all unipotence degrees 
involved in $\tilde \pi_{1}$ are $\leq r$ and that all unipotence degrees involved in 
$\tilde \pi_{2}$ are $>r$. Then $G=H_{1}C\o(H_{1})$.
\item[$(4)$]
Assume $G=U_{1}U_{2}$ where $U_{1}$ is a normal nilpotent subgroup such 
that $U_{1}=U_{\tilde \pi_{1}}(U_{1})$, 
will all unipotence degrees involved in $\tilde \pi_{1}$ infinite, and 
$U_{2}=U_{\tilde \pi_{1}}(U_{2})$, where all unipotence degrees 
$r$ involved in $\tilde \pi_{2}$ satisfy  $0<r<\infty$. 
Then $U_{2}\leq C(U_{1})$. 
\end{itemize}
\end{fait}
\proof
The first item, which was the main key tool in \cite{Cherlin05}, 
is a mere application of Fact \ref{ActionOnFiniteSets} together with the fact 
that Pr\"ufer $p$-ranks of decent tori are finite for any prime $p$ 
\cite{BorovikPoizat90}. 

The second item is \cite[Proposition 2.10]{FreconJaligot07}. See also 
\cite[\S 3]{JaligotFrecon} and \cite[\S 4]{Burdges0SylowPaper} for earlier 
versions of the same fact. 

For the third item, we notice that if $\tilde p \in \tilde \pi_{2}$ and 
$\Sigma$ is any definable connected $\tilde p$-subgroup of $H_{2}$, then 
$H_{1}\cdot \Sigma$ is nilpotent by the second point, and both factors 
commute by our assumption on the unipotence degrees involved and 
Fact \ref{StrucNilpGroups} (2). In particular $U_{\tilde p}(H_{2})\leq C\o(H_{1})$ 
and as $H_{2}={\<U_{\tilde p}(H_{2})~|~\tilde p \in \tilde \pi_{2}\>}$, 
our claim follows. 

For the last item we refer to \cite[Lemma 4.3]{Burdges0SylowPaper} for the fact 
that an $(\infty,r)$-group, with $0<r<\infty$, which normalizes a $p$-unipotent group 
must centralize it. This is essentially a corollary of \cite[Corollary 8]{MR1833472}. 
Then one can argue as in the third point.  
\qed

\bigskip
Fact \ref{ActionOnLessUnip} has as a general consequence 
the existence of a very good approximation of semisimple subgroups of algebraic 
groups in the context of groups of finite Morley rank. 
If $\tilde \pi$ is a set of unipotence parameters, a {\em Carter $\tilde \pi$-subgroup} 
of a group of finite Morley rank is a definable connected nilpotent subgroup 
$Q_{\tilde \pi}$ such that $U_{\tilde\pi}(N(Q))=Q$. 
A {\em Carter} subgroup of a group of finite Morley rank is a definable connected 
nilpotent subgroup $Q$ such that $N\o(Q)=Q$. By Fact \ref{quotindec2} 
this corresponds to a Carter $\tilde{\pi}$-subgroup for the set $\tilde{\pi}$ 
of {\em all} unipotence parameters, or merely the set of unipotence parameters 
admitted by the ambient group. 

The existence of Carter subgroups in arbitrary groups of finite Morley rank, 
which appeared in \cite{JaligotFrecon}, has been looked for by the second 
author originally in the context of minimal connected simple groups in order to 
generalize \cite{CherlinJaligot2004}. 
It follows essentially from Fact \ref{ActionOnLessUnip}, by considering 
$\tilde p$-subgroups from the least to the most unipotent. 

\begin{fait}\label{FactExistCarter}
{\bf \cite[Theorem 3.3]{FreconJaligot07}}
Let $G$ be a group of finite Morley rank and $\tilde \pi$ a set of unipotence parameters. 
Let $r$ be the smallest unipotence degree involved in $\tilde \pi$. Then any Sylow 
$(p,r)$-subgroup of $G$ is contained in a Carter $\tilde \pi$-subgroup of $G$. 
\end{fait}

A definable subset $X$ of a group $G$ of finite Morley rank is {\em generous} 
in $G$ if the union $X^{G}$ of its $G$-conjugates is generic in $G$. In simple 
algebraic groups maximal tori are generous. In groups of finite Morley rank we only 
have equivalent conditions to this property. 

\begin{fait}\label{CarterGeneiffgendis}
{\bf \cite[Corollary 3.8]{Jaligot06}}
Let $G$ be a group of finite Morley rank and $Q$ a Carter subgroup of
$G$. Then the following are equivalent.
\begin{itemize}
\item[$(1)$]
$Q$ is generous in $G$.
\item[$(2)$]
There exists a definable generic subset $Y$ of $Q$ such that, for each
$y\in Y$, $Q$ is the unique maximal definable connected nilpotent
subgroup containing $y$. 
\item[$(3)$]
$Q$ is generically disjoint from its conjugates. 
\item[$(4)$]
There exists a definable generic subset of $Q$ all of whose elements are contained 
in only finitely many conjugates of $Q$. 
\end{itemize} 
\end{fait}

At the opposite of semisimple groups, we pass now to the approximations of 
unipotent subgroups. 
We denote by $F(G)$ the {\em Fitting} subgroup of any group $G$, i.e. 
the subgroup generated by all normal nilpotent subgroups. It is always 
definable and nilpotent in the finite Morley rank case 
\cite[Theorem 7.3]{BorovikNesin(Book)94}. 
A consequence of Fact \ref{ActionOnLessUnip} dual to 
Fact \ref{FactExistCarter} is the following. 

\begin{fait}\label{FactUnipHeaviest}
Let $H$ be a connected solvable group of finite Morley rank and $\tilde{p}=(p,r)$ 
a unipotence parameter with $r>0$. Assume $d_{p}(H)\leq r$. Then 
$U_{\tilde p}(H)\leq F\o(H)$. 
\end{fait}
\proof
See \cite[Lemma 2.11]{FreconJaligot07}, 
and \cite[Theorem 2.16]{Burdges03} for the original version. 
It suffices to use Fact \ref{ActionOnLessUnip} $(2)$ and $(4)$ to conclude 
that $F\o(H) \cdot U_{\tilde p}(H)$ is nilpotent, and then to use 
the fact that $H/F\o(H)$ is abelian (Fact \ref{GConSolvG/FGdivAb} below). 
\qed

\bigskip
We note that the assumption $r>0$ is necessary in Fact \ref{FactUnipHeaviest}. 
In the standard Borel subgroup $B$ of $\PSL_{2}$ in positive characteristic, 
$d_{\infty}(B)=0$, but maximal tori of $B$ are not in the unipotent radical of $B$. 

Unipotent subgroups are usually not generous in linear algebraic groups, 
and thus in general more difficult to detect. Every {\em nontrivial} subgroup 
$U_{\tilde p}(H)$ as in Fact \ref{FactUnipHeaviest} is generally a 
good approximation of unipotent radical, at least much finer than the Fitting 
subgroup. We will need even finer approximations when considering 
locally solvable groups of finite Morley rank, notably the property of being homogeneous 
and central in the Fitting subgroup. This issues from the minimal 
subgroups used originally in \cite{jal3prep}, after the 
considerable reworking in \cite{Deloro05, Deloro07}. 

Recall that, for every connected solvable group $H$ of finite Morley rank, a unipotence 
parameter $\tilde{q}=(q,d)$ is {\em maximal} for $H$ if 
$d(H)=0$ whenever $d=0$, or $d_{q}(H)=d$ otherwise. By 
Lemma \ref{LemGenSmallMaxUnipParam}, a nontrivial connected solvable 
group $H$ is a good torus if and only if its unique maximal unipotence parameter 
is $(\infty,0)$. Otherwise, 
maximal unipotence parameters are all the $(p,\infty)$ such that $U_{p}(H)\neq 1$ and 
the $(\infty,d)$ with $d\geq 1$ and $d_{\infty}(H)=d$ if it exists. 

\begin{defi}
Let $H$ be a connected solvable group of finite Morley rank. A subgroup $U$ of 
$H$ is {\em soapy} (resp. {\em characteristically soapy}) in $H$ if the two following 
conditions hold. 
\begin{itemize}
\item[$(1)$]
$U$ is a nontrivial definable connected subgroup of $Z(F\o(H))$, 
$\tilde{q}$-homoge\-neous for some unipotence parameter $\tilde q$ 
maximal for $H$. 
\item[$(2)$]
$U$ is normal (resp. definably characteristic) in $H$.
\end{itemize}
\end{defi}

We haven't found a better name for these subgroups. 
We will see in Section \ref{SectionConsOnSoapySubgp} that 
in locally$\o$ solvable$\o$ groups these subgroups have 
a strong tendency to escape from intersections of distincts Borel subgroups, like 
unipotent subgroups in $\PSL_{2}$ and like a soap between to hands. 
Another not less serious reason for this name is that these groups were born 
near Marseilles, which is famous for its soap. 

We could also specify a set of maximal unipotence parameters for $H$, 
and define these interesting subgroups as products of the present ones. In practice 
only one unipotence parameter will suffice for us. 

The next lemma says that the existence of soapy subgroups is not essentially 
weaker than that of characteristically soapy subgroups. 

\begin{lem}\label{LemFromSoapyToCharSoapy}
Let $H$ be a connected solvable group of finite Morley rank and $\tilde q$ a 
unipotence parameter maximal for $H$. If $H$ contains a 
$\tilde{q}$-homogeneous soapy subgroup, 
then it contains a $\tilde q$-homogeneous characteristically soapy subgroup as well. 
\end{lem}
\proof
If $\tilde q=(\infty,0)$ then $H$ is a good torus, and $H$ itself is the desired group. 

In general one can proceed as follows. 
Let $U$ be a $\tilde q$-homogeneous soapy subgroup of $H$. Let 
$\tilde U$ be the subgroup of $Z(F\o(H))$ generated by all 
$\tilde q$-homogeneous soapy subgroups of $H$. It is nontrivial, definable and 
connected as the product of finitely many soapy subgroups by Zilber's generation lemma, 
and one sees easily that it is $\tilde q$-homogeneous with Fact \ref{quotindec2} 
(see also \cite[Corollary 3.5]{FreconUnipotence}). 
It is clearly definably characteristic in $H$. 
Hence $\tilde U$ is characteristically soapy in $H$. 
\qed

\bigskip
We finish this section with a general criterion for building 
characteristically soapy subgroups. 

\begin{lem}\label{getwise1}
Let $H$ be a connected solvable group of finite Morley rank and $\tilde q$ a unipotence 
parameter maximal for $H$. If $U_{\tilde q}(Z(F\o(H)))$ is not central in 
$H$, then $H$ contains a $\tilde q$-homogeneous characteristically soapy subgroup. 
\end{lem}
\proof
Set $U=[U_{\tilde{q}}(Z(F\o(H))),H]$. By assumption $U$ is nontrivial. 
It is a definable connected homogeneous $\tilde{q}$-subgroup by 
Fact \ref{FactHomogenization}, contained in $Z(F\o(H))$ as the latter 
is normal in $H$, and obviously definably characteristic in $H$. 
\qed

\subsection{Conjugacy theorems}

As far as unipotence theory in concerned, there are two general conjugacy theorems 
in groups of finite Morley rank. The first one has a nontrivial content only 
in presence of divisible torsion. 

\begin{fait}\label{ConjdecentTori}
{\bf \cite{Cherlin05}}
Let $G$ be a group of finite Morley rank. Then $C\o(T)$ is generous in $G\o$ 
for every definable decent torus $T$ of $G\o$, and maximal definable decent 
tori of $G\o$ are $G\o$-conjugate.  
\end{fait}

The following corollary of Fact \ref{ConjdecentTori} 
has been known for a long time in presence of $2$-divisible 
torsion \cite[Lemma 10.22]{BorovikNesin(Book)94}. 

\begin{cor}[Control of fusion]\label{CorControlFusion}
Let $G$ be a group of finite Morley rank, $p$ a prime, and $T$ a $p$-torus of $G$. 
If $X$ and $Y$ are two $G$-conjugate subsets of $G$ such that 
$C_{T}(X)$, $C_{T}(Y)$, and $C(Y)$ all have the same Pr\"{u}fer $p$-ranks, then 
$Y=X^{g}$ for some $g$ conjugating $C\o_{T}(X)$ to $C\o_{T}(Y)$. 
In particular if $T$ is a maximal $p$-torus of $G$ then any two 
$G$-conjugate subsets of $C(T)$ are $N(T)$-conjugate. 
\end{cor}
\proof
First notice that there are always {\em maximal} $p$-tori, by 
finiteness of the Pr\"ufer $p$-rank \cite{BorovikPoizat90} 
and compactness. 

Assume $Y=X^{g}$ for some $g\in G$. Then $C\o_{T}(X)^{g}$ and 
$C\o_{T}(Y)$ are both contained in 
the definable subgroup $C\o(Y)$. By Fact \ref{ConjdecentTori} and the assumption, 
$C\o_{T}(X)^{g}=C\o_{T}(Y)^{\gamma}$ for some $\gamma \in C\o(Y)$. 
Then $g\gamma^{-1}$ conjugates $C\o_{T}(X)$ to $C\o_{T}(Y)$ 
and as $Y^{\gamma}=Y=X^{g}$ the 
element $g\gamma^{-1}$ conjugates $X$ to $Y$. 

When $X$ and $Y$ are two subsets of $C(T)$ and $T$ is maximal 
we can apply the preceding and 
the new element $g$ conjugating $X$ to $Y$ will now normalize $T$. 
\qed

\bigskip
There is no reason why an arbitrary group of finite Morley rank should 
contain nontrivial torsion as in Fact \ref{ConjdecentTori}. However 
the next general conjugacy theorem relies on an assumption which is likely 
to be true in general \cite[\S4]{Jaligot06}. 

\begin{fait}\label{FactConjGenCarter}
{\bf \cite{Jaligot06}}
Let $G$ be a group of finite Morley rank. Then generous Carter subgroups of 
$G$ are generous in $G\o$ and $G\o$-conjugate. 
\end{fait}

In our study of locally solvable groups of finite Morley rank we will 
of course use much more conjugacy theorems where they are much more aboundant, 
that is in solvable groups. 

\subsection{Solvable groups}

\begin{fait}\label{GConSolvG/FGdivAb}
{\bf \cite{Nesin90-b}}
Let $H$ be a connected solvable group of finite Morley rank. 
Then $H/F\o(H)$ is divisible abelian. 
\end{fait}

\begin{fait}\label{HallGpConInConSolv}
{\bf (\cite[Corollaire 7.15]{Frecon00}, \cite{BorovikNesin92})}
Let $H$ be a connected solvable group of finite Morley rank, and $\pi$ any set 
of prime numbers. Then Hall $\pi$-subgroups of $H$ are connected. 
\end{fait}

\begin{fait}\label{StructpSylSolvGps}
{\bf \cite{BorovikPoizat90}}
Let $p$ be a prime and $S$ a $p$-subgroup of a solvable group of finite Morley rank, 
or more generally a locally finite $p$-subgroup of any group of finite Morley rank. Then 
\begin{itemize}
\item[$(1)$]
$S\o$ is a central product of a $p$-torus and a $p$-unipotent subgroup.
\item[$(2)$]
If $S$ is infinite and of bounded exponent, then $Z(S)$ contains infinitely many 
elements of order $p$. 
\end{itemize}
\end{fait}

\begin{lem}\label{FacptTorFitUpPerp}
Let $H$ be a connected solvable group of finite Morley rank and $p$ a prime. 
If $U_{p}(H)=1$, then the Sylow $p$-subgroup of $F(H)$ is central in $H$. 
\end{lem}
\proof
Assume $U_{p}(H)=1$, and let $S$ denote the Sylow $p$-subgroup of $F(H)$. 
By Fact \ref{StrucNilpGroups}, $S$ is the product of a finite 
$p$-subgroup and of a $p$-torus. As each of these two subgroups 
is normal in $H$, each is central in $H$, by 
Facts \ref{ActionOnFiniteSets} and \ref{ActionOnLessUnip} $(1)$ respectively.  
\qed

\bigskip
The following fact gradually appeared in \cite{Wagner94}, \cite{Frecon00}, 
and \cite[3.5]{CherlinJaligot2004}. 

\begin{fait}\label{Carterconjugateandselfnorm}
{\bf \cite[Theorem 3.11]{FreconJaligot07}}
Let $H$ be a connected solvable group of finite Morley rank. 
Then Carter subgroups of $H$ are generous, conjugate and self-normalizing.
\end{fait}

\begin{cor}\label{ConseCarterSelfNorm}
Let $G$ be a group of finite Morley rank, $Q$ a Carter subgroup and $\sigma$ an 
element normalizing $Q$ and not in $Q$. Then $\sigma \notin C\o(X)$ 
for every $X\subseteq Z(Q)$ such that $C\o(X)$ is solvable. In 
particular such $\sigma$ and $X$ cannot be in the same definable connected 
abelian subgroup. 
\end{cor}
\proof
Assuming the contrary, then $\sigma \in N_{C\o(X)}(Q)=Q$ by 
the selfnormalization given in Fact \ref{Carterconjugateandselfnorm}, a contradiction. 
For the second point we simply notice that otherwise $\sigma \in C\o(X)$. 
\qed

Following \cite[\S4-5]{FreconJaligot07} there are nice links between 
Carter $\tilde{\pi}$-groups and covering properties in connected solvable 
groups of finite Morley rank, the so-called connected subformation theory.  
In particular one knows that the collection $\cal N$ of connected nilpotent 
groups of finite Morley rank is a connected subformation. 
The main link between Carter subgroup theory 
and subformation theory in connected solvable groups is then a guarantee that 
Carter subgroups of a connected solvable group $G$ of finite Morley rank are 
$\cal N$-covering subgroups of $G$, which provides the following important result.

\begin{fait}\label{Carter=NProj}
{\bf \cite[Proposition 5.1]{FreconJaligot07}}
In any connected solvable group of finite Morley rank, Carter subgroups are exactly 
the $\cal N$-covering subgroups and the $\cal N$-projectors.
\end{fait}

Here the properties of $\cal N$-covering subgroups and $\cal N$-projectors which 
interest us are that these groups cover all nilpotent connected sections containing them. 

\begin{fait}\label{TheoConjPiCarter}
{\bf \cite[Theorem 5.8]{FreconJaligot07}}
Let $G$ be a connected solvable group of finite Morley rank and 
$\tilde \pi$ a set of unipotence parameters. Then 
Carter $\tilde \pi$-subgroups are exactly 
the ${\cal N}_{\tilde \pi}$-projectors and the ${\cal N}_{\tilde \pi}$-covering 
subgroups of $G$, and are in particular conjugate. 
\end{fait}

We note that when $\tilde{\pi}$ is a single unipotence parameter, 
Carter ${\tilde{p}}$-subgroup coincide with Sylow $\tilde p$-subgroups 
\cite[\S 3.2]{FreconJaligot07}, so that Sylow $\tilde p$-subgroups are conjugate 
in connected solvable groups of finite Morley rank. There is also structural information 
concerning Carter $\tilde{\pi}$-subgroups of connected solvable groups 
of finite Morley rank \cite[Corollary 5.9]{FreconJaligot07}, and we will use this only 
with $\tilde{\pi}=\{ \tilde{p} \}$. 

\begin{fait}\label{StrucSylowUSolv}
{\bf (\cite[Corollary 5.11]{FreconJaligot07}, 
\cite[Theorem 6.7]{Burdges0SylowPaper})}
Let $G$ be a connected solvable group of finite Morley rank. Then 
the Sylow $\tilde p$-subgroups of $G$ are exactly the subgroups of the form 
$U_{\tilde p}(G')U_{\tilde p}(Q)$ for some Carter subgroup $Q$ of $G$. 
\end{fait}

If $G$ is a group of finite Morley rank, we denote by 
$$O_{p'}(H)$$
the largest normal definable connected subgroup without $p$-torsion. It exists by ascending 
chain condition on definable connected subgroups and elementary properties 
of lifting of torsion \cite{BorovikNesin92}. 

The following facts will be useful when dealing with $p$-strongly embedded subgroups 
in Section \ref{SectionPruferRanks} below. 

\begin{fait}\label{FaitH/Op'HDivAbelian}
{\bf (Compare with \cite[Lemma 3.2]{CherlinJaligot2004})}
Let $H$ be a connected solvable group of finite Morley rank 
such that $U_{p}(H)=1$. Then $H/O_{p'}(H)$ is divisible abelian. 
\end{fait}
\proof
Dividing by $O_{p'}(H)$, we may assume it is trivial and we want 
to show that $H$ is divisible abelian. 

Let $F=F\o(H)$. As $O_{p'}(H)=1$, $O_{p'}(F)=1$ as well, 
and $U_{q}(H)=1$ for any prime $q$ different from $p$. By assumption 
$U_{p}(H)=1$ also, and $F$ is divisible by Fact \ref{StrucNilpGroups}. 
As $F'$ is torsion-free, by \cite[Theorem 2.9]{BorovikNesin(Book)94} 
or Fact \ref{StrucNilpGroups} and Corollary \ref{Centrestildep}, 
it must be trivial by assumption. Hence $F$ is divisible abelian. 

To conclude it suffices to show that $F$ is central in $H$, as then $H$ 
is nilpotent, hence equal to $F$, and hence divisible abelian, as desired. 
Let $h$ be any element of $H$; we want to show that $[h,F]=1$. But 
$[h,F]$ is torsion-free, as the torsion subgroup of $F$ is central 
in $H$ by Fact \ref{ActionOnLessUnip} $(1)$, or using 
Fact \ref{FactHomogenization}. Hence $[h,F]\leq O_{p'}(H)=1$, as desired. 
\qed

\begin{fait}\label{BiGenerationWithpElts}
{\bf \cite[Fact 3.7]{Burdges03}}
Let $H$ be a solvable group of finite Morley rank without elements of order $p$ 
for some prime $p$. 
Let $E$ be a finite elementary abelian $p$-group acting definably on $H$. 
Then 
$$H=\<C_{H}(E_{0})~|~E_{0}\leq E,~[E:E_{0}]=p\>.$$
\end{fait}

\begin{lem}\label{Bipgeneration}
Let $H$ be a connected solvable group of finite Morley rank 
such that $U_{p}(H)=1$ for some prime $p$. Assume $H$ contains an elementary abelian 
$p$-group $E$ of order $p^{2}$. Then 
$$H=\<C\o_{H}(E_{0})~|~E_{0}
\mbox{~is a cyclic subgroup of order $p$ of~}E\>.$$
\end{lem}
\proof
By assumption and Facts \ref{HallGpConInConSolv} and \ref{StructpSylSolvGps}, 
Sylow $p$-subgroups of $H$ are $p$-tori. Hence $E$ is in a maximal $p$-torus of $H$, 
which is included in a Carter subgroup $Q$ of $H$ by Fact \ref{FactExistCarter}. 
By Fact \ref{FaitH/Op'HDivAbelian}, 
$H/O_{p'}(H)$ is abelian. As Carter subgroups cover all abelian quotients 
in connected solvable groups of finite Morley rank by Fact \ref{Carter=NProj}, 
$H=O_{p'}(H)\cdot Q$. As $E\leq Z(Q)$, it suffices to show that 
$$O_{p'}(H)=\<C\o_{O_{p'}(H)}(E_{0})~|~E_{0}
\mbox{~is a cyclic subgroup of order $p$ of~}E\>.$$
But the generation by the full centralizers is given by 
Fact \ref{BiGenerationWithpElts}, and these centralizers are connected 
by \cite[Fact 3.4]{Burdges03}. 
\qed

\begin{cor}\label{CorBipgeneration}
Let $H$ be a connected solvable group of finite Morley rank with a toral 
subgroup $E$ of order $p^2$ for some prime $p$. Then 
$$H=\<C\o_{H}(E_{0})~|~E_{0}
\mbox{~is a cyclic subgroup of order $p$ of~}E\>.$$
\end{cor}
\proof
For a connected nilpotent group of finite Morley rank $L$, we define 
the ``complement" $C_{p}(L)$ of $U_{p}(L)$, namely the product 
of all factors of $L$ as in Fact \ref{StrucNilpGroups} $(2)$, except $U_p(L)$. 

Now if $H$ is any connected solvable group of finite Morley rank and 
$Q$ a Carter subgroup of $H$, then $H$ is the product of the definable 
connected subgroup $C_{p}(Q)C_{p}(F\o(H))$ with the normal definable 
connected subgroup $U_{p}(H)$, and the first factor has trivial 
$p$-unipotent subgroups. 

In our particular case, $E$ is by torality contained in a $p$-torus, and the latter 
is contained in a Carter subgroup $Q$ of $H$. By 
Facts \ref{HallGpConInConSolv} and \ref{StructpSylSolvGps} $E$ 
centralizes the normal definable connected subgroup $U_{p}(H)$, 
so it suffices to show the generation by centralizers$\o$ in 
$C_{p}(Q)C_{p}(F\o(H))$. But this follows from 
Lemma \ref{Bipgeneration}. 
\qed

\subsection{Genericity}

\begin{lem}\label{HUnifCovFinHBounded}
Let $H$ be a connected solvable group of finite Morley rank generically 
covered by a uniformly definable family of finite subgroups. 
Then $H$ is nilpotent and of bounded exponent. 
\end{lem}
\proof
We first note that any group generically covered by a uniformly definable 
family of finite groups is generically of bounded exponent. In fact, 
by elimination of infinite quantifiers \cite[Proposition 2.2]{MR1955252}, 
there is a uniform bound on the cardinals of the finite groups involved. 

Now $H/F\o(H)$ is divisible abelian by Fact \ref{GConSolvG/FGdivAb}. 
As Pr\"ufer $p$-ranks are finite for each prime $p$, 
there is a finite subgroup of $H/F\o(H)$ containing all images modulo $F\o(H)$ 
of the finite groups. This shows by generic covering that $H/F\o(H)$ is trivial. 
Hence $H$ is nilpotent. Now it suffices to use the generic covering again 
and Fact \ref{StrucNilpGroups} $(2)$ with Fact \ref{quotindec2} $(1)$. 
\qed

\bigskip
The following lemma has its roots in \cite[Lemme 2.13]{jal3prep} 
(see \cite[Fact 2.36]{CherlinJaligot2004}). 

\begin{lem}\label{LemGenXcapMGenInXGInv}
Let $G$ be a connected group of finite Morley rank and $X$ a non\-empty 
definable $G$-invariant subset of $G$. If $M$ is a definable subgroup of $G$ 
such that $X\cap M$ is generic in $X$, then $X\cap M$ contains a definable 
$G$-invariant subset generic in $X$.  
\end{lem}
\proof
By assumption $X$ is a union of $G$-conjugacy classes. By assumption also, 
$X\cap M$ is nonempty. 

Let $Y_{1}$ be a definable generic subset of $X\cap M$ consisting of elements of 
$X\cap M$ whose $G$-conjugacy classes have traces on $X\cap M$ of 
constant ranks. Let $Y_{2}$ be a definable generic subset of $Y_{1}$ consisting of 
elements of $Y_{1}$ whose $G$-conjugacy classes in $G$ have constant 
ranks. Both exist as we have, by definability of the rank,  
finite definable partitions in each case. 
Now $Y_{2}$ is generic in $Y_{1}$ which is generic in $X\cap M$, so $Y_{2}$ is 
generic in $X\cap M$ and in $X$. Replacing $X$ by $Y_{2}^{G}$, one can 
thus assume that $G$-conjugacy classes in $G$ of elements of $X$, as well as their 
traces on $M$, are of constant ranks. We also have then that $x^{G}\cap M$ is 
nonempty for any $x$ in $X$. 

Now, as $X$ is the union of the $G$-conjugacy classes of its elements in $X\cap M$ 
and reduced to the situation where all relevant ranks are constant, 
the assumption that $X\cap M$ is generic in $X$ implies easily by additivity 
of the rank that $x^{G}\cap M$ is generic in $x^{G}$ for any $x$ in $X$. 

Let $N=\bigcap_{g\in G}M^{g}$. By descending chain condition on definable 
subgroups, $N=M^{g_{1}}\cap \cdots \cap M^{g_{n}}$ 
for finitely many elements $g_{1}$, ..., $g_{n}$ of $G$. 
As $G$ is connected, $x^{G}$, which is in definable bijection with 
$G/C(x)$, has Morley degree $1$ for any $x$ in $X$. 
By taking conjugates one also has $x^{G}\cap M^{g_{i}}$ 
generic in $x^{G}$ for each $x$ in $X$ and each $g_{i}$. 
Hence $x^{G}\cap N$, which 
can be written as 
$$(x^{G}\cap M^{g_{1}}) \cap \cdots \cap (x^{G}\cap M^{g_{n}}),$$ 
is also generic in $x^{G}$, for any $x$ in $X$. Now the fact that all ranks 
involved are constant implies that $X\cap N$ is generic in $X$ as well. 

But $X\cap N$ is $G$-invariant as both sets involved are. Hence $X\cap N$ 
is the desired definable $G$-invariant subet of $X\cap M$ generic in $X$. 
\qed

\section{Locally solvable groups}\label{SectLocSolvGps}

\subsection{Fundations}\label{SectFund}

\begin{defi}\label{DefLocalSolv}
We say that a group of finite Morley rank is 
\begin{itemize}
\item[$(1)$]
{\em locally solvable} if $N(A)$ is solvable for each nontrivial definable abelian 
subgroup $A$. 
\item[$(2)$]
{\em locally solvable$\o$} if $N(A)$ is solvable for each nontrivial definable 
abelian connected subgroup $A$. 
\item[$(3)$]
{\em locally$\o$ solvable} if $N\o(A)$ is solvable for each nontrivial definable 
abelian subgroup $A$. 
\item[$(4)$]
{\em locally$\o$ solvable$\o$} if $N\o(A)$ is solvable for each nontrivial definable 
abelian connected subgroup $A$. 
\end{itemize}
\end{defi}

\begin{lem}\label{FirstRemOnLocSolvGps}
Let $G$ be a group of finite Morley rank. 
\begin{itemize}
\item[$(1)$]
If $G$ satisfies one of the Definitions \ref{DefLocalSolv} $(1)$, $(2)$, $(3)$, or $(4)$, 
then so does any definable subgroup of $G$. 
\item[$(2)$]
If $G$ is locally solvable, then is it locally solvable$\o$ and locally$\o$ solvable, and if 
$G$ has any of the two latter properties, then it is locally$\o$ solvable$\o$. 
\end{itemize}
\end{lem}
\proof
Obvious. 
\qed

\begin{defi}
Let $G$ be a group of finite Morley rank and $H$ a subgroup of $G$. We say that 
a subgroup $L$ of $G$ is 
\begin{itemize}
\item[$(1)$]
{\em $H$-local} if $L\leq N(H)$. 
\item[$(2)$]
{\em $H$-local$\o$} if $L\leq N\o(H)$. 
\end{itemize}
\end{defi}

Then we say that a subgroup $L$ is 
{\em local} if it is $H$-local for some subgroup $H$, and 
{\em local$\o$} if it is $H$-local$\o$ for some subgroup $H$. 
We can give conditions a priori stronger, but actually equivalent, to 
Definitions \ref{DefLocalSolv} $(1)$--$(4)$ in terms of local subgroups. 

\begin{lem}\label{LemGenCharLocSolvGps}
Let $G$ be a group of finite Morley rank. Then $G$ is 
\begin{itemize}
\item[$(1)$]
locally solvable if and only if $X$-local subgroups subgroups are solvable for every 
nontrivial solvable subgroup $X$. 
\item[$(2)$]
locally solvable$\o$ if and only if $X$-local subgroups are solvable for every 
infinite solvable subgroup $X$. 
\item[$(3)$]
locally$\o$ solvable if and only if $X$-local$\o$ subgroups are solvable for every 
nontrivial solvable subgroup $X$. 
\item[$(4)$]
locally$\o$ solvable$\o$ if and only if $X$-local$\o$ subgroups are solvable for every 
infinite solvable subgroup $X$. 
\end{itemize}
\end{lem}
\proof
Clearly the right conditions are stronger than the left ones. 

Assume now a left condition, and suppose $X$ is some nontrivial solvable 
subgroup of $G$, and $L$ is an $X$-local subgroup, i.e. $L\leq N(X)$. 
Then $L$ normalizes the definable hull $H(X)$ of $X$, and its connected component 
$H\o(X)$ as well. Now a classical corollary of Zilber's generation lemma on 
derived subgroups (Fact \ref{Fait[H,X]DefCon}) implies that the last nontrivial term 
of the derived series of $H(X)$, as well as $H\o(X)$, is definable . 
It is abelian by definition, and as it is 
characteristic in $H(X)$ (resp. $H\o(X)$), it is normalized by $L$. Then ones 
sees in each case which has to be considered that the latter is solvable by the 
left condition. 
\qed

\bigskip
Nontrivial solvable groups $H$ of finite Morley rank contain certain nontrivial definable 
characteristic Sylow $\tilde p$-subgroups or Sylow $p$-subgroups, by 
Fact \ref{StrucNilpGroups} applied in $F(H)$. Hence $H$-local subgroups normalize 
nontrivial $\tilde p$-groups or $p$-groups, so that our definitions are 
coherent with the notion due to Alperin of local subgroup in finite group theory 
\cite{Thompson68}, as subgroups normalizing nontrivial $p$-subgroups. 
Before stating this a little bit more precisely in the locally$\o$ solvable$\o$ case, we 
look at quotients. 

\begin{lem}
Let $G$ be a group of finite Morley rank and $N$ a definable normal solvable subgroup. 
\begin{itemize}
\item[$(1)$]
If $G$ is locally solvable, then so is $G/N$. 
\item[$(2)$]
If $G$ is locally solvable$\o$, then so is $G/N$. 
\item[$(3)$]
If $G$ is locally$\o$ solvable, then so is $G/N$. 
\item[$(4)$]
If $G$ is locally$\o$ solvable$\o$, then so is $G/N$. 
\end{itemize}
\end{lem}
\proof
We denote by $\overline G$ the quotients by $N$. 

$(1)$. Let $\overline A$ be a nontrivial definable abelian subgroup of $\overline G$. 
The preimage of $N_{\overline G}(\overline A)$ normalizes $AN$, which is solvable 
and nontrivial, and hence it is solvable by local solvability of $G$. As $N$ is solvable, 
$N_{\overline G}(\overline A)$ is also solvable. 

$(2)$. One can proceed as in $(1)$, taking $A$ infinite modulo $N$, 
and looking at the normalizer of $(AN)\o$. 

$(3)$. One can proceed as in $(1)$, taking connected components of normalizers 
throughout. 

$(4)$. It suffices to mix the two preceding cases. 
\qed

\bigskip
We continue with trivial remarks. In a group of finite Morley rank we call 
{\em Borel} subgroup any maximal definable connected solvable subgroup. 

\begin{lem}\label{CharactBorelXLocal}
Let $G$ be a locally$\o$ solvable$\o$ group of finite Morley rank. 
Then a subgroup $B$ is a Borel subgroup if and only if $B$ is a 
maximal $X$-local$\o$ subgroup for some infinite solvable subgroup $X$. 
Furthermore $X$ can be chosen to be any, and has to be an 
infinite normal subgroup of $B$. 
\end{lem}
\proof
Let $B$ be a Borel subgroup of $G$. Then $B\leq N\o(X)\leq N\o(H(X))$ 
for any infinite normal subgroup $X$ of $B$, and $N\o(H(X))$ is solvable by 
local$\o$ solvability$\o$ of $G$. Hence we get equality by maximality 
of $B$, and hence $B$ is a connected $X$-local$\o$ subgroup. 
If $B$ is contained in a $Y$-local$\o$ subgroup $L$ of $G$ for some infinite solvable 
subgroup $Y$, then $B\leq N\o(H(Y))$ and as $N\o(H(Y))$ is solvable 
by local$\o$ solvability$\o$ one gets $B=N\o(H(Y))$ again by maximality 
of $B$. As $B\leq L\leq N\o(H(Y))$, $B=L$. 

Let now $B$ be a maximal $X$-local$\o$ subgroup of $G$ for some infinite 
solvable subgroup $X$. By local$\o$ solvability$\o$ $B$ is contained in a Borel subgroup 
$B_{1}$. Now $B_{1}$ is $Y$-local$\o$ for some infinite solvable 
subgroup $Y$, and the maximality of $B$ implies that $B=B_{1}$. 

Now if a Borel subgroup $B$ normalizes an infinite solvable subgroup $X$, 
then $X\cdot B$ is solvable, as well as its definable hull, and by maximality 
$H\o(X)\leq B$, and $H\o(X)$ is an infinite normal subgroup of $B$. 
\qed

\begin{lem}\label{LemEquivCondL0R0NotSolv}
Let $G$ be a locally$\o$ solvable$\o$ group of finite Morley rank. 
Then the following are equivalent. 
\begin{itemize}
\item[$(1)$]
$N\o(A)<G\o$ for each nontrivial definable connected abelian subgroup $A$. 
\item[$(2)$]
$G\o$ is not solvable. 
\item[$(3)$]
$G\o$ has two distinct Borel subgroups. 
\end{itemize}
\end{lem}
\proof
If $G\o$ is solvable, then $G\o\leq N\o(A)$ where $A$ is the last 
nontrivial term of the derived series of $G$, which is definable and connected by 
Fact \ref{Fait[H,X]DefCon}. Hence the first condition implies the second one. 

If $G\o$ has two distinct Borel subgroups, then clearly $G\o$ cannot be solvable. 

Finally, assume $G\o=N\o(A)$ for some nontrivial definable connected abelian 
subgroup $A$. 
By local$\o$ solvability$\o$ $G\o$ is then solvable, and hence cannot have two distinct 
Borel subgroups. Hence the last condition implies the first one. 
\qed

\bigskip
Lemma \ref{LemEquivCondL0R0NotSolv} can be refined as follows in the 
locally$\o$ solvable case. 

\begin{lem}
Let $G$ be a locally$\o$ solvable group of finite Morley rank. 
Then the following are equivalent. 
\begin{itemize}
\item[$(1)$]
$N\o(A)<G\o$ for each nontrivial definable abelian subgroup $A$ of $G$. 
\item[$(2)$]
$G\o$ is not solvable. 
\item[$(3)$]
$G\o$ has two distinct Borel subgroups. 
\end{itemize}
\end{lem}
\proof
As in Lemma \ref{LemEquivCondL0R0NotSolv}. 
If $N\o(A)=G\o$ for some nontrivial definable abelian subgroup $A$ 
of $G$, then $G\o$ is now solvable by local$\o$ solvability. 
\qed

\bigskip
In $\PSL_{2}$, normalizers$\o$ of unipotent subgroups correspond to Borel subgroups. 
The following is a first approximation of this in locally$\o$ solvable$\o$ 
groups. 

\begin{lem}\label{NUnipHeavyBorel}
Let $G$ be a locally$\o$ solvable$\o$ group of finite Morley rank. Assume that for 
$q$ prime or infinite $d_{q}(G)\geq 1$, 
and let $U$ be a Sylow $(q,d_{q}(G))$-subgroup of $G$. 
Then $N\o(U)$ is a Borel subgroup of $G$. 
\end{lem}
\proof
By local$\o$ solvability$\o$ of $G$, $N\o(U)\leq B$ for some Borel subgroup $B$. 
Now Fact \ref{FactUnipHeaviest} implies $U\leq F\o(B)$, and in particular 
$B\leq N\o(U)$ by maximality of $U$. Hence $N\o(U)=B$ is a Borel subgroup of $G$. 
\qed

\begin{lem}\label{LemU(B)SylOfG}
Let $G$ be a locally$\o$ solvable$\o$ group of finite Morley rank, 
$\tilde p=(p,r)$ a unipotence parameter with $r>0$, and $B$ a Borel subgroup 
of $G$ such that $d_{p}(B)=r$. 
Then $U_{\tilde p}(B)$ is a Sylow $\tilde p$-subgroup of $G$. 
\end{lem}
\proof
By Fact \ref{FactUnipHeaviest}, $U:=U_{\tilde p}(B)$ is in 
$F\o(B)$, and in particular is a $\tilde p$-group. It is obviously 
definably characteristic in $B$. 
If $U<V$ for some Sylow $\tilde p$-subgroup $V$ of $G$, then 
$U<U_{\tilde p}(N_{V}(U))$ by normalizer condition, Fact \ref{NormCond}. 
But as $N\o(U)$ is solvable by local$\o$ solvability$\o$ of $G$, and contains $B$, 
it is $B$ by maximality of $B$. 
Hence $U<U_{\tilde p}(N_{V}(U))\leq U_{\tilde p}(B)=U$, a contradiction. 
\qed

\bigskip
When $r=0$ Lemma \ref{LemU(B)SylOfG} fails. For example, in the 
standard Borel subgroup $B$ of $\PSL_{2}$ over a pure algebraically closed field 
of positive characteristic, $U_{(\infty,0)}(B)=B$. However the lemma becomes 
true for $r=0$ if one assumes that the absolute unipotence degree of $B$ satisfies 
$d(B)=0$. 

\subsection{Semisimple groups}

Obviously with locally solvable groups one becomes quickly interested 
in normal solvable subgroups. 

\begin{fait}\label{FaitExistSolvRad}
{\bf \cite[Theorem 7.3]{BorovikNesin(Book)94}}
Let $G$ be a group of finite Morley rank. Then $G$ has a largest normal solvable 
subgroup, which is definable. It is denoted by $R(G)$ and 
called the {\em solvable radical} of $G$. 
\end{fait}

\begin{defi}
Let $G$ be a group of finite Morley rank. We say that 
\begin{itemize}
\item[$(1)$]
$G$ is {\em semisimple} if $R(G)=1$, or equivalently if $N(A)<G$ for each nontrivial 
abelian subgroup $A$ of $G$. 
\item[$(2)$]
$G$ is {\em semisimple$\o$} if $R\o(G)=1$, or equivalently if $N(A)<G$ for each 
nontrivial connected abelian subgroup $A$ of $G$.
\end{itemize} 
\end{defi}

Of course, if $G$ is any group of finite Morley rank, then $G/R(G)$ is semisimple 
and $G/R\o(G)$ is semisimple$\o$, as solvable-by-solvable groups are solvable. 

\begin{fait}\label{FactGConFinCentG/Z}
{\bf \cite[Lemma 6.1]{BorovikNesin(Book)94}}
Let $G$ be a connected group of finite Morley rank with a finite center. 
Then $G/Z(G)$ has a trivial center. 
\end{fait}

The following fact has certainly been implicit in previous arguments, and we 
just state it precisely. 

\begin{fait}\label{GConFinRadG/RGSemi}
Let $G$ be a connected group of finite Morley rank with $R(G)$ finite. 
Then $R(G)=Z(G)$ and $G/R(G)$ is semisimple. 
\end{fait}
\proof
The connected group $G$ acts by conjugation on its finite solvable radical 
$R(G)$, and thus by Fact \ref{ActionOnFiniteSets} $R(G)\leq Z(G)$. 
As the center is always contained in the solvable radical one gets $R(G)=Z(G)$. 
The semisimplicity of $G/R(G)$ is always true. 
\qed

\begin{lem}\label{LemHNonSolvSemio}
Let $G$ be a group of finite Morley rank and 
$H$ a nonsolvable definable connected subgroup of $G$. 
\begin{itemize}
\item[$(1)$] 
If $G$ is locally$\o$ solvable$\o$, then $H$ is semisimple$\o$, $R(H)=Z(H)$ is finite 
and $H/R(H)$ is semisimple. 
\item[$(2)$]
If $G$ is locally$\o$ solvable, then $H$ is semisimple. 
\end{itemize}
\end{lem}
\proof
This is obvious by definitions and Fact \ref{GConFinRadG/RGSemi}. 
\qed

\subsection{New configurations}\label{SectionSpeedLimits}

All the work concerning minimal connected simple groups of finite Morley rank generalizes 
identically to the case of locally$\o$ solvable groups of finite Morley rank. The 
reason is that in the study of minimal connected simple groups every argument is based 
on the consideration of normalizers$\o$ of nontrivial subgroups $X$. 
If such a subgroup 
$X$ is finite, then its normalizer$\o$ coincides with its centralizer$\o$ by 
Fact \ref{ActionOnFiniteSets}. 

When dealing with the more general class of locally$\o$ solvable$\o$ groups 
centra\-lizers$\o$ of elements of finite order might be nonsolvable. 
In the present papers we try to concentrate exclusively on the more 
general class of locally$\o$ solvable$\o$ groups, and hence new phenomena 
can appear. In the present section we try to give an overview of the new 
pathological configurations which might occur in this context. 
We see these new configurations as some kind of ``speed limits" when 
generalizing arguments from the 
minimal connected simple/locally$\o$ solvable case to the more general 
locally$\o$ solvable$\o$ case. 

Recall from Lemma \ref{FirstRemOnLocSolvGps} that 
$$\{ \mbox{locally$\o$ solvable groups} \}\subsetneq 
\{ \mbox{locally$\o$ solvable$\o$ groups} \},$$
the inclusion being strict. The main (and unique) 
example in the algebraic category of a connected 
group which is locally$\o$ solvable$\o$ but not locally$\o$ solvable is 
$\SL_{2}(K)$, with $K$ an algebraically closed field of characteristic different from $2$: 
its solvable radical consists of a cyclic group of order $2$. 

In the context of groups of finite Morley rank there might be other configurations 
occuring, contradicting even the latter property a priori. 
In what follows we just make a list of potential pathological configurations of 
connected locally$\o$ solvable$\o$ groups of finite Morley rank which are not 
locally$\o$ solvable, and which remain at the end of our classification. 

A {\em full Frobenius} group is a group $G$ with a proper subgroup $H$ such that 
$$H\mbox{ is malnormal in }G\mbox{ and }{G=H^{G}}$$
and we often use the sentence ``$H<G$ is a full Frobenius group" to specify the subgroup 
$H$. The existence of such groups of finite Morley rank is the main obstacle 
to the Algebraicity Conjecture for simple groups of finite Morley rank. We just 
record a few basic properties of such groups, if they exist. 

\begin{fait}\label{FactFullFrobGps}
{\bf \cite[Propositions 3.3 and 3.4]{Jaligot01}}
Let $H<G$ be a full Frobenius group, with $G$ of finite Morley rank and connected. Then 
\begin{itemize}
\item[$(1)$]
$C(x)\leq H$ and is infinite for each nontrivial element $x$ of $H$. 
\item[$(2)$]
$H$ is definable in the pure group $G$ and connected. 
\item[$(3)$]
$HgH\cap Hg^{-1}H=\emptyset$ for any element $g$ in $G\setminus H$. 
\item[$(4)$]
$\rk(G)\geq 2\rk(H)+1$.  
\item[$(5)$]
There exists a nontrivial definable simple subgroup $\tilde{G}$ of $G$ such that 
$(H\cap \tilde{G})<\tilde{G}$ is a full Frobenius group.  
\end{itemize}
\end{fait}

We often call a group $G$ as in Fact \ref{FactFullFrobGps} and with $H$ nilpotent 
a {\em bad} group. (This notion is floppy.) In any case these groups have no involutions, 
and hence their torsion can involve only odd primes. 

We view the following potential configuration of locally$\o$ solvable$\o$ group, 
or any of its natural variations, as a kind of ``universal conterexample" 
to the algebraic case as far as torsion is concerned. Elements belonging to a decent 
torus are called {\em toral}. 

\begin{config}\label{BadConfig1}
$G$ is a connected locally$\o$ solvable$\o$ group of finite Morley rank with a 
proper (definable connected) subgroup $H$ such that 
\begin{itemize}
\item[$(1)$]
$H<G$ is a full Frobenius group. 
\item[$(2)$]
$R(H)=Z(H)$ is finite and nontrivial, consisting of $p$-toral elements of $H$ for some 
prime $p$. 
\item[$(3)$]
$H/Z(H)$ is a full Frobenius group for some proper definable connected solvable 
subgroup $B/Z(H)$.  
\item[$(4)$]
$B/Z(H)$ has nontrivial $p$-unipotent subgroups, for some prime $p$ dividing 
$|Z(H)|$, and also nontrivial $q$-unipotent subgroups for other primes $q$.  
\end{itemize}
\end{config}

A group $G$ as in Configuration \ref{BadConfig1} would have {\em $p$-mixed type}, 
i.e. containing both nontrivial $p$-tori and $p$-unipotent subgroups, and have 
nontrivial $q$-unipotent subgroups for several primes $q$. 

In $\SL_{2}$, a generic element belongs to a maximal torus, and in particular 
to the connected component of its centralizer. Here is another potential new 
pathological phenomenon with locally$\o$ solvable$\o$ groups. 

\begin{config}\label{BadConfig2}
$G$ is a connected locally$\o$ solvable$\o$ group of finite Morley rank with a 
proper (definable connected) subgroup $B$ such that. 
\begin{itemize}
\item[$(1)$]
$B<G$ is a full Frobenius group. 
\item[$(2)$]
$B$ is a nilpotent group such that, for $x$ generic in $B$, $x\notin C\o(x)$. 
\end{itemize}
\end{config}

A generic element $x$ of a group $G$ as in Configuration \ref{BadConfig2} 
would satisfy $x\notin C\o(x)$. We note that examples of connected nilpotent groups 
$B$ of finite Morley rank as in clause $(2)$ of Configuration \ref{BadConfig2} 
are provided by \cite[\S3.2.3]{BorovikNesin(Book)94} or the Baudisch 
$2$-nilpotent group \cite{Baudisch96}. With such subgroups $B$ a group 
$G$ as in Configuration \ref{BadConfig2} would be locally$\o$ solvable, 
but if $G$ had the prescribed property modulo a nontrivial finite center, 
then it would not be locally$\o$ solvable. 

Even with involutions and algebraic subgroups one can imagine the following 
configuration which seems to remain open at the end 
of our second paper \cite{DeloroJaligotII}. 

\begin{config}\label{ConfigSL2Proper}
$G$ is a connected locally$\o$ solvable$\o$ group of finite Morley rank 
with an involution $i$ such that $C(i)<G$ and 
$C(i)\simeq \SL_{2}(K)$ 
for some algebraically closed field $K$ of characteristic different from $2$. 
\end{config}

In \cite{CherlinJaligot2004} all nonalgebraic configurations are known 
to have nongenerous Borel subgroups. Even assuming all Borel subgroups generous 
does not seem to be helpful in \cite{DeloroJaligotII} 
toward finding a contradiction in 
Configuration \ref{ConfigSL2Proper}. This is a major new phenomenon 
possibly occuring in the locally$\o$ solvable$\o$ case as opposed 
to the minimal connected simple/locally$\o$ solvable one. 

\subsection{Local$\o$ solvability/solvability$\o$}\label{SectionL0RvsL0R0}

In Section \ref{SectionSpeedLimits} we saw certain speed limits when 
considering generalizations to the wider class of 
locally$\o$ solvable$\o$ groups, which usually rely on the existence 
of certain semisimple$\o$ but not semisimple groups. 
We nevertheless intend in this section to start dealing with these aspects in the 
general class of locally$\o$ solvable$\o$ groups of finite Morley rank, bearing 
in mind the speed limits of Section \ref{SectionSpeedLimits}. 
For this purpose it is useful to study systematically 
subgroups of the form $C\o(x)$ in locally$\o$ solvable$\o$ groups. 
When such a subgroup is not solvable it has a 
finite solvable radical, which is then the center, and its quotient modulo the center is 
semisimple. This boils down to the study of semisimple locally$\o$ solvable$\o$ 
groups. 

We start with some generalities. 

\begin{lem}\label{LemCentNonSolvFinite}
Let $G$ be a locally$\o$ solvable$\o$ group of finite Morley rank. 
If $H$ is a nonsolvable definable connected subgroup of $G$, then 
$C_{G}(H)$ is finite. 
\end{lem}
\proof
Assume $C\o_{G}(H)$ infinite. Then it contains a nontrivial 
definable connected solvable subgroup $B$ 
by Lemma \ref{LemGenSmallMaxUnipParam}. 
We then have $H\leq C\o(B) \leq N\o(B)$, which must be solvable 
by local$\o$ solvability$\o$ of $G$. 
\qed

\bigskip
In a locally$\o$ solvable$\o$ group $G$ of finite Morley rank, we call a subset 
$X$ {\em exceptional} in $G$ if $C\o(X)$ is nonsolvable. Such sets are finite 
by Lemma \ref{LemCentNonSolvFinite}, and as $C(X)=C(\<X\>)$ any such 
subset $X$ can always be identified with the finite subgroup it generates. 

Dually, we call a definable connected subgroup $H$ {\em exceptional} in $G$ 
if $H$ is nonsolvable. Then $C(H)$ centralizes the nonsolvable definable 
connected subgroup $H$ and is an exceptional subset of $G$. 

We denote by ${\cal E}_{f}$ and ${\cal E}_{s}$ the set of finite exceptional 
subgroups of $G$ and the set of nonsolvable definable connected subgroups respectively 
($\cal E$ stands for ``e"xceptional, $f$ for ``f"inite, and $s$ for ``s"emisimple). 
Both sets are nonempty if and only if $G\o$ is nonsolvable. Of course both 
sets are naturally ordered by inclusion. 

Taking centralizers$\o$ $C\o(\cdot)$ from ${\cal E}_{f}$ to ${\cal E}_{s}$ and 
centralizers $C(\cdot)$ from ${\cal E}_{s}$ to ${\cal E}_{f}$ defines a 
{\em Galois connection} between ${\cal E}_{f}$ and ${\cal E}_{s}$ 
(see \cite{Birkhoff67}). That is, and following a similar exposition 
in \cite{AltBorCher(Book)}, they satisfy the following properties. 

\begin{lem}\
\begin{itemize}
\item[$(1)$]
The mappings $C\o$ and $C$ are order-reversing. 
\item[$(2)$]
If $X\in {\cal E}_{f}$ then $X\leq C(C\o(X))$ and if 
$H\in {\cal E}_{s}$ then $H\leq C\o(C(H))$. 
\end{itemize}
\end{lem}

As in any Galois connection, this has the following consequence. 

\begin{prop}
Let $X\in {\cal E}_{f}$ and $H\in {\cal E}_{s}$. Then 
$C\o(X)=C\o(C(C\o(X)))$ and $C(H)=C(C\o(C(H)))$. 
\end{prop}

If we denote for $X$ in ${\cal E}_{f}$ and $H$ in ${\cal E}_{s}$ 
$$\overline{X}=C(C\o(X))\mbox{~and~}\overline{H}=C\o(C(H)),$$
then the two operations $\overline{\phantom{H}}$ are {\em closure operations} on 
${\cal E}_{f}$ and ${\cal E}_{s}$ respectively. That is, they satisfy the following. 

\begin{cor}\
\begin{itemize}
\item[$(1)$]
For $X\in {\cal E}_{f}$ and $H\in {\cal E}_{f}$, we have 
$X\leq \overline{X}=\overline{\overline{X}}$ and 
$H\leq \overline{H}=\overline{\overline{H}}$.
\item[$(2)$]
Monotonicity: For $X_{1}\subseteq X_{2}$ in ${\cal E}_{f}$ and 
$H_{1}\leq H_{2}$ in ${\cal E}_{s}$, we have 
$\overline{X_{1}}\leq \overline{X_{2}}$ and 
$\overline{H_{1}}\leq \overline{H_{2}}$. 
\end{itemize}
The closed elements of ${\cal E}_{f}$ and ${\cal E}_{s}$ are those of the 
form $\overline{X}$ and $\overline{H}$ respectively. 
\end{cor}

One can also refine Lemma \ref{LemCentNonSolvFinite} by giving a uniform 
bound on cardinals of elements of ${\cal E}_{f}$. We first note the following 
general fact. 

\begin{lem}\label{LemmaUnifBoundCCX}
Let $G$ be a group of finite Morley rank. Then there exists a natural number $m$ 
such that, for any subset $X$ of $G$, $|C(X)|\leq m$ or $C(X)$ 
is infinite. 
\end{lem}
\proof
As $G$ is stable, it satisfies the Baldwin-Saxl chain condition 
\cite[\S1.3]{Poizat(book)87}. 
This means that there exists a fixed integer $k$ such 
that, for every subset $X$ of $G$, 
$C(X)=C(x_{1},\cdots , x_{k})$ for some elements $x_{1}$, ..., $x_{k}$ 
of $X$, and the family of all subgroups of the form $C(X)$ is uniformly definable 
(by a formula without parameters). 

Now the uniform bound $m$ on the cardinals of the finite sets 
of the family is provided by elimination of infinite quantifiers 
\cite[Proposition 2.2]{MR1955252}. 
\qed

\begin{lem}\label{LemUnifBoundException}
Let $G$ be a locally$\o$ solvable$\o$ group of finite Morley rank. Then there exists a 
natural number $m$ bounding uniformly the cardinals of finite exceptional subsets of $G$. 
\end{lem}
\proof
Let $m$ be as in Lemma \ref{LemmaUnifBoundCCX}. 
If $X$ is a finite exceptional subset of $G$, then $C\o(X)$ is 
nonsolvable and $C(C\o(X))$ is finite by Lemma \ref{LemCentNonSolvFinite}. 
As $C(C\o(X))$ is a finite centralizer, its cardinal is uniformly bounded by $m$. 
Now $X\subseteq C(C\o(X))$, and thus the cardinal of 
$X$ is uniformly bounded by $m$. 
\qed

\bigskip
If $G$ is a locally$\o$ solvable$\o$ group of finite Morley rank, we call 
{\em exception index} and denote by $e(G)$ the maximal integer $m$ such 
that $G$ has an exceptional nonsolvable definable connected subgroup centralizing 
a subset $X$ with $m$ elements. Notice that $X$ coincides with $\< X \>$, so 
that $e(G)$ is the largest cardinal of an exceptional subgroup in ${\cal E}_{f}$. 

Maximal exceptional subgroups of ${\cal E}_{f}$ correspond to 
minimal exceptional subgroups of ${\cal E}_{s}$, and vice-versa. 
A case of particular interest is the following. 

\begin{lem}
Minimal nontrivial exceptional subgroups are cyclic of prime order. 
\end{lem}
\proof
Obvious. 
\qed

\bigskip
One can clarify the structure of elements of ${\cal E}_{s}$ as follows. 

\begin{lem}\label{LemHExceptStrutGen}
Let $G$ be a locally$\o$ solvable$\o$ group of finite Morley rank and 
$H$ an exceptional nonsolvable definable connected subgroup. 
Then $H$ is semisimple$\o$, $R(H)=Z(H)$ is finite and $H/R(H)$ is semisimple. 
\end{lem}
\proof
Lemma \ref{LemHNonSolvSemio} $(1)$. 
\qed

\bigskip
Of course, the sets of closed sets in ${\cal E}_{f}$ and ${\cal E}_{s}$ are at most 
reduced to $\{1\}$ and $\{G\o\}$ in the locally$\o$ solvable case 
(in case $G\o$ is nonsolvable, and empty otherwise). 

The next lemma seems to be the only way to get locally$\o$ solvable 
groups out of locally$\o$ solvable$\o$ ones. 

\begin{lem}\label{HExceptHRHSemisimple}
Let $G$ be a locally$\o$ solvable$\o$ group of finite Morley rank and 
$H$ a nonsolvable definable connected subgroup exceptional in $G$, 
which is minimal with respect to this property. Then $H/R(H)$ is locally$\o$ solvable. 
\end{lem}
\proof
The conclusions of Lemma \ref{LemHExceptStrutGen} are valid in $H$ 
and we will use them freely. Denote by $\overline{\phantom{H}}$ the quotients 
by $R(H)$, and let $A$ be the preimage in $H$ of a nontrivial definable 
abelian subgroup $\overline{A}$ of $\overline{H}$. Of course 
$A$ is a definable solvable subgroup of $H$. 

Let $N$ be the preimage of $N\o_{\overline{H}}(\overline{A})$ in $H$. 
We have $N\leq N(A)$ and $\overline{N}\o=\overline{N\o}$, so that 
$N=N\o R(H)$. 

If $A$ is infinite modulo $R(H)$, then $A$ is infinite as well, 
and as $N\o\leq N\o(A)$ we get $N\o$ solvable by local$\o$ solvability$\o$ of $G$ 
and Lemma \ref{FirstRemOnLocSolvGps} $(1)$. Then $N=N\o R(H)$ is solvable, 
as well as $N\o_{\overline{H}}(\overline{A})$. 

If $A$ is finite modulo $R(H)$, then $A$ is finite as well as $R(H)$ is. Now 
$N\o$ acts on the nontrivial finite group $A$, and therefore centralizes it 
by Fact \ref{ActionOnFiniteSets}. The minimality of $H$ yields 
$N\o$ solvable or $N\o=H$. In the first case one concludes that 
$N=N\o R(H)$ is solvable, as well as $N\o_{\overline{H}}(\overline{A})$. 
The second case implies that $A\leq Z(H)=R(H)$, and is thus impossible as 
$A$ is nontrivial modulo $R(H)$. 
\qed

\bigskip
Lemma \ref{HExceptHRHSemisimple} seems to be a very rough indication 
that the new locally$\o$ solvable$\o$ groups which are not locally$\o$ solvable 
are more or less as in Configuration \ref{BadConfig1}. 

We also note that exceptional nonsolvable definable connected subgroups attached 
to a nontrivial finite exceptional subgroup are of finite index in their normalizers. 

\begin{lem}
Let $G$ be a locally$\o$ solvable$\o$ group of finite Morley rank. 
If $X$ is an exceptional finite subset of $G$, then $N\o(C\o(X))=C\o(X)$.
\end{lem}
\proof
$N\o(C\o(X))$ normalizes $C(C\o(X))$, which is finite and contains $X$. So it 
centralizes $X$ by Fact \ref{ActionOnFiniteSets}, and we are done. 
\qed

\bigskip
A natural question is to know whether exceptional finite subsets $X$ are contained 
in their attached exceptional nonsolvable definable connected subgroups, i.e. 
whether $X\subseteq C\o(X)$. This would follow from the more general, 
but similarly natural, question to know whether nonsolvable definable 
connected subgroups are selfnormalizing. This is the kind of problem 
which seems optimistically trackable when $C\o(X)$ is generous in the 
ambient group, since the intensive experience on Weyl groups from 
\cite{CherlinJaligot2004}, and we will get positive answers in the most 
interesting situations in Section \ref{SectionLocAnalandGen} below. 

We are now going to look more closely at the interesting case in which an 
exceptional finite subgroup $X$ of ${\cal E}_{f}$ satisfies $X\leq C\o(X)$. 
In this case $X\leq Z(C\o(X))$, and $X$ is in particular an abelian finite 
subgroup. Typical finite abelian groups belonging to the connected component 
of their centralizers are the finite subgroups of decent tori. 
(And this is in general not true around groups of bounded exponent, as noticed 
after Configuration \ref{BadConfig2}.)

\begin{lem}\label{LemExcepEltsOfDecentTori}
Let $G$ be a locally$\o$ solvable$\o$ group of finite Morley rank and $T$ a 
maximal definable decent torus of $G$. Then the union of elements of 
${\cal E}_{f}$ contained in $T$ is finite and invariant under any 
automorphism of $G$ leaving $T$ invariant. 
\end{lem}
\proof
For the finiteness we can use Lemma \ref{LemUnifBoundException} to get a uniform 
bound, at most the exception index $e(G)$ of $G$, on the cardinals of the finite groups 
involved. Then, as Pr\"ufer $p$-ranks are finite for any prime $p$ in a decent torus, 
subgroups of order at most $e(G)$ must be contained in a finite subgroup of $T$.  

The second point is obvious. 
\qed

\bigskip
A question, which might be difficult, is to know whether the union in 
Lemma \ref{LemExcepEltsOfDecentTori} is necessarily a (finite) 
subgroup of $T$, and is itself exceptional. If this were the case, then 
calling this group $E$, one would have a nonsolvable group 
$C\o(E)/R(C\o(E))$ where nontrivial toral elements are not exceptional anymore. 
This is a desirable property for certain questions such as bounding Pr\"ufer ranks, 
as we will see later, in our treatment of odd type groups \cite{DeloroJaligotII}. 
This desirable property can however be obtained as follows. 

\begin{lem}\label{LemKillExceptionsInTori}
Let $G$ be a connected nonsolvable locally$\o$ solvable$\o$ group of finite Morley 
rank, and $T$ a maximal definable decent torus of $G$. 
Then $G$ has an exceptional nonsolvable definable connected subgroup $H$ 
containing $T$ and such that $C\o(\overline{t})$ is solvable for 
any nontrivial toral element $\overline{t}$ of $H/R(H)$. 
\end{lem}
\proof
Let $X$ be a maximal exceptional finite subgroup of $T$. Then $H=C\o(X)$ 
is nonsolvable. As $X\leq T$ and $T$ is abelian and connected, 
$X\leq T\leq H$. 

Let now $\overline{t}$ be a nontrivial toral element of 
$\overline{H}=H/R(H)$. By pullback of decent tori, 
Fact \ref{quotindec2} $(2)$ or rather \cite[Lemma 3.1]{frecontame}, and 
Fact \ref{ConjdecentTori}, we may assume $t$ in $TR(H)$, i.e. 
$t=t'r$ for some $t'\in T$ and some $r\in R(H)$. As in 
Lemma \ref{HExceptHRHSemisimple}, one sees that the preimage of 
the centralizer$\o$ of $t$ modulo $R(H)$ cannot be nonsolvable: otherwise its 
connected component would centralize $t=t'r$, and as $r\in R(H)=Z(H)$ it would 
centralize $t'$, so that $X\<t'\>$ would be an exceptional finite subgroup of $T$ 
containing $X$ properly, a contradiction. This finishes our proof. 
\qed

\bigskip
Before moving ahead we close the present section by 
describing more precisely the set of exceptional subsets of 
a decent torus $T$ as in Lemma \ref{LemKillExceptionsInTori}, or more generally 
of an {\em arbitrary} subset $T$ of a locally$\o$ solvable$\o$ group $G$ 
of finite Morley rank. 

First we naturally consider the notion of closure relative to $T$. 
For $X$ an exceptional subset of $T$, we say that $X$ is closed {\em in} $T$ if 
$X = C_T(C\o (X))$. Of course the notion of relative closedness is robust. 

\begin{rem}\label{RemRelClosure}
Any set of the form $C_{T}(C\o(X))$ is closed in $T$. 
\end{rem}
\proof 
As $X\subseteq C_{T}(C\o(X))\subseteq \overline{X}$, 
$\overline{X}=C(C\o(C_{T}(C\o(X))))$ by taking the closure in $G$, and 
$C_{T}(C\o(C_{T}(C\o(X))))=C_{T}(C\o(X))$ by taking the intersection with $T$. 
\qed

\medskip
The poset of exceptional subsets of $T$ is best described as follows by the notion 
of minimal extensions of closed subsets. We say that $(X_{1},X_{2})$ is 
a {\em minimal} extension of closed sets of $T$ if $X_{1}\subsetneq X_{2}$ are 
two exceptional subsets of $T$ closed in $T$ and any closed subset $Y$ of $T$ such that 
$X_{1}\subseteq Y\subseteq X_{2}$ is either $X_{1}$ or $X_{2}$. The 
relation ``$(X_{1},X_{2})$ is a minimal extension of closed sets of $T$" 
defines an oriented graph on the set of closed sets of $T$, which is clearly irreflexive, 
antisymmetric, and loop-free, that is without cycles preserving the orientation 
(but possibly with cycles not preserving the orientation). 
We call this graph the {\em graph of exceptional subsets of $T$}. Its main properties are the 
following. 

\begin{lem}\label{extensionsminimales}
Let $G$ be a locally$\o$ solvable$\o$ group of finite Morley rank and $T$ an arbitrary 
subset of $G$. 
\begin{itemize}
\item[$(1)$]
Assume $(X_{1},X_{2})$ is a minimal extension in the graph of exceptional subsets of 
$T$ and $Y$ is a subset such that $X_{1}\subsetneq Y \subseteq X_2$. 
Then $C_T (C\o (Y)) = X_2$. Moreover $C\o(X_{2})<C\o(X_{1})$. 
\item[$(2)$]
Assume $(X,X_{1})$ and $(X,X_{2})$ are two minimal extensions in the graph 
of exceptional subsets of $T$. Then either $X_{1}=X_{2}$ or $X_{1}\cap X_{2}=X$. 
\end{itemize}
\end{lem}

\proof
$(1)$. 
$X_{1}\subsetneq Y\subseteq C_{T}(C\o(Y))\subseteq X_{2}$ and as 
$C_{T}(C\o(Y))$ is closed in $T$ by Remark \ref{RemRelClosure} 
it must be $X_{2}$ by minimality of the extension 
$(X_{1},X_{2})$. 

The claim that $C\o(X_{2})<C\o(X_{1})$ follows merely from 
the fact that $X_{1}\neq X_{2}$ are closed in $T$. 

$(2)$. 
Let $Y=X_1 \cap X_2$. If $X\subsetneq Y$, then the first point implies that 
$X_{1}=C_T(C\o (Y))=X_2$. 
\qed

\bigskip
Finally, we note that the graph of exceptional subsets of $T$ as in 
Lemma \ref{extensionsminimales} always has a ``minimal" element, 
namely $T\cap Z(G)$, and ``maximal" elements, corresponding to the maximal 
traces on $T$ of exceptional sets in ${\cal E}_{f}$, which are of cardinal at most 
$e(G)$. We also note that the graph has a {\em finite height}: the length of a maximal chain 
of exceptional closed sets in $T$ is at most $e(G)$. 

When $T$ is a nilpotent divisible subgroup of $G$ (for example as in 
Lemmas \ref{LemExcepEltsOfDecentTori} and \ref{LemKillExceptionsInTori}), 
then exceptional subsets of $T$ are necessarily in a same decent torus 
(the maximal decent torus of the definable hull of $T$) and by 
Lemma \ref{LemExcepEltsOfDecentTori} applied in this decent torus the 
graph of exceptional subsets of $T$ is finite. 

\subsection{Genericity}

\begin{fait}\label{DecentTorusInGenCarter}
{\bf (Compare with \cite[Theorem 7.3]{FreconJaligot07})}
Let $G$ be a locally$\o$ solv\-able$\o$ group of finite Morley rank with 
a nontrivial decent torus $T$, and $Q$ a Carter subgroup of $G$ containing $T$. 
Then $Q$ is generous in $G\o$, and 
$T\leq \tilde{T} \leq Q$ for some maximal definable decent torus 
$\tilde{T}$ of $G$. 
\end{fait}
\proof
The existence of $Q$ is guaranteed by Fact \ref{FactExistCarter}, 
as decent tori are of minimal unipotence degree. 

By Fact \ref{ConjdecentTori}, $C\o(T)$ is generous in $G\o$. 
Now $C\o(T)$ is solvable by local$\o$ solvability$\o$ of $G$, and the Carter 
subgroup $Q$ is generous in $C\o(T)$ by Fact \ref{Carterconjugateandselfnorm}. 
It follows that $Q$ is generous in $G\o$ by the transitivity of generosity 
provided in \cite[Lemma 3.9]{Jaligot06}. 

Doing the same argument as above for a maximal definable decent torus 
$\tilde T$ containing $T$, one gets a generous Carter subgroup $\tilde Q$ of $G\o$ 
containing $\tilde T$, and as generous Carter subgroups are conjugate by 
Fact \ref{FactConjGenCarter} one gets that $Q$ contains a maximal definable decent 
torus, which necessarily contains $T$. 
\qed

\bigskip
We record here an application of Lemma \ref{LemGenXcapMGenInXGInv} 
in the case of locally solvable groups of finite Morley rank. This will be the 
{\em cl\'e de vo\^{u}te} for a concentration argument in one of the most 
prominent theorem on odd type groups in \cite{DeloroJaligotII}. 

\begin{lem}\label{Concentrationimpliessolvability}
{\bf (Compare with \cite[Corollaire 2.4]{Deloro05})}
Let $G$ be a group of finite Morley rank and $X$ a nonempty 
definable $G\o$-invariant subset of $G\o$. Let $M$ be a definable solvable subgroup 
of $G\o$ such that $X\cap M$ is generic in $X$. 
\begin{itemize}
\item[$(1)$]
If $G$ is locally$\o$ solvable and $X\neq {\{1\}}$, then $G\o$ is solvable. 
\item[$(2)$]
If $G$ is locally$\o$ solvable$\o$ and $X$ is infinite, then $G\o$ is solvable. 
\end{itemize}
\end{lem}
\proof
$(1)$. Let $Y$ be the definable $G\o$-invariant subset of $X\cap M$ 
generic in $X$ provided by Lemma \ref{LemGenXcapMGenInXGInv}. 
As $X$ is nonempty, $Y$ is also nonempty, and $G\o=N\o(\<Y\>)$. Now 
$\<Y\>$ is a subgroup of $M$, and hence is solvable. If it is nontrivial, then 
$G\o$ must be solvable by Lemma \ref{LemGenCharLocSolvGps} $(3)$. 
Otherwise, $\{1\}$ is a generic subset of $X$, and $X$ must be finite. Hence 
$X$ is a finite set of finite conjugacy classes, with one nontrivial by 
assumption. This nontrivial finite $G\o$-conjugacy class must be central in $G\o$ 
by Fact \ref{ActionOnFiniteSets}, and as $G\o$ has then a nontrivial center 
it must again be solvable by local$\o$ solvability. 

$(2)$. One argues in the same way. Now, as $X$ is infinite, $Y$ is also infinite 
by genericity. As $G\o=N\o(\<Y\>)$ is $\<Y\>$-local$\o$ with $\<Y\>$ infinite 
and solvable, as contained in $M$, 
Lemma \ref{LemGenCharLocSolvGps} $(4)$ now gives the solvability of $G\o$. 
\qed

\section{Local analysis}\label{SectLocAnal}

We now proceed to the {\em local analysis} of locally solvable groups of finite 
Morley rank, that is the analysis of intersections of their (most interesting) subgroups. 

In Section \ref{SectionJaligLemma} we deal with a series of results 
which correspond to the Bender method in finite group theory. 
In general these lemmas say in our context that sufficiently unipotent subgroups 
of locally solvable 
groups of finite Morley rank are disjoint, like unipotent subgroups in $\PSL_{2}$ or 
$\SL_{2}$. They are the main tool for analyzing locally solvable groups, notably 
the only trick involving unipotence in the recognition of $\PSL_{2}$ 
in the algebraic parts of our second paper \cite{DeloroJaligotII}. 
The original form was first proved in the context of minimal connected 
simple groups in the unpublished \cite{jal3prep}. It was in a form embryonal 
in characteristic $0$ compared to the one provided later by the general abstract 
unipotence theory of Burdges as in Section \ref{SectionAbstractUnipTh}, 
but both in positive and null characteristic. 
Then they appeared in the tame context in \cite[Section 3.4]{CherlinJaligot2004} where 
they were treated essentially as the positive characteristic case, i.e. involving no 
particular graduation in the unipotence theory. 
The positive characteristic case was recalled as the outline of 
\cite{BurdgesJalLemma}, and later the characteristic $0$ case appeared in 
\cite[\S 3.2]{Deloro05} for the recognition of $\PSL_{2}$ in characteristic 
different from $2$. We are going to give forms of these lemmas entirely 
uniform in the unipotence degrees, in particular independent of the characteristics, 
and in the most general context of locally$\o$ solvable$\o$ groups. 

Section \ref{SectMaxPaire} will then concern the situation in which 
a particular consequence of such uniqueness theorems fails. 
This is a priori a possibility when the subgroups considered 
are not unipotent enough relative to the ambient group. 
The pathological situation appearing can be analyzed somehow 
by replacing the maximality in terms of unipotence degrees by a maximality 
for inclusion concerning a pair of Borel subgroups involved. The endless, but 
very precise, description resulting is the bulk of \cite{BurdgesJalLemma}, 
and in the context of locally$\o$ solvable$\o$ groups we will 
follow the exposition of this paper {\em verbatim}. This full description 
will be applied one time in a nonalgebraic situation in \cite{DeloroJaligotII} and 
that's why we need to restate, slightly more generally but in its full detail, 
this analysis from \cite{BurdgesJalLemma}. 

\subsection{Uniqueness Theorem}\label{SectionJaligLemma}

\subsubsection{The main theorem}

The following Uniqueness Theorem is our analog of the Bender method in groups 
of finite Morley rank and is the main tool for analyzing locally$\o$ solvable$\o$ groups of 
finite Morley rank. There are various forms of this theorem but the present one 
seems to be the most relevant, at least for our applications in \cite{DeloroJaligotII}. 
Its consequences on Borel subgroups in Section \ref{SectConseqBorelSngps} below 
will be the closest analogs of the Uniqueness Theorem 
of Bender in finite group theory \cite{MR0288180} \cite{MR0262351} 
\cite[\S5-7]{MR0407127} \cite[Chapter II]{benderglauberman1994}

\begin{theo}\label{UniquenessLemma}
Let $G$ be a locally$\o$ solvable$\o$ group of finite Morley rank, 
$\tilde p=(p,r)$ a unipotence parameter with $r>0$, 
and $U$ a Sylow $\tilde p$-subgroup of $G$. 
Assume that $U_{1}$ is a nontrivial definable $\tilde p$-subgroup of $U$ containing 
a nonempty (possibly trivial) subset $X$ of $G$ such that $d_{p}(C\o(X))\leq r$. 
Then $U$ is the unique Sylow $\tilde p$-subgroup of $G$ 
containing $U_{1}$, and in particular $N(U_{1})\leq N(U)$. 
\end{theo}

Before the proof, a few remarks.
\begin{itemize}
\item[$(1)$]
If $G\o$ is solvable and $r=d_{p}(G)>0$, then assuming that 
$F\o(G\o)$ has a nontrivial normal definable $\tilde{p}$-subgroup $U_{1}$ 
one gets with Theorem \ref{UniquenessLemma} (applied with $X=1$ for example!) 
that $G\o$ has a unique Sylow $\tilde{p}$-subgroup, which is thus normal and 
contained in $F\o(G\o)$. In the event of the absence of such a subgroup $U_{1}$ one 
easily gets the same conclusion with Fact \ref{ActionOnLessUnip} $(3)$ and $(4)$. 
Hence in some sense Theorem \ref{UniquenessLemma} can be seen as a generalization 
from solvable groups to locally$\o$ solvable$\o$ groups of Fact \ref{FactUnipHeaviest}. 
\item[$(2)$]
The nontriviality of $U_{1}$ is needed in Theorem \ref{UniquenessLemma}, 
as in a hypothetic $\tilde p$-homogeneous semisimple 
bad group the trivial subgroup would be contained in infinitely many conjugates of the 
Sylow $\tilde p$-subgroup. 
\item[(3)]
Theorem \ref{UniquenessLemma} fails if $\tilde{p}=(\infty,0)$. 
For exemple if $G$ is of the form $T \times (U\rtimes T)$, with 
$U$ $p$-unipotent for some prime $p$ and $T$ a good torus, whose second 
copy acts faithfully on $U$, then $d_{\infty}(G)=0$, so that all assumptions of 
Theorem \ref{UniquenessLemma} are satisfied with $U_{1}$ the central copy 
of $T$, but the latter is contained in infinitely many conjugates of the maximal 
good torus $T\times T$. We will give in Lemma \ref{UniquenessLemmar=0} below 
a version of Theorem \ref{UniquenessLemma} specific for the unipotence parameter 
$\tilde{p}=(\infty,0)$, by replacing $d_{p}(C\o(X))$ by the absolute unipotence 
degree $d(C\o(X))$ but with no more local solvability assumption. 
\end{itemize}

After these comments we pass to the proof of Theorem \ref{UniquenessLemma}. 

\bigskip
\proof
Assume $V$ is a Sylow $\tilde p$-subgroup of $G$ distinct from $U$ and 
containing $U_{1}$, and chosen so as to maximize the rank of 
$U_{\tilde p}(U\cap V)$. Let $T$ denote $U\cap V$. 
As $U_{1}\leq T$, 
the subgroup $T$ is infinite. As $T$ is nilpotent, 
$N:=N\o(T)$ is solvable by local$\o$ solvability$\o$ of $G$ 
and Lemma \ref{LemGenCharLocSolvGps} $(4)$. 
Notice that $T<U$, as otherwise $U=(U\cap V)\leq V$ and 
$U=V$ by maximality of $U$. Similarly $T<V$, as otherwise 
$V=(U\cap V)\leq U$ and $V=U$ 
by maximality of $V$. In particular by normalizer condition, Fact \ref{NormCond}, 
$U_{\tilde{p}}(T)<U_{\tilde p}(N_{U}(T))$ and 
$U_{\tilde{p}}(T)<U_{\tilde p}(N_{V}(T))$.

We claim that $d_{p}(N)=r$. 
If $d_{p}(N)>r$, then $r<\infty$, $p=\infty$, and 
$N$ contains a nontrivial Sylow $(\infty,r')$-subgroup $\Sigma$ with $r'>r$. 
Notice that $d_{p}(T)\leq r$ by Corollary \ref{Centrestildep} 
and Fact \ref{StrucNilpGroups} and our assumption that the subset $X$ of 
$T$ satisfies $d_{p}(C(X))\leq r$. 
Then $T\cdot \Sigma$ is nilpotent by 
Fact \ref{ActionOnLessUnip} $(3)$ and $(4)$, and 
$T$ commutes with $\Sigma$ by Fact \ref{StrucNilpGroups}. In particular 
$\Sigma$ commutes with $X$ and $d_{p}(C\o(X))\geq r'>r$, a contradiction to our 
assumption. 
Hence $d_{p}(N)\leq r$, and as $N$ contains $U_{\tilde p}(N_{U}(T))$ 
(or $U_{\tilde p}(N_{V}(T))$) which is nontrivial and 
of unipotence degree $r$ we get $d_{p}(N)=r$. 

By Fact \ref{FactUnipHeaviest} and the assumption that $r\geq 1$ 
we get $U_{\tilde p}(N)\leq F\o(N)$. In particular $U_{\tilde p}(N)$ is 
nilpotent, and contained in a Sylow $\tilde p$-subgroup $\Gamma$ of $G$. Now 
$U_{1}\leq U_{\tilde{p}}(T)<U_{\tilde p}(N_{U}(T))\leq \Gamma$, 
so our maximality assumption on $V$ implies that $\Gamma=U$. In particular 
$U_{\tilde p}(N_{V}(T))\leq \Gamma =U$. But then 
$U_{\tilde p}(T)<U_{\tilde p}(N_{V}(T))\leq 
U_{\tilde p}(U\cap V)=U_{\tilde{p}}(T)$, 
a contradiction which finishes the proof of our first statement. 

The inclusion $N(U_{1})\leq N(U)$ follows from the uniqueness. 
\qed

\bigskip
We conclude with a version of Theorem \ref{UniquenessLemma} specific for the 
unipotence parameter $\tilde p=(\infty,0)$, and which indeed does not rely on 
local solvability. 

\begin{lem}\label{UniquenessLemmar=0}
Let $G$ be a group of finite Morley rank, $T$ a maximal definable 
decent torus of $G$, and $x$ an element of $T$ such that $C\o(x)$ is solvable 
and $d(C\o(x))\leq 0$. Then $T$ is the unique maximal definable decent torus 
of $G$ containing $x$, and in particular $N(\<x\>)\leq N(T)$. 
\end{lem}
\proof
By assumption and Lemma \ref{LemGenSmallMaxUnipParam} $(2)$, $C\o(x)$ is a 
good torus. As $x\in T$ and $T$ is connected abelian, $T\leq C\o(x)$, and 
$T=C\o(x)$ by maximality of $T$. Now any maximal 
definable decent torus containing $x$ must be in $C\o(x)$ for the same reason, 
hence in $T$, and hence equal to $T$ by maximality of $T$. Furthermore, 
$N(\<x\>)\leq N(C\o(x))=N(T)$. 
\qed

\subsubsection{Consequences on Borel subgroups}\label{SectConseqBorelSngps}

Applied to the case of Borel subgroups Theorem \ref{UniquenessLemma} 
has the following corollaries.  These can be seen as absolute 
approximations in the context of locally$\o$ solvable$\o$ groups of 
finite Morley rank of the fact that any unipotent subgroup of 
$\PSL_{2}$ belongs to a unique Borel subgroup of the ambient group.

\begin{cor}\label{corUnicityBorelMax}
Let $G$ be a locally$\o$ solvable$\o$ group of finite Morley rank, 
$\tilde p=(p,r)$ a unipotence parameter with $r>0$, and 
$B$ a Borel subgroup of $G$ such that $d_{p}(B)=r$. 
Let $U_{1}$ be a nontrivial definable $\tilde p$-subgroup of $U_{\tilde p}(B)$ 
containing a nonempty subset $X$ such that $d_{p}(C\o(X))\leq r$. Then 
$U_{\tilde p}(B)$ is the unique Sylow $\tilde p$-subgroup of $G$ 
containing $U_{1}$, and in particular $N(U_{1})\leq N(U_{\tilde p}(B))=N(B)$. 
Furthermore, $B$ is the unique Borel subgroup containing $U_{1}$ and 
admitting $\tilde p$ as a unipotence parameter maximal in its characteristic. 
\end{cor}
\proof
The fact that $U_{\tilde p}(B)$ is a Sylow $\tilde p$-subgroup of $G$ is 
Lemma \ref{LemU(B)SylOfG}. The uniqueness of $U_{\tilde p}(B)$ among 
Sylow $\tilde p$-subgroups containing $U_{1}$, as well as the inclusion 
$N(U_{1})\leq N(U_{\tilde p}(B))$, is then Theorem \ref{UniquenessLemma}.  

Let now $B_{1}$ be a Borel subgroup of $G$ containing $U_{1}$ and admitting 
$\tilde p=(p,r)$ as maximal in its characteristic. Notice that 
$U_{\tilde p}(B_{1})$ is a Sylow $\tilde p$-subgroup of $G$ by 
Lemma \ref{LemU(B)SylOfG}. As it contains $U_{1}$, 
Theorem \ref{UniquenessLemma} now implies 
$U_{\tilde p}(B_{1})=U_{\tilde p}(B)$. Now the normalizers$\o$ of these (equal) 
groups are solvable by local$\o$ solvability$\o$ of $G$, contain $B_{1}$ and $B$ 
respectively, hence are equal to $B_{1}$ and $B$ respectively by maximality, and 
are equal. Hence $B_{1}=B$, as desired for our second claim. 
\qed

\bigskip
$\PSL_{2}$ in positive characteristic offers a counterexample to 
Corollary \ref{corUnicityBorelMax} when $r=0$. It 
suffices to consider for $U_{1}$ a maximal torus of the standard Borel subgroup 
$B$, so that $N(T)\nleq N(B)$ and $T\leq B^{w}$ where $w$ is a nontrivial Weyl group 
element associated to $T$. For the case $r=0$ we refer to 
Lemma \ref{UniquenessLemmar=0}. 

Corollary \ref{corUnicityBorelMax} takes the following form when $(p,r)$ is maximal 
in its characteristic over the whole ambient group $G$. 

\begin{cor}\label{corUnicityBorelSuperMax}
Let $G$ be a locally$\o$ solvable$\o$ group of finite Morley rank, 
$\tilde p=(p,r)$ a unipotence parameter with $r>0$ such that $d_{p}(G)=r$. 
Let $B$ be a Borel subgroup of $G$ such that $d_{p}(B)=r$. Then 
$U_{\tilde p}(B)$ is a Sylow $\tilde p$-subgroup of $G$, 
and if $U_{1}$ is a nontrivial definable $\tilde p$-subgroup of $B$, then 
$U_{\tilde p}(B)$ is the unique Sylow $\tilde p$-subgroup of $G$ 
containing $U_{1}$, $N(U_{1})\leq N(U_{\tilde p}(B))=N(B)$, and $B$ is the unique 
Borel subgroup of $G$ containing $U_{1}$. 
\end{cor}
\proof
Let $X=1$. Then $d_{p}(C\o(X))=d_{p}(G)=r$, 
so Corollary \ref{corUnicityBorelMax} applies with $X=1$. Notice that 
$\tilde p$ is maximal in its characteristic for any Borel subgroup admitting it, 
and that when $U_{1}$ is a nontrivial definable $\tilde p$-group then any Borel 
subgroup containing it admits $\tilde p$. 
\qed

\bigskip
As for Corollary \ref{corUnicityBorelMax}, $\PSL_{2}$ in positive characteristic 
offers a counterexample when $r=0$ in Corollary \ref{corUnicityBorelSuperMax}. 

The preceding uniqueness theorems are often used as follows to ``fusion" Borel subgroups 
sharing too unipotent subgroups. 

\begin{lem}\label{LemFusionBorel1}
Let $G$ be a locally$\o$ solvable$\o$ group of finite Morley rank. Assume that, 
for $i=1$ and $2$, $\tilde {p_{i}}=(p_{i}, r_{i})$ are two unipotence parameters 
with $r_{i} >0$ and $B_{i}$ are two distinct Borel subgroups of $G$ such that 
$d_{p_{i}}(B_{i})=r_{i}$. Then there is no Borel subgroup $B_{3}$ of $G$ 
such that $d_{p_{i}}(B_{i}\cap B_{3})=d_{p_{i}}(B_{3})=r_{i}$ and 
$d_{p_{i}}(C\o(U_{\tilde p_{i}}(B_{i}\cap B_{3})))\leq r_{i}$. 
\end{lem}
\proof
Assume the contrary. Applying Corollary \ref{corUnicityBorelMax} 
with $U_{1}=X=U_{\tilde p_{i}}(B_{1}\cap B_{3})$ implies that $B_{1}=B_{3}$, 
and with $U_{1}=X=U_{\tilde p_{i}}(B_{2}\cap B_{3})$ that $B_{2}=B_{3}$. 
Hence $B_{1}=B_{2}$, a contradiction. 
\qed

\bigskip
We finish with a version of Lemma \ref{LemFusionBorel1} concerning the case 
in which the unipotence degrees $r_{i}$'s are maximized 
over the whole ambient group. 

\begin{lem}\label{LemFusionBorel2}
Let $G$ be a locally$\o$ solvable$\o$ group of finite Morley rank. Assume that, 
for $i=1$ and $2$, $\tilde p_{i}=(p_{i}, r_{i})$ are two unipotence parameters 
with $r_{i} >0$ and $B_{i}$ are two distinct Borel subgroups of $G$ such that 
$d_{p_{i}}(G)=d_{p_{i}}(B_{i})=r_{i}$. Then there is no Borel subgroup 
$B_{3}$ of $G$ such that $d_{p_{i}}(B_{i}\cap B_{3})=r_{i}$. 
\end{lem}
\proof
Under the stated assumptions $d_{p_{i}}(B_{i})=r_{i}$. 
If there was a contradicting Borel subgroup $B_{3}$, then 
$d_{p_{i}}(B_{i}\cap B_{3})=r_{i}=d_{p_{i}}(B_{i})$ and 
$d_{p_{i}}(C\o(U_{\tilde{p_{i}}}(B_{i}\cap B_{3})))\leq r_{i}$, 
a contradiction to Lemma \ref{LemFusionBorel1}. 
\qed

Again $\PSL_{2}$ in positive characteristic offers counterexamples to 
Lemmas \ref{LemFusionBorel1} and \ref{LemFusionBorel2} 
when $r_{i}=0$, as we may take for $B_{1}$ and $B_{2}$ two distinct conjugates 
of the standard Borel subgroup $B$ and for $B_{3}$ any of these two. 

\subsubsection{Consequences on Fitting subgroups}\label{SectConsFittingSbgps}

The first paragraph of the proof of the following lemma appeared as 
\cite[Corollary 2.2]{BurdgesJalLemma}. 

\begin{lem}\label{FB1capFB2torsionfree}
Let $G$ be a locally$\o$ solvable$\o$ group of finite Morley rank. 
If $B_{1}$ and $B_{2}$ are two distinct Borel subgroups and 
$X$ denotes ${F(B_{1})\cap F(B_{2})}$, then 
$X\o$ is torsion free, $X=X\o\times S$ for a finite subgroup $S$, and for any 
subgroup $S_{1}$ of $X$ $C\o(S_{1})$ is nonsolvable if and only if 
$S_{1}\leq S$. 
\end{lem}
\proof
Assume $X\o$ not torsion free. Then it contains 
a nontrivial decent torus $T$ or a nontrivial $p$-unipotent subgroup $U$. In the first 
case, $T\leq {Z(B_{1})\cap Z(B_{2})}$ by Fact \ref{ActionOnLessUnip} $(1)$, 
$B_{1}=N\o(T)=B_{2}$ by local$\o$ solvability$\o$ and 
Lemma \ref{CharactBorelXLocal}, a contradiction. In the second case 
Corollary \ref{corUnicityBorelSuperMax} 
with $\tilde p=(p,\infty)$ and $X=U$ yields $B_{1}=B_{2}$, again a contradiction. 

We have now $X=X\o\times S$ for some finite subgroup $S$ of $X$ by 
Fact \ref{StrucNilpGroups}. 

Let $S_{1}$ be a subgroup of $X$. If $S_{1}\nleq S$, then $S_{1}$ contains an 
element of the form $s\cdot x$ for some $s$ in $S$ and some nontrivial 
element $x$ in $X\o$. As $X\o$ is torsion-free, $x$ as infinite order, as well as 
$s\cdot x$, and $C\o(S_{1})\leq C\o(H(s\cdot x))\leq 
N\o(H\o(s\cdot x))$, which is solvable by local$\o$ solvability$\o$ 
of $G$. Hence $C\o(S_{1})$ nonsolvable implies $S_{1}\leq S$. 

We now want to show that if $S_{1}\leq S$, then $C\o(S_{1})$ is nonsolvable. 
It suffices to do it for $S$, so we assume toward a contradiction $C\o(S)$ solvable. 
Let $B_{3}$ be a Borel subgroup of $G$ containing $C\o(S)$. 
Notice that the finite nilpotent group $S$ is the product of its Sylow $p$-subgroups. 
If $p_{1}$ and $p_{2}$ are two (not necessarily distinct) primes dividing the 
order of $S$, then we claim that one cannot have 
$U_{p_{1}}(B_{1})\neq 1$ and $U_{p_{2}}(B_{2})\neq 1$. 
Assume the contrary. As Sylow subgroups for 
primes different from $p_{1}$ in $F(B_{1})$ commute 
with $U_{p_{1}}(B_{1})$ (by Fact \ref{StrucNilpGroups} $(1)$!), 
$U_{p_{1}}(B_{1}\cap  C\o(S))$ is nontrivial by 
Fact \ref{StructpSylSolvGps} $(2)$. 
Similarly, $U_{p_{2}}(B_{2}\cap  C\o(S))$ is nontrivial. 
Now Lemma \ref{LemFusionBorel2} gives a contradiction, which proves our claim. 
It follows that all nontrivial $p$-unipotent subgroups of $B_{1}$ or $B_{2}$, for 
$p$ dividing the order of $S$, are on one side, say they are all in $B_{1}$. 
Notice then that all $p$-unipotence blocks of $B_{2}$, 
for $p$ dividing the order of $S$, are trivial. In particular $S\leq Z(B_{2})$ 
by Lemma \ref{FacptTorFitUpPerp}. Hence $B_{2}\leq C\o(S)\leq B_{3}$, 
$B_{3}=B_{2}$, and $C\o(S)=B_{2}$. 
Hence one cannot have $C\o(S)=B_{1}$, as $B_{1}\neq B_{2}$. 
Hence $S$ is not central in $B_{1}$. By Lemma \ref{FacptTorFitUpPerp}, 
there is a prime $p$ dividing the order of $S$ and such that $U_{p}(B_{1})\neq 1$. 
As above, $U_{p}(C\o_{B_{1}}(S))$ is nontrivial by Fact \ref{StrucNilpGroups} $(1)$ 
and Fact \ref{StructpSylSolvGps} $(2)$, 
and Corollary \ref{corUnicityBorelSuperMax} gives then 
$B_{1}=B_{2}$, a contradiction. 
\qed

\bigskip
A subgroup $S$ as in Lemma \ref{FB1capFB2torsionfree} 
could for example be the subgroup $Z(H)$ in the hypothetic 
Configuration \ref{BadConfig1}. 

We mention, parenthetically, that it is a version of the following lemma 
which has been baptized ``Jaligot's lemma" in \cite[\S2]{BurdgesJalLemma} 
(see \cite[\S3.4]{CherlinJaligot2004} and 
\cite[Lemma 2.1]{BurdgesJalLemma}). 

\begin{lem}\label{JaligLemma}
Let $G$ be a locally$\o$ solvable$\o$ group of finite Morley rank. Assume that, 
for $i=1$ and $2$, $\tilde p_{i}=(p_{i}, r_{i})$ are two unipotence parameters 
such that $d_{p_{i}}(G)=r_{i}$, and $B_{i}$ are two distinct Borel 
subgroups such that $d_{p_{i}}(B_{i})=r_{i}$. If $X$ denotes 
$F(B_{1})\cap F(B_{2})$, then $X$ is finite and 
$C\o(S_{1})$ is nonsolvable for any subgroup $S_{1}$ of $X$. 
\end{lem}
\proof
Assume $X\o$ non-trivial. By local$\o$ solvability$\o$ of $G$, 
$N\o(X)$ is solvable and hence contained in a Borel subgroup $B_{3}$ of $G$. 
As $X\o$ is torsion-free by Lemma \ref{FB1capFB2torsionfree}, the 
assumption that $d_{p_{i}}(B_{i})=r_{i}$ implies $r_{i}>0$ for each $i$. 
By Fact \ref{FactUnipHeaviest}, $U_{\tilde p_{i}}(B_{i})$ is in $F\o(B_{i})$, 
and by Fact \ref{StrucNilpGroups}, $\Gamma_{i}$, the last nontrivial iterated term of the 
descending central series of $U_{\tilde p_{i}}(B_{i})$, is central in $F\o(B_{i})$. 
Hence $\Gamma_{i}\leq N\o(X)\leq {B_{3}\cap B_{i}}$. 
Now by assumption each $\Gamma_{i}$ is nontrivial, and a 
$\tilde p_{i}$-group by Corollary \ref{Centrestildep}. 
Corollary \ref{corUnicityBorelSuperMax} implies that each $\Gamma_{i}$ 
is contained in a unique Borel subgroup of $G$, which gives $B_{3}=B_{1}$ 
and $B_{3}=B_{2}$, contradicting the assumption that $B_{1}\neq B_{2}$. 
Hence $X$ is finite. 

Our last claim is contained in Lemma \ref{FB1capFB2torsionfree}. 
\qed

\bigskip
In absence of local$\o$ solvability one might have $F(B_{1})\cap F(B_{2})$ 
(finite and) nontrivial in Lemma \ref{JaligLemma}, as for example in 
Configuration \ref{BadConfig1} again. 

\subsubsection{Consequences on soapy subgroups}\label{SectionConsOnSoapySubgp}

We continue as in Sections \ref{SectConseqBorelSngps} and \ref{SectConsFittingSbgps}
with consequences of the Uniqueness Theorem \ref{UniquenessLemma}, now 
on soapy subgroups. All these properties 
make us think of a soap sliding between two hands, exactly like 
a unipotent subgroup which cannot be contained in two distinct 
Borel subgroups in $\PSL_{2}$. 
The following lemmas will be used in our most critical computations in 
\cite{DeloroJaligotII}. 

\begin{lem}\label{propstronglywisegp}
Let $G$ be a locally$\o$ solvable$\o$ group of finite Morley rank, 
$B_{1}$ and $B_{2}$ two Borel subgroups each having a soapy subgroup 
$U_{1}$ and $U_{2}$ respectively. Then 
\begin{itemize}
\item[$(1)$]
$B_{1}$ is unique among Borel subgroups of $G$ containing $U_{1}$ and 
admitting the unipotence parameter of $U_{1}$ as maximal. 
\item[$(2)$]
If $[U_{1},U_{2}]=1$, then $B_{1}=B_{2}$. 
\end{itemize}
\end{lem}
\proof

$(1)$. By local$\o$ solvability$\o$ of $G$, $N\o(U_{1})$ is solvable. As 
$U_{1}$ is normal in $B_{1}$, the maximality of $B_{1}$ implies $N\o(U_{1})=B_{1}$. 
If the unipotence parameter of $U_{1}$ is $(\infty,0)$, then $B_{1}$ is a good 
torus, as well as any Borel subgroup admitting $(\infty,0)$ as maximal. So any 
such Borel subgroup is contained in $C\o(U_{1})=B_{1}$, and thus equal to 
$B_{1}$. Otherwise, as $C\o(U_{1})\leq N\o(U_{1})$, the first item 
is a mere application of Corollary \ref{corUnicityBorelMax}. 

$(2)$. Again $N\o(U_{1})=B_{1}$ and similarly $N\o(U_{2})=B_{2}$. 
Hence $U_{1}$, $U_{2} \leq {B_{1}\cap B_{2}}$ under the assumption 
that $U_{1}$ and $U_{2}$ commute. If $U_{1}$ is a good torus, 
then as for the first item $B_{1}$ is a good torus as well, as well as its subgroup $U_{2}$, 
and similarly $B_{2}$ also. We then get 
$B_{2}\leq C\o(U_{1})\leq N\o(U_{1})=B_{1}$, and equality of $B_{1}$ and $B_{2}$. 
One concludes symmetrically when $U_{2}$ is a good torus, 
so one can assume that both $U_{1}$ and $U_{2}$ are not good tori. 
As $U_{1}$, $U_{2}\leq{B_{1}\cap B_{2}}$, 
Corollary \ref{corUnicityBorelSuperMax} gives $B_{1}=B_{2}$ or 
$\max(d(U_{1}),d(U_{2}))<\infty$. In any case 
Corollary \ref{corUnicityBorelMax} gives $B_{1}=B_{2}$. 
\qed

\bigskip
The following lemma allows one to build soapy subgroups in presence of 
two Borel subgroups. 

\begin{lem}\label{getwise2}
Let $G$ be a locally$\o$ solvable$\o$ group of finite Morley rank, 
$B_{1}$ and $B_{2}$ two Borel subgroups, and $U_{1}$ a soapy subgroup 
of $B_{1}$. If $U_{1} \leq B_{2}$, then $B_{2}$ contains a 
characteristically soapy subgroup. 
\end{lem}
\proof
If $B_{1}=B_{2}$, then $U_{1}$ is a soapy subgroup of $B_{2}$ and we may 
use Lemma \ref{LemFromSoapyToCharSoapy}. 

Assume now $B_{1}\neq B_{2}$, and let $\tilde q_{1}$ be the unipotence parameter 
attached to $U_{1}$. Let $\tilde q_{2}$ be a unipotence parameter maximal for 
$B_{2}$. If $\tilde q_{2}=(\infty,0)$, then $B_{2}$ is a good torus, as well 
as $U_{1}$, as well as $B_{1}$, and then one concludes as usual by 
local$\o$ solvability$\o$ that $B_{1}=B_{2}$. Hence $\tilde q_{2}$ is not 
$(\infty,0)$. If $U_{\tilde q_{2}}(Z(F\o(B_{2})))$ is not central in 
$B_{2}$, then we may apply Lemma \ref{getwise1}. 

So now assume toward a contradiction 
$U_{\tilde q_{2}}(Z(F\o(B_{2})))$ central in $B_{2}$. 
In particular $U_{\tilde q_{2}}(Z(F\o(B_{2})))\leq C\o(U_{1})\leq 
N\o(U_{1})=B_{1}$ by local$\o$ solvability$\o$ of $G$. By 
Corollary \ref{corUnicityBorelSuperMax}, $\tilde q_{1}$ and 
$\tilde q_{2}$ do not represent subgroups of bounded exponent, as 
$B_{1}\neq B_{2}$. The maximality of $\tilde q_{1}$ for $B_{1}$ and 
of $\tilde q_{2}$ for $B_{2}$ then yields $\tilde q_{1}=\tilde q_{2}$. 
But Corollary \ref{corUnicityBorelMax} gives the uniqueness of $B_{2}$ among 
Borel subgroups containing $U_{\tilde q_{2}}(Z(F\o(B_{2})))$ and admitting 
$\tilde q_{2}$ as maximal. Thus $B_{1}=B_{2}$, a contradiction in the last 
case under consideration. 
\qed

\subsubsection{Consequences on Carter subgroups}

Theorem \ref{UniquenessLemma} also gives information on Carter subgroups 
possessing a subgroup sufficiently unipotent relatively to the ambient group. 

\begin{lem}\label{LemCarterWithHeavy}
Let $G$ be a locally$\o$ solvable$\o$ group of finite Morley rank, 
$Q$ a Carter subgroup of $G$ and $\tilde p=(p,r)$ a unipotence parameter admitted 
by $Q$. Assume $Q$ contains a nontrivial definable central $\tilde p$-subgroup $U_{1}$ 
with a nonempty subset $X$ such that $d_{p}(C(X))\leq r$. 
Then exactly one of the following three cases occur. 
\begin{itemize}
\item[$(1)$]
$Q$ is a generous Carter subgroup. 
\item[$(2)$]
For $g$ generic in $Q$, $d_{p}(C\o(g))>r$. 
\item[$(3)$]
The generic element of $Q$ is exceptional, and in particular any element of 
$Q$ has order at most $e(G)$. 
\end{itemize}
\end{lem}
\proof
Notice that the assumption together with Corollary \ref{Centrestildep} 
and Fact \ref{StrucNilpGroups} $(2)$ 
implies that $\tilde{p}$ is maximal in its characteristic for $Q$. 

If $Q$ is generous, then $C\o(g)\leq Q$ holds for $g$ generic in $Q$ 
by Fact \ref{CarterGeneiffgendis} 
(see also \cite[Lemma 3.10]{Jaligot06}), so cases $(2)$ and $(3)$ 
cannot occur. 

Assume $Q$ not generous in $G$. By Fact \ref{DecentTorusInGenCarter}, 
$Q$ contains no nontrivial good torus, and thus $r>0$ as $\tilde p$ is maximal 
in its characteristic for $Q$. By Theorem \ref{UniquenessLemma}, $U_{1}$ is 
contained in a unique Sylow $\tilde p$ subgroup of $G$, say $U$, and 
$Q\leq N\o(U)$. Notice that $N\o(U)$ is solvable by 
local$\o$ solvability$\o$ of $G$. By condition $(4)$ in 
Fact \ref{CarterGeneiffgendis} a generic element $g$ of $Q$ is in infinitely 
many conjugates of $Q$. 

Suppose toward a contradiction $d_{p}(C\o(g))\leq r$ and $C\o(g)$ solvable. 
Then $\tilde p \neq (\infty,0)$ is a unipotence parameter maximal 
in its characteristic for the definable connected 
solvable subgroup $C\o(g)$. It follows that $C\o(g)$ contains a unique 
Sylow $\tilde p$-subgroup by Fact \ref{FactUnipHeaviest}, 
which is necessarily a $\tilde p$-subgroup of 
$U$ as it contains $U_{1}$. 
If $\gamma$ is an element 
of $G$ such that $g\in Q^{\gamma}$, then $U_{1}$ and 
$U_{1}^{\gamma}$ are both contained in $U_{\tilde p}(C\o(g))$, and 
by uniqueness applied now to $U_{1}^{\gamma}$ one 
gets $U=U^{\gamma}$. Hence all $G$-conjugates of $Q$ containing $g$ 
are actually $N\o(U)$-conjugate. But now $Q$ is generous in 
the definable connected solvable subgroup $N\o(U)$, and thus a generic 
element of $Q$ is in a unique $N\o(U)$-conjugate of $Q$ by 
Fact \ref{CarterGeneiffgendis}. This is a contradiction. 
Hence when $Q$ is not generous one of the two cases $(2)$ or $(3)$ must occur. 

Notice that in case $(3)$ a generic element of $Q$, being exceptional, has order 
at most $e(G)$, and then the exponent of $Q$ is bounded by $e(G)$ by 
Fact \ref{StrucNilpGroups} $(2)$. 

It just remains to show that cases $(2)$ and $(3)$ cannot occur simultaneously. 
But in case $(2)$ $r$ cannot be $\infty$, and in case $(3)$ it must. 
\qed

\bigskip
Of course,  by Corollary \ref{Centrestildep}, 
Lemma \ref{LemCarterWithHeavy} applies when 
$d_{p}(G)=d_{p}(Q)=r$. In particular a nongenerous Carter subgroup which 
is not divisible must be as in case $(3)$ of Lemma \ref{LemCarterWithHeavy}. 

\subsection{Uniqueness Theorem and cosets}\label{SectionLocAnalandGen}

In \cite{CherlinJaligot2004} arguments pending on cosets and generosity 
were developed intensively for determining Weyl groups in groups of finite Morley rank 
in the specific case of minimal connected simple groups. 
This systematic approach was strongly inspired by the 
seminal work of Nesin in the context of bad groups \cite{Nesin89-a}. 
These arguments generally split into two parts. 
Cosets corresponding to an undesirable Weyl element are usually shown to be both 
generous {\em and} nongenerous in the ambient group, and then the coset as well as 
the unexpected Weyl element does not exist. 
Local properties of small groups 
often allow one to prove that some cosets are generous, as this is done intensively 
in \cite{CherlinJaligot2004}, and normally this is contradictory by itself. 

In the light of the fine analysis of generous sets of \cite{Jaligot06} and in continuation 
of this work, these {\em coset arguments} have certain generalizations, 
and what follows is part of it. 
It is however worth recalling these arguments in the specific 
context of locally$\o$ solvable$\o$ groups. The interest is both to put in a 
uniform format this essential content of \cite{CherlinJaligot2004} 
(here with the appeal to \cite{Jaligot06}), 
and to see how the specific local analysis of small groups originates 
further such arguments. In the process we will also encounter an interesting 
pathological configuration. 

As far as generosity is concerned, the fine analysis of conjugacy classes in 
\cite{Jaligot06} definitively provided the right understanding concerning generosity. 

\begin{fait}\label{FactGPermRanks}
{\bf \cite{Jaligot06}}
Let $(G,\Omega)$ be a permutation group of finite Morley rank in which the Morley rank 
is additive (or a ranked permutation group), $H$ a definable 
subset of $\Omega$, and assume that for 
$r$ between $0$ and $\rk(G/N(H))$ the definable set $H_{r}$, consisting of those elements 
of $H$ contained in a set of conjugates of $H$ of rank exactly $r$, is nonempty. Then 
$$\rk({H_{r}}^{G})=\rk(G)+\rk(H_{r})-\rk(N(H))-r.$$ 
\end{fait}
\proof
This is essentially the content of the fine analysis of conjugacy classes of 
\cite[\S2.2]{Jaligot06}. Here the geometric proof for this mentioned later 
by Cherlin yields this equality exactly as in \cite[\S2.3]{Jaligot06}. One uses 
the additivity of the Morley rank, or of the rank function if the structure is ranked 
as in the axioms of \cite{BorovikNesin(Book)94}, for computing the ranks 
of the set of flags in the naturally associated geometry. 
\qed

\bigskip
We say that a connected locally$\o$ solvable$\o$ 
group $H$ of finite Morley rank is {\em sick} if it contains 
a generous Carter subgroup, $H$ contains no nontrivial decent torus, 
the generous Carter subgroup has a nontrivial $\pi$-unipotent subgroup, 
$H$ does not conjugate its maximal $p$-unipotent subgroups for any $p\in \pi$, 
and for any such $p$ $N\o_{H}(U)$ is a Carter subgroup of bounded exponent of $H$ 
for some maximal $p$-unipotent subgroup $U$ of $H$. 

\begin{theo}\label{TheoProperCosetsNonGen}
Let $G$ be a locally$\o$ solvable$\o$ group of finite Morley rank with a generous 
Carter subgroup. Let $H$ be a definable connected subgroup of $G$ and 
$x$ an element of $N_{G\o}(H)$ not in $H$. Assume that $H$ is solvable or that $H$ 
contains a generous Carter subgroup of the ambient 
group and is not sick. Then $xH$ is not generous in $G$. 
\end{theo}

This essential content of \cite{CherlinJaligot2004} was proved locally, 
usually for $H$ a Carter subgroup or indeed 
the centralizer$\o$ of a torus of the ambient group, and we reformat 
it in its natural form here, replacing the applications of 
\cite[Proposition 3.11]{CherlinJaligot2004} in that paper 
by Section \ref{SectionJaligLemma} here and arguing directly for the production 
of bounded exponent torsion of \cite[\S3.3]{CherlinJaligot2004}. 

\bigskip

\proof
Assume towards a contradiction $xH$ is generous in $G$. 

Notice that for the suitable $r$ as in Fact \ref{FactGPermRanks} such that 
$X_{r}^{G}$ is generic in $G$, where $X$ denotes the coset $xH$, 
one has $\rk(X_{r})-\rk(N(X))=r\geq 0$, hence 
$\rk(H)\leq \rk(N(X))\leq \rk(X_{r})\leq \rk(xH)=\rk(H)$, and thus 
$r=0$, $\rk(H)=\rk(N(xH))=\rk([xH]_{0})$, and a generic element of 
$xH$ is in only finitely many conjugates of $xH$. Also $N\o(xH)=H$. 
(The argument of this paragraph is of course general.) 

A generic element $w$ of $xH$ is also generic in $G$. 
By Fact \ref{FactGPermRanks}, $w$ is in only finitely many conjugates of $xH$, 
and this implies as in \cite[Fundamental Lemma 3.3]{Jaligot06} that 
$C\o(w)\leq N\o(xH)=H$. By Fact \ref{CarterGeneiffgendis} 
(see also \cite[Lemma 3.10]{Jaligot06}), $C\o(w)\leq Q_{w}$, where $Q_{w}$ denotes 
the unique conjugate of the generous Carter subgroup $Q$ containing $w$. 

(The next paragraph corresponds to the local applications 
of \cite[\S3.3]{CherlinJaligot2004} in that paper, though things may be stated in 
somewhat reversed ways there. The paragraph following it 
will then concern the application of the uniqueness theorems of 
Section \ref{SectionJaligLemma}, and this was usually done 
first in the sequence of argumentations in \cite{CherlinJaligot2004} via 
\cite[Proposition 3.11]{CherlinJaligot2004}. Actually 
\cite[Proposition 3.11]{CherlinJaligot2004} provided trivial intersections at 
the level of subgroups, and then cosets consisting generically of bounded exponent 
elements, and then \cite[\S3.3]{CherlinJaligot2004} gave bounded exponent subgroups.) 

Let $n$ be the order of $x$ modulo $H$. By assumption $n>1$. 
As $w\in xH$ and $x$ normalizes $H$, $H(w)\leq {\<x\>}H$ and in the definable 
hull $H(w)$, $w$ has order a nontrivial multiple of $n$ modulo $H\o(w)$. 
This shows that the generic element of $Q_{w}$ has the property of having 
order a nontrivial multiple of $n$ modulo the connected 
component of its definable hull. Hence $Q_{w}$ contains by 
Fact \ref{StrucNilpGroups} a nontrivial definable connected subgroup of bounded 
exponent (whose elements are generically of order the above multiple of $n$). 

The generic element $w$ of $Q_{w}$ centralizes a nontrivial definable 
connected abelian subgroup of exponent $n$ by 
Facts \ref{StrucNilpGroups} $(2)$ and \ref{StructpSylSolvGps} $(2)$.  
By Corollary \ref{corUnicityBorelMax}, 
$C\o(w)$ is contained in a unique Borel subgroup, say $B_{w}$, and 
$C\o(w)\leq Q_{w}\leq B_{w}$. 

Now $wH$ is generous in any definable connected subgroup containing it. 
This is a general fact, for which one can proceed as in \cite[Lemma 3.9 b.]{Jaligot06}. 
Indeed the property of the generic element of $xH$ of being contained in finitely many 
conjugates of $xH$ is obviously preserved when passing to definable subgroup, 
and this suffices with Fact \ref{FactGPermRanks} and the fact that 
$\rk(N(xH))=\rk(xH)$. 

When $H$ is solvable, we have as $C\o(w)\leq H$ also that $H\leq B_{w}$. 
In particular $wH$ is generous in the connected solvable group $B_{w}$, 
and this is ridiculous. One can argue, being inside a connected solvable group. 
One can also argue noticing that $N\o_{Q_{w}}(\<w\>(H\cap Q_{w}))$ normalizes 
$w(H\cap Q_{w})$ by Fact \ref{ActionOnFiniteSets}, hence normalizes 
$wH$ as in \cite[Fundamental Lemma 3.3]{Jaligot06}, so it is in 
$N\o(wH)=H$, and the normalizer condition in connected nilpotent groups gives 
$w\in \<w\>(H\cap Q_{w})=Q_{w}$, a contradiction as $Q_{w}$ is connected and 
$w$ is not in $(H\cap Q_{w})$. 

This finishes our rearrangement of \cite{CherlinJaligot2004} when $H$ is 
solvable, and when $H$ contains a generous Carter subgroup $Q$ one can proceed as 
follows. 

By the preceding case one may assume $H$ nonsolvable. By a Frattini Argument following 
from the conjugacy of generous Carter subgroups, \cite[Corollary 3.13]{Jaligot06}, 
we may suppose that $x$ normalizes the generous Carter subgroup $Q$ of $H$. 

Let $\pi$ be the set of primes involved in the bounded exponent 
part of the generous Carter subgroup. Assume that $H$ conjugates its maximal 
$p$-unipotent subgroups for some prime $p$ in $\pi$, or that $Q$ contains a nontrivial 
decent torus. In the first case one may assume after $H$-conjugacy, with 
Fact \ref{FactUnipHeaviest} and Lemma \ref{LemFusionBorel2} that 
$Q$ and $Q_{w}$ are in a common Borel subgroup of $G$. 
Similarly, if $Q$ and $Q_{w}$ contain a nontrivial decent torus, one may 
assume after $H$-conjugacy that $Q$ and $Q_{w}$ are the centralizer$\o$ of 
their common maximal decent torus, and hence that they are by 
local$\o$ solvability$\o$ in a same Borel subgroup of $G$. 
If we denote by $B$ this Borel subgroup in both cases, then we get 
$Q\leq (H\cap B)\o < \<w\>(H\cap B)\o \leq B$, 
and this is impossible by a Frattini argument as the Carter subgroup $Q$ is selfnormalizing 
in $B$. 

This leaves us with the case in which $Q$ has a nontrivial $\pi$-unipotent 
subgroup, $H$ does not conjugate its maximal $p$-unipotent subgroups for any 
$p\in \pi$, and $H$ contains no nontrivial decent torus. Hence in this 
pathological situation $H$ has all the symptoms of sickness, except maybe the last one. 
But this will be seen in Lemma \ref{LemmaNontrivialpUnipoExists} below 
(whose proof will be independent). 
\qed

\bigskip
Theorem \ref{TheoProperCosetsNonGen} represents coset arguments 
of \cite{CherlinJaligot2004} for dealing with Weyl groups. 
We note that its proof actually provided the following much more general 

\begin{fait}\label{FactCosetsNeverGen}
{\bf \cite{JaligotGenerixCosets}}
Let $G$ be a group of finite Morley rank in which the generic element of 
$G\o$ is in a connected nilpotent subgroup, and let $H$ be a definable 
subgroup of $G\o$. Then $H\setminus H\o$ is not generous in $G$. 
\end{fait}

The main consequence of Fact \ref{FactCosetsNeverGen} 
is the following, a general fact in which its conclusion is true, which 
also recasts some corresponding consequences as in \cite{CherlinJaligot2004} 
somehow in their original content. 

\begin{fait}\label{FactConseqGenWeylGpsProperCent}
{\bf \cite{JaligotGenerixCosets}}
Let $G$ be a group of finite Morley rank, $n$ a natural number, 
$H$ a definable connected generous subgroup with the property that, for $h$ generic 
in $H$, $h$ is in a connected nilpotent subgroup of $H$ and $h^{n}$ is also 
generic in $H$, and assume $w$ is an element of $G\o$ of finite order $n$ normalizing $H$ 
without being inside. Then $C_{H}(w)<H$.
\end{fait}

Groups of finite Morley rank with a generous Carter subgroup not divisible or not 
abelian can be dealt with the Bender method, the results 
of Sections \ref{SectionJaligLemma} and \ref{SectMaxPaire} here, 
and otherwise Fact \ref{FactConseqGenWeylGpsProperCent} applies. 
In a locally$\o$ solvable$\o$ context and in presence 
of a generous Carter subgroup the situation will be considered in a separate 
paper with the results of Section \ref{SectBoundExpCase} below. 
Here we merely mention the following basic commutation principle relevant 
for Weyl groups, and which builds upon \cite[Lemme 3.1]{Deloro05}. 

\begin{lem}\label{LemCommutPrincipleForWeyl}
Let $G=NQ$ be a group, with $N$ and $Q$ two subgroups and $N$ normal. 
Assume $\sigma$ is an automorphism of $G$ normalizing $Q$ and 
fixing $N$ pointwise. Then 
\begin{itemize}
\item[$(1)$]
$N$ and $\<[\sigma,Q]\>$ commute. 
\item[$(2)$]
If $Q=\<[\sigma,Q]\>C_{Q}(\sigma)$ and $N\leq N(C_{Q}(\sigma))$, then 
$N\leq N(Q)$. 
\end{itemize}
\end{lem}
\proof
For any element $q$ in $Q$ and $h$ in $N$ one has 
$$h^{[\sigma,q]}=h^{q^{-1}\sigma q}=(h^{q^{-1}})^{\sigma q}=h^{q^{-1}q}=h$$
and thus $h\in C([\sigma,q])$. Hence $N\leq C(\<[\sigma,Q]\>)$, 
the general commutation principle of \cite[Lemme 3.1]{Deloro05}. 
The second item follows. 
\qed 

\bigskip
In particular, if in Lemma \ref{LemCommutPrincipleForWeyl} $(2)$ 
$G=NQ$ has finite Morley rank and $Q$ is a Carter subgroup, then 
$N\o\leq C\o_{Q}(\sigma)$. 

We finish this section with one word about centralizers of definable connected 
exceptional subgroups. If $G$, $H$, and $n$ are as in 
Fact \ref{FactConseqGenWeylGpsProperCent}, with $G$ 
locally$\o$ solvable$\o$ and $H$ an exceptional 
definable connected nonsolvable subgroup of $G$, then $C(H)$ is 
finite by Lemma \ref{LemUnifBoundException}, and if $x$ is an 
element of $G\o$ in this finite centralizer and of order $n$, then 
Fact \ref{FactConseqGenWeylGpsProperCent} implies that $x$ is in $Z(H)$. 

\subsection{Maximal pairs of Borel subgroups}\label{SectMaxPaire}

When the absolute maximality assumptions concerning unipotence degrees fail in 
Lemma \ref{JaligLemma} one might have (or rather cannot exclude) pairs 
of Borel subgroups whose Fitting subgroups have an infinite intersection. This 
situation has been studied intensively in \cite{BurdgesJalLemma}. 
In what follows, not only we claim no originality compared to this paper, 
but also we will tend to follow it word by word. 
The only differences will appear in the notation used for unipotence parameters and 
in a special care needed for dealing here with our weakest assumption of 
local$\o$ solvability$\o$. Some additional results from \cite{Deloro05} will 
be mentionned in the process. 

\begin{defi}
Let $G$ be a group of finite Morley rank, $B_{1}$ and $B_{2}$ two 
distinct Borel subgroups. We say that $(B_{1},B_{2})$ is a {\em maximal pair} 
(of Borel subgroups) if the definable connected subgroup 
$(B_{1}\cap B_{2})\o$ is maximal for inclusion among all definable connected subgroups 
of the form $(L_{1}\cap L_{2})\o$, with $L_{1}$ and $L_{2}$ two 
distinct Borel subgroups of $G$. 
\end{defi}

\begin{hyp}\label{HypMaxPairXneq1}
{\bf \cite[Hypothesis 3.2]{BurdgesJalLemma}}
We assume the following configuration:
\begin{itemize}
\item[$(1)$]
$G$ is a locally$\o$ solvable$\o$ group of finite Morley rank. 
\item[$(2)$]
$(B_{1},B_{2})$ is a maximal pair of Borel subgroups of $G$. 
\item[$(3)$]
$[F(B_{1})\cap F(B_{2})]\o$ is nontrivial. 
\end{itemize}
\end{hyp}

\begin{notation}
{\bf \cite[Notation 3.3]{BurdgesJalLemma}}
We let 
\begin{itemize}
\item[$(1)$]
$H=(B_{1}\cap B_{2})\o$. 
\item[$(2)$]
$X=F(B_{1})\cap F(B_{2})$. 
\item[$(3)$]
$r'=d_{\infty}(X)$. 
\end{itemize}
\end{notation}

Recall that $X\o$ is torsion free by Lemma \ref{FB1capFB2torsionfree}. In particular 
$0<r'<\infty$. 
In particular $0<d_{\infty}(B_{1})<\infty$ and $0<d_{\infty}(B_{2})<\infty$. 

Notice that by Lemma \ref{JaligLemma} one cannot have 
$d(B_{1})=d(B_{2})=\infty$. So at least one of the two Borel subgroups 
$B_{1}$ and $B_{2}$, say $B_{i}$, has no bounded exponent subgroup. 
In particular $0<d(B_{i})<\infty$. The other Borel subgroup $B_{i+1}$ 
might satisfy $0<d(B_{i+1})\leq \infty$ (this latter inequality will be shown 
to be also strict in the analysis below). 

\subsubsection{Homogeneity of $X$}

We observe that $H'\leq X\trianglelefteq H$. We will show the asymmetry of 
the situation, i.e. $d_{\infty}(B_{1})\neq d_{\infty}(B_{2})$. We may assume 
in any case that 

\begin{hyp}\label{HypdB1<dB2}
$d_{\infty}(B_{2})\leq d_{\infty}(B_{1})$. 
\end{hyp}

\bigskip
and we will indeed show that $d_{\infty}(B_{2})<d_{\infty}(B_{1})$. 
Notice that $d(H)=d_{\infty}(H)$. 

\begin{lem}\label{HypdB2leqdB1}
{\bf \cite[Lemma 3.5]{BurdgesJalLemma}}
$d_{\infty}(H)<d_{\infty}(B_{1})$. 
\end{lem}
\proof
As there is no nontrivial $p$-unipotent subgroup in $H$, $d(H)<\infty$. 

Suppose toward a contradiction $d_{\infty}(H)\geq d_{\infty}(B_{1})$. 
As $H\leq B_{1}$, $d_{\infty}(H)\leq d_{\infty}(B_{1})$ in any case, 
so our assumption becomes $d(H)=d(B_{1})$. 
Since $d(H)\leq d_{\infty}(B_{2})\leq d_{\infty}(B_{1})$ 
by Hypothesis \ref{HypdB2leqdB1}, all these unipotence degrees are equal to a 
certain $d$, and 
$U_{(\infty, d)}(H)\leq U_{(\infty , d)}(B_{1})\cap U_{(\infty ,d)}(B_{2})$. 
As $G$ is locally$\o$ solvable$\o$, $N\o(U_{(\infty,d)}(H))$ is solvable, and thus 
contained in a Borel subgroup $B_{3}$ of $G$. 

Now we contradict the fact that $B_{1}\neq B_{2}$. 

If $U_{(\infty,d)}(H)=U_{(\infty,d)}(B_{i})$ for some $i=1$ or $2$, 
then by local$\o$ solvability$\o$ and maximality of $B_{i}$, 
$B_{i}=N\o(U_{(\infty,d)}(H))\leq B_{3}$, 
and $B_{i}=B_{3}$. 

If $U_{(\infty,d)}(H)<U_{(\infty,d)}(B_{i})$ for some $i=1$ or $2$, then, 
as $d_{\infty}(B_{i})=d\geq 1$, $U_{(\infty,d)}(B_{i})\leq F\o(B_{i})$ and is in 
particular nilpotent, Fact \ref{NormCond} gives 
$$U_{(\infty,d)}(H)<
U_{(\infty,d)}(N_{U_{(\infty,d)}(B_{i})}(U_{(\infty,d)}(H))) \leq B_{3}.$$
Since $U_{(\infty,d)}(H)\trianglelefteq H$, we must get $H<(B_{i}\cap B_{3})\o$. 
By maximality of $H$ we get $B_{i}=B_{3}$. 

As $B_{1}\neq B_{2}$, $U_{(\infty,d)}(H)$ is proper in one of the two 
subgroups $U_{(\infty,d)}(B_{i})$, and not in both. 
In any case we get $B_{1}=B_{3}=B_{2}$, a 
contradiction. 
\qed

\begin{lem}\label{LemdH=dB2}
{\bf \cite[Lemma 3.6]{BurdgesJalLemma}}
$d(H)=d_{\infty}(B_{2})$. 
\end{lem}
\proof
Suppose toward a contradiction $d(H)<d_{\infty}(B_{2})$. 

By local$\o$ solvability$\o$, $N\o(U_{(\infty,r')}(X))$ is solvable, 
and contained in a Borel 
subgroup $B_{3}$ of $G$. Since $U_{(\infty,r')}(X)\trianglelefteq H$, 
$H$ is contained 
in $B_{3}$. Since $d(H)< d_{\infty}(B_{i})$ for $i=1$ and $2$ by 
Lemma \ref{HypdB2leqdB1} and Hypothesis \ref{HypdB1<dB2}, 
Fact \ref{StrucNilpGroups} gives 
$U_{(\infty,d_{\infty}(B_{i}))}(B_{i})\leq C\o(U_{(\infty,r')}(X))\leq B_{3}$. 
Hence 
$H<(B_{i}\cap B_{3})\o$ and $B_{i}=B_{3}$ by maximality of $H$ for 
$i=1$ and $2$, a contradiction to $B_{1}\neq B_{2}$. 
\qed

\begin{cor}\label{corAsymMaxPair}
$d(H)=d_{\infty}(H)=d_{\infty}(B_{2})<d_{\infty}(B_{1})$. 
\end{cor}

\begin{prop}\label{HnonabB1B2unique}
{\bf \cite[Proposition 3.7]{BurdgesJalLemma}}
If $H$ is nonabelian, then $B_{1}$ and $B_{2}$ are the only 
Borel subgroup containing $H$. 
\end{prop}
\proof
Suppose there is a Borel subgroup $B_{3}$ distinct from $B_{1}$ and $B_{3}$ and 
containing $H$. The maximality of $H$ yields 
$H=(B_{1}\cap B_{3})\o=(B_{2}\cap B_{3})\o$. 
Since $1\neq H'\leq F\o(B_{3})$, the maximal pairs 
$(B_{1},B_{3})$ and $(B_{2},B_{3})$ satisfy Hypothesis \ref{HypMaxPairXneq1}. 
Since $d_{\infty}(H)<d_{\infty}(B_{1})$ by Lemma \ref{HypdB2leqdB1}, 
$d_{\infty}(H)=d_{\infty}(B_{3})$ by 
Lemma \ref{LemdH=dB2} applied to the maximal pair $(B_{1},B_{3})$. 
But since $d_{\infty}(H)=d_{\infty}(B_{2})$ by Lemma \ref{LemdH=dB2}, 
$d_{\infty}(H)<d_{\infty}(B_{3})$ by Lemma \ref{HypdB2leqdB1} 
applied to the maximal pair 
$(B_{2},B_{3})$. This is a contradiction. 
\qed

\bigskip
As $d_{\infty}(B_{2})<d_{\infty}(B_{1})$, the 
Borel subgroups are not conjugate, a point we exploit in the next lemma. 

\begin{lem}\label{FBiNotinH}
{\bf \cite[Lemma 3.8]{BurdgesJalLemma}}
$F\o(B_{i})\nleq H$ for $i=1$ and $2$. 
\end{lem}
\proof
Since $d_{\infty}(H)<d_{\infty}(B_{1})$ by 
Lemma \ref{HypdB2leqdB1}, $F\o(B_{1})\nleq H$. 
Suppose toward a contradiction $F\o(B_{2})\leq H$. Then $H\trianglelefteq B_{2}$ 
by Fact \ref{GConSolvG/FGdivAb}, and $H\leq B_{1}\cap B_{1}^{g}$ for some 
$g\in B_{2}\setminus N(B_{1})$. By maximality of $H$, $(B_{1},B_{1}^{g})$ is a 
maximal pair, and Corollary \ref{corAsymMaxPair} applied to this maximal pair 
gives a contradiction. 
\qed

\begin{lem}\label{LemX1GenNX1}
{\bf \cite[Lemma 3.9]{BurdgesJalLemma}}
If $X_{1}$ is an infinite definable subgroup of $X$ normal in $H$, then 
$N\o(X_{1})\leq B_{1}$.  
\end{lem}
\proof
By local$\o$ solvability$\o$ $N\o(X_{1})$ is solvable, and hence contained in a Borel 
subgroup $B_{3}$ of $G$. By assumption $H\leq B_{3}$. Since 
$d_{\infty}(H)<d_{\infty}(B_{1})$ 
by Lemma \ref{HypdB2leqdB1}, Fact \ref{StrucNilpGroups} yields 
$U_{(\infty,d_{\infty}(B_{1}))}(B_{1})\leq C\o(X_{1})\leq B_{3}$. Thus 
$H<(B_{1}\cap B_{3})\o$ and by maximality of $H$ we get $B_{1}=B_{3}$. 
In particular $N\o(X_{1})\leq B_{1}$. 
\qed

\begin{cor}
$[X\cap Z(F(B_{2}))]\o=1$. 
\end{cor}

\begin{theo}\label{TheoXHomogene}
{\bf \cite[Theorem 3.10]{BurdgesJalLemma}}
$X\o$ is a homogeneous $(\infty,r')$-sub\-group. 
\end{theo}
\proof
Recall that $X\o$ is torsion-free. Suppose toward a contradiction that  
$U_{(\infty,r)}(X)$ is nontrivial for some $0\leq r<r'$. By 
Fact \ref{StrucNilpGroups} and Lemma \ref{LemX1GenNX1} 
$$F\o(B_{2})\leq C\o(U_{(\infty,r')}(X))C\o(U_{(\infty,r)}(X))\leq B_{1}$$
and it follows that $F\o(B_{2})\leq H$. This contradicts Lemma \ref{FBiNotinH}. 
Hence $X\o$ is homogeneous in the maximal unipotence parameter in its 
characteristic, that is $(\infty,r')$. 
\qed

\subsubsection{Fitting subgroup of $B_{2}$}

We delineate now $F\o(B_{2})$, and in particular determine which of its factors 
are contained in $H$. 

\begin{lem}\label{FactdeFB2dansH}
{\bf \cite[Lemmas 3.11 and 3.12]{BurdgesJalLemma}}
$F\o(B_{2})$ is divisible (in particular $d(B_{2})=d_{\infty}(B_{2})$) 
and $U_{(\infty,r)}(F\o(B_{2}))\leq Z(H)$ 
when $0\leq r\leq d(B_{2})$ and $r\neq r'$.
\end{lem}
\proof
If $d(B_{2})=\infty$, then $U_{p}(B_{2})$ is nontrivial for some 
prime $p$, and contained in $N\o(X\o)\leq B_{1}$ by 
Lemma \ref{LemX1GenNX1}. This is a contradiction to 
Lemma \ref{FB1capFB2torsionfree} or Lemma \ref{JaligLemma}. 
Hence $d(B_{2})=d_{\infty}(B_{2})$ and $F\o(B_{2})$ is divisible. 

By Theorem \ref{TheoXHomogene}, Fact \ref{StrucNilpGroups}, and 
Lemma \ref{LemX1GenNX1}, 
$U_{(\infty,r)}(F\o(B_{2}))\leq C\o(X\o)\leq N\o(X\o)\leq B_{1}$, 
and hence each of these groups is contained in $H$. As each such group is nilpotent and 
normalized by the subgroup $H$ of $B_{2}$, each such group is in $F\o(H)$. 
Now Fact \ref{CommuttildpGp} and Theorem \ref{TheoXHomogene} give 
$$[H,U_{(\infty,r)}(F\o(B_{2}))]\leq U_{(\infty,r)}(H')\leq U_{(\infty,r)}(X\o)=1.$$
\qed

\begin{lem}\label{Fr'B2notinH}
{\bf \cite[Lemma 3.13]{BurdgesJalLemma}}
$U_{(\infty,r')}(F\o(B_{2}))$ is not contained in $H$ and not abelian. 
\end{lem}
\proof
By Fact \ref{StrucNilpGroups} $F\o(B_{2})$ is generated by its 
Sylow $\tilde p$-subgroups. But by Lemma \ref{FBiNotinH} $F\o(B_{2})\nleq H$, 
so Lemma \ref{FactdeFB2dansH} implies $U_{(\infty,r')}(F\o(B_{2}))\nleq H$. 

Since $N\o(X\o)\leq B_{1}$ by Lemma \ref{LemX1GenNX1}, 
$U_{(\infty,r')}(F\o(B_{2}))$ cannot be abelian. 
\qed

\bigskip
Now we can deduce from the two preceding lemmas 
that the unipotence degree $r'$ is uniquely determined by the structure of $B_{2}$. 

\begin{cor}
{\bf \cite[Corollary 3.14]{BurdgesJalLemma}}
$U_{(\infty,r)}(F\o(B_{2}))$ is not abelian if and only if $r=r'$. 
\end{cor}

\begin{lem}\label{Lem3.15}
{\bf \cite[Lemma 3.15]{BurdgesJalLemma}}
$U_{(\infty,r)}(B_{2})\leq F\o(B_{2})$ for every $r>r'$. 
\end{lem}
\proof
Let $Q$ be a definable $(\infty,r)$-subgroup of $B_{2}$. 
By Fact \ref{ActionOnLessUnip} and our assumption $r>r'$, 
$U_{(\infty,r')}(F\o(B_{2}))\cdot Q$ is nilpotent. It follows by 
Fact \ref{StrucNilpGroups} that $Q$ centralizes $U_{(\infty,r')}(F\o(B_{2}))$, and 
in particular $Q$ centralizes its subgroup $X\o$. Hence $Q\leq N\o(X\o)\leq B_{1}$ 
by Lemma \ref{LemX1GenNX1}, and $Q\leq H$. By Lemma \ref{FactdeFB2dansH}, 
$Q$ centralizes all factors of $F\o(B_{2})$, except maybe the one of unipotence parameter 
$(\infty,r')$. But as $U_{(\infty,r')}(F\o(B_{2}))\cdot Q$ is nilpotent, 
$F\o(B_{2})\cdot Q$ is nilpotent, and as it is normal in $B_{2}$ by 
Fact \ref{GConSolvG/FGdivAb}, we deduce that $Q\leq F\o(B_{2})$, as desired. 
\qed

\subsubsection{Structure of $H$}

\begin{lem}\label{U0r'H<FB2}
{\bf \cite[Lemma 3.16]{BurdgesJalLemma}}
$U_{(\infty,r')}(H)\leq F\o(B_{2})$. In particular $U_{(\infty,r')}(H)$ is nilpotent 
and the unique Sylow $(\infty,r')$-subgroup of $H$. 
\end{lem}
\proof
Let $Q$ be any definable $(\infty,r')$-subgroup of $H$. 
By Fact \ref{ActionOnLessUnip}, the group 
$U_{(\infty,r')}(F\o(B_{2}))\cdot Q$ is nilpotent. 
For any integer $r\neq r'$, $Q$ centralizes 
$U_{(\infty,r)}(F\o(B_{2}))$ by Lemma \ref{FactdeFB2dansH}. Hence 
$F\o(B_{2})\cdot Q$ is nilpotent by Fact \ref{StrucNilpGroups}. 
By Fact \ref{GConSolvG/FGdivAb}, this product is normal in $B_{2}$, and hence 
it must be contained in $F\o(B_{2})$. 

In particular $U_{(\infty,r')}(H)$ is nilpotent, and the unique Sylow $(\infty,r')$-subgroup 
of $H$. 
\qed

\bigskip
This has a consequence purely in $B_{2}$, as seen in \cite{Deloro05}. 

\begin{lem}\label{CarterB1B2UB2Nilp}
{\bf \cite[Lemme 3.11]{Deloro05}}
If a Carter subgroup of $H$ is also a Carter subgroup of $B_{2}$, then 
$U_{(\infty,r')}(B_{2})$ is nilpotent, included in $F\o(B_{2})$, and the unique 
Sylow $(\infty,r')$-subgroup of $B_{2}$. 
\end{lem}
\proof
Let $Q$ be a Carter subgroup of $H$, which is also a Carter subgroup of $B_{2}$. 
Then $U_{(\infty,r')}(F\o(B_{2}))\cdot U_{(\infty,r')}(Q)$ is a Sylow 
$(\infty,r')$-subgroup of $B_{2}$ by Fact \ref{StrucSylowUSolv}. By conjugacy of such 
subgroups in $B_{2}$, Fact \ref{TheoConjPiCarter}, 
it suffices to show that it is contained 
in $F\o(B_{2})$. But the first factor clearly is, 
and the second also by Lemma \ref{U0r'H<FB2}. 
\qed

\bigskip
We return to the structure of $H$. 

\begin{notation}
Let $Y=U_{(\infty,r')}(H)$ be the unique Sylow $(\infty,r')$-subgroup of $H$. 
It is normal in $H$, and in $F\o(H)$. 
\end{notation}

\bigskip
We find that $Y$ has properties antisymmetric to those of $X$. 

\begin{lem}\label{Lem3.17}
{\bf \cite[Lemma 3.17]{BurdgesJalLemma}}
$N\o(Y)\leq B_{2}$ and $X\o<Y$. 
In addition $U_{(\infty,r')}(N_{F(B_{2})}(Y))\nleq H$. 
\end{lem}
\proof
Let $P=U_{(\infty,r')}(F\o(B_{2}))$. Then $Y\leq P$ by Lemma \ref{U0r'H<FB2}. 
By Lemma \ref{Fr'B2notinH} $P\nleq H$, so $Y<P$. By Normalizer Condition, 
Fact \ref{NormCond}, $Y<U_{(\infty,r')}(N\o_{P}(Y))$. Now $X\o<Y$ by 
Lemma \ref{LemX1GenNX1}, and $N\o(Y)\leq B_{2}$ by maximality of $H$. 
\qed

\begin{theo}\label{ThmSylIntAb}
{\bf \cite[Theorem 3.18]{BurdgesJalLemma}}
Every definable connected nilpotent subgroup of $H$ is abelian. 
\end{theo}
\proof
By Fact \ref{StrucNilpGroups} it suffices to show that any Sylow $\tilde p$-subgroup 
of $H$ is abelian. As this is true for decent tori and there is no nontrivial $p$-unipotent 
subgroup in $H$, it suffices to show this when $\tilde p=(p,r)$, with $1\leq r<\infty$. 
For $r\neq r'$, $U_{(\infty,r)}(H')=1$ by Theorem \ref{TheoXHomogene}, so a 
Sylow $(\infty,r)$-subgroup of $H$ must be abelian by Fact \ref{CommuttildpGp}. It 
remains to show that the unique Sylow $(\infty,r')$-subgroup of $H$, 
Lemma \ref{U0r'H<FB2}, is also abelian. For if $Y'$ is not trivial, 
then $N\o_{B_{2}}(Y)\leq N\o_{B_{2}}(Y')\leq B_{1}$ 
by Lemma \ref{LemX1GenNX1}, contradicting Lemma \ref{Lem3.17}. 
This completes the proof. 
\qed

\begin{lem}\label{BenderHNonAbNH=H}
{\bf \cite[Lemma 3.19]{BurdgesJalLemma}}
If $H$ is not abelian, then $N\o(H)=H$. 
\end{lem}
\proof 
Lemma \ref{Lem3.17} implies that $N\o(H)\leq N\o(Y)\leq B_{2}$. When $H'\neq 1$, 
Lemma \ref{LemX1GenNX1} implies also that $N\o(H)\leq N\o(H')\leq B_{1}$. 
\qed

\subsubsection{Structure of $B_{1}$}

\begin{lem}\label{LemFB1DivBender}
{\bf \cite[Lemma 3.20]{BurdgesJalLemma}}
$F\o(B_{1})$ is divisible (and in particular $d(B_{1})=d_{\infty}(B_{1})$)
and $U_{(\infty,0)}(F\o(B_{1}))\leq Z\o(H)$. 
\end{lem}
\proof
$Y$ is an $(\infty,r')$-group, with $0<r'<\infty$, normalizing an hypothetic 
$p$-unipotent subgroup of $B_{1}$. By Fact \ref{ActionOnLessUnip} $(4)$, 
it centralizes all of them. As $N\o(Y)$ is solvable by local$\o$ solvability$\o$ 
of $G$, one gets if $U_{p}(B_{1})\neq 1$ for some prime $p$ that 
$N\o(Y)\leq B_{1}$ by Lemma \ref{corUnicityBorelSuperMax}. But 
$N\o(Y)\leq B_{2}$ also by Lemma \ref{Lem3.17}, which gives a nontrivial 
$p$-unipotent subgroup in $H$, a contradiction to Lemma \ref{corUnicityBorelSuperMax} 
or Lemma \ref{FB1capFB2torsionfree}. Hence $F\o(B_{1})$ is divisible and 
$d(B_{1})=d_{\infty}(B_{1})$. 

As the maximal definable decent torus of $F\o(B_{1})$ is central in $B_{1}$, it 
centralizes $Y$. Hence Lemma \ref{Lem3.17} gives 
$U_{(\infty,0)}(F\o(B_{1}))\leq N\o(Y)\leq B_{2}$, and hence 
$U_{(\infty,0)}(F\o(B_{1}))$ is contained in $H$. But the latter normalizes the 
former, and hence centralizes it by Fact \ref{ActionOnLessUnip} (1). 
\qed

\begin{lem}\label{Lem3.21}
{\bf \cite[Lemma 3.21]{BurdgesJalLemma}}
$X\o=U_{(\infty,r')}(F\o(B_{1}))$, and also 
$B_{1}=N\o(X\o)$. 
\end{lem}
\proof
By Fact \ref{ActionOnLessUnip}, $U_{(\infty,r')}(F\o(B_{1}))\cdot Y$ is nilpotent. 
Lemma \ref{Lem3.17} gives the inclusion 
$N\o_{U_{(\infty,r')}(F\o(B_{1}))\cdot Y}(Y)\leq H$. 
So $U_{(\infty,r')}(F\o(B_{1}))\leq Y$ by Normalizer Condition, 
Fact \ref{NormCond}. By Lemma \ref{U0r'H<FB2}, $Y\leq F(B_{2})$, 
so $U_{(\infty,r')}(F\o(B_{1})) \leq X\o$. But the converse to the latter inclusion 
holds by Theorem \ref{TheoXHomogene}. 

Our last claim follows by local$\o$ solvability$\o$ of $G$ 
and maximality of $B_{1}$. 
\qed

\begin{cor}\label{Cor3.22}
{\bf \cite[Corollary 3.22]{BurdgesJalLemma}}
$U_{(\infty,r')}(F\o(B_{1}))$ is abelian, and $F\o(B_{1})\leq C\o(X\o)$. 
\end{cor}
\proof
By Lemma \ref{Lem3.21} and Theorem \ref{ThmSylIntAb} 
$U_{(\infty,r')}(F\o(B_{1}))$ is abelian, and contained in $C\o(X\o)$. 
For an integer $r\neq r'$, $U_{(\infty,r)}(F\o(B_{1}))\leq C\o(X\o)$ by 
Fact \ref{StrucNilpGroups}. So our last claim follows. 
\qed

\begin{notation}
We let $Q$ denote a Carter subgroup of $H$. 
\end{notation}

\begin{lem}\label{Lem3.23}
{\bf \cite[Lemma 3.23]{BurdgesJalLemma}}
$U_{(\infty,r')}(Q)=U_{(\infty,r')}(Z(H))$, and this group is not trivial. 
\end{lem}
\proof
By Lemma \ref{Lem3.17}, $U_{(\infty,r')}(H/H')$ is not trivial. So 
$U_{(\infty,r')}(Q)$ is not trivial by Facts \ref{GConSolvG/FGdivAb} and 
\ref{Carter=NProj}. 

By Theorem \ref{ThmSylIntAb}, $Q$ and $Y$ are abelian. By 
Lemma \ref{U0r'H<FB2}, $U_{(\infty,r')}(Q)\leq Y$. So 
$U_{(\infty,r')}(Q)$ centralizes both $Q$ and the subgroup $H'$ of $Y$. 
So $U_{(\infty,r')}(Q)\leq Z(H)$ by Fact \ref{Carter=NProj}. 
Conversely, $Z\o(H)\leq Q$. 
\qed

\begin{theo}\label{TheoQCarterB1}
{\bf \cite[Theorem 3.24]{BurdgesJalLemma}}
$N\o(U_{(\infty,r')}(Q))\leq B_{2}$. So $N\o(Q)\leq B_{2}$, and $Q$ is a Carter 
subgroup of $B_{1}$. 
\end{theo}
\proof
We first show that $N\o(U_{(\infty,r')}(Q))\leq B_{2}$. By Lemma \ref{Lem3.17}, 
$N\o(Y)\leq B_{2}$. So we may assume that $U_{(\infty,r')}(Q)<Y$, and hence 
$H$ is not abelian by Lemma \ref{Lem3.23}. So $B_{1}$ and $B_{2}$ are the only 
Borel subgroups of $G$ containing $H$ by Proposition \ref{HnonabB1B2unique}. 
By Lemma \ref{Lem3.23}, $H\leq N\o(U_{(\infty,r')}(Q))$. 
By local$\o$ solvability$\o$ of $G$ 
the latter group is solvable. If it contains $H$ properly, then it can grow 
only in one Borel $B_{1}$ or $B_{2}$, and must agree with $H$ on the other. 
By Lemma \ref{Lem3.21}, $N\o(X\o)=B_{1}$. 
Since $Y=X\o\cdot U_{(\infty,r')}(Q)$ by 
Fact \ref{StrucSylowUSolv}, $N\o_{B_{1}}(U_{(\infty,r')}(Q))\leq N\o(Y)\leq B_{2}$ 
by Lemma \ref{Lem3.17}. So $N\o(U_{(\infty,r')}(Q))\leq B_{2}$. 

It follows that $N\o(Q)\leq N\o(U_{(\infty,r')}(Q))\leq B_{2}$, and 
$N\o_{B_{1}}(Q)\leq N\o_{H}(Q)=Q$, so that $Q$ is a Carter subgroup of $B_{1}$. 
\qed

\bigskip
We show now that $r'$ is the only unipotence degree $\geq 1$ (and in fact 
$\geq 0$ as well) appearing in both 
$F(B_{1})$ and $F(B_{2})$. 

\begin{lem}\label{lem3.25}
{\bf \cite[Lemma 3.25]{BurdgesJalLemma}}
$U_{(\infty,r)}(F\o(B_{1}))=1$ for any $r\neq r'$ with $1\leq r \leq d(B_{2})$.
\end{lem}
\proof
Let $T=U_{(\infty,r)}(F\o(B_{1}))$. We claim that $T\leq H$. First suppose 
that $d(B_{2})=r'$. Then $T\cdot Y$ is nilpotent by Fact \ref{ActionOnLessUnip}, and 
$Y$ centralizes $T$ by Fact \ref{StrucNilpGroups}. So $T\leq N\o(Y)\leq B_{2}$ 
by Lemma \ref{Lem3.17}, and $T\leq H$. Next, suppose that $d(B_{2})>r'$. Then 
$U_{(\infty,d(B_{2}))}(B_{2})\leq Z(H)$ by Lemma \ref{FactdeFB2dansH}. 
By Fact \ref{StrucNilpGroups}, $U:=T\cdot U_{(\infty,d(B_{2}))}(B_{2})$ is 
nilpotent. If $r\neq d(B_{2})$, then $T\leq C\o(U_{(\infty,d(B_{2}))}(B_{2}))$ 
by Fact \ref{StrucNilpGroups}, and $T\leq H$. So we may assume that $r=d(B_{2})$. 
If $T\nleq B_{2}$, then by Normalizer Condition, Fact \ref{NormCond}, 
$U_{(\infty,r)}(B_{2})<U_{(\infty,r)}(N\o_{U}(U_{(\infty,r)}(B_{2})))$, a 
contradiction to the fact that $B_{2}=N\o(U_{(\infty,r)}(B_{2}))$ 
by local$\o$ solvability$\o$ of $G$. Thus $T\leq H$. 

Since $T\leq H$, and $U_{(\infty,r)}(H')=1$ by Theorem \ref{TheoXHomogene}, 
$T$ is contained in a Carter subgroup of $H$ by Fact \ref{StrucSylowUSolv}. 
Now $T\leq Q$ because $T\trianglelefteq H$ and Carter subgroups are conjugate 
in $H$. Clearly $T\leq F\o(H)$ too. By Fact \ref{Carter=NProj} and 
Theorem \ref{ThmSylIntAb} $H=F\o(H)Q\leq C\o(T)$, and hence $T\leq Z(H)$. 

Now consider the case where $r>r'$. Then $U_{(\infty,r')}(F\o(B_{2}))\cdot T$ 
is nilpotent by Fact \ref{ActionOnLessUnip}, and both factors commute by 
Fact \ref{StrucNilpGroups}. If $T\neq 1$, then $B_{1}=N\o(T)$ 
by local$\o$ solvability$\o$ of $G$ 
and $U_{(\infty,r')}(F\o(B_{2}))\leq N\o(T)\leq B_{1}$, a contradiction to 
Lemma \ref{Fr'B2notinH}. Thus $T=1$. 

Finally consider the case where $r<r'$. Since $T\leq Z(H)$, $TY$ is abelian by 
Theorem \ref{ThmSylIntAb}. Recall that $Y\leq F(B_{2})$ by 
Lemma \ref{U0r'H<FB2}. Let $P$ denote the group 
$U_{(\infty,r')}(N_{F(B_{2})}(Y))$. Then 
$[x,h]\in Y$ for any $x\in X\o$ and any $h\in P$, and hence $[x,h]=[x,h]^{t}=[x,h^{t}]$ 
for any $t\in T$. So $[h^{-1},t]=hh^{-t}\in C(X\o)$. Now $[P,T]\leq Y$ by 
Lemma \ref{LemX1GenNX1} and Fact \ref{CommuttildpGp}. 
Since $P$ is nilpotent, and $T$ commutes with $Y$, the product $TP$ is nilpotent. 
By Fact \ref{StrucNilpGroups}, $P\leq N\o(T)$, which is equal to $B_{1}$ 
by local$\o$ solvability$\o$ of $G$ if 
$T\neq 1$. This contradiction to Lemma \ref{Lem3.17} shows that $T=1$. 
\qed

\bigskip
As a result, $r'$ is also uniquely determined by $B_{1}$. 

\begin{cor}\label{Cor3.26}
{\bf \cite[Corollary 3.26]{BurdgesJalLemma}}
$r'$ is the minimal unipotence degree $1\leq r<\infty$ 
such that $F(B_{1})$ admits the unipotence parameter $(\infty,r)$. 
\end{cor}

\begin{cor}
{\bf \cite[Corollary 3.27]{BurdgesJalLemma}}
For $1\leq r\leq d(B_{2})$, a Sylow $(\infty,r)$-subgroup of $H$ is a Sylow 
$(\infty,r)$-subgroup of $B_{1}$. 
\end{cor}
\proof
By Lemmas \ref{lem3.25} and \ref{Lem3.21}, $U_{(\infty,r)}(F\o(B_{1}))\leq H$. 
Since $Q$ is a Carter subgroup of $B_{1}$ by Theorem \ref{TheoQCarterB1}, 
the subgroup 
$$U_{(\infty,r)}(F\o(B_{1}))\cdot U_{(\infty,r)}(Q)$$
of $H$ is a Sylow $(\infty,r)$-subgroup of $H$ and of $B_{1}$ by 
Theorem \ref{StrucSylowUSolv}. One concludes then by conjugacy of 
Sylow $(\infty,r)$-subgroups. 
\qed

\subsubsection{Nonabelian intersections}

\begin{rem}\label{RemarkTorXToral}
$\Tor(X)$ is toral and in $Z(B_{1})\cap Z(B_{2})$, and $C\o(X)=C\o(X\o)$. 
\end{rem}
\proof
Let $S$ be the (finite) torsion subgroup of $X$, as in Lemma \ref{FB1capFB2torsionfree}. 
As $d(B_{1})<\infty$ and $d(B_{2})<\infty$, $S$ is a toral subgroup of 
$B_{1}$ and $B_{2}$, and in $Z(B_{1})\cap Z(B_{2})$ by 
Lemma \ref{FacptTorFitUpPerp}. 

By Lemma \ref{LemX1GenNX1}, $C\o(X)\leq C\o(X\o)\leq B_{1}$, and as 
$X=X\o \times S$ with $S\leq Z(B_{1})$, $C\o(X)=C\o(X\o)$. 
\qed

\begin{lem}\label{lem3.28}
{\bf \cite[Lemma 3.28]{BurdgesJalLemma}}
The subgroup $C\o(X\o)$ is not nilpotent. 
If $H$ is not abelian, then $B_{1}$ is the unique 
Borel subgroup of $G$ containing $C\o(X\o)$. 
\end{lem}
\proof
By Lemma \ref{LemX1GenNX1}, $C\o(X\o)\leq B_{1}$. By Lemma \ref{U0r'H<FB2} 
and Theorem \ref{ThmSylIntAb}, $U_{(\infty,r')}(Q)\leq C\o(X\o)$. 
By Fact \ref{StrucNilpGroups} and the fact that $d(B_{1})\neq r'$ 
(Corollary \ref{corAsymMaxPair}), 
$U_{(\infty,d(B_{1}))}(B_{1})\leq C\o(X\o)$ too. 
By Theorem \ref{TheoQCarterB1} and Fact \ref{StrucNilpGroups}, 
$U_{(\infty,d(B_{1}))}(B_{1})\cdot U_{(\infty,r')}(Q)$ is not nilpotent. 
So $C\o(X\o)$ is not nilpotent. 

Suppose now $H$ not abelian. Suppose then toward a contradiction that a Borel subgroup 
of $G$ distinct from $B_{1}$ contains $C\o(X\o)$. So there is a maximal pair 
$(B_{3},B_{4})$ which contains $C\o(X\o)$. We may assume $d(B_{3})\geq d(B_{4})$. 
Let $K=[C\o(X\o)]'$. By Corollary \ref{Cor3.22}, $F\o(B_{1})\leq C\o(X\o)$. So 
$C\o(X\o)\trianglelefteq B_{1}$ by Fact \ref{GConSolvG/FGdivAb}. Now 
$N\o(K)=B_{1}$ by local$\o$ solvability$\o$ of $G$ and maximality of $B_{1}$. 
Since $K\leq{[F(B_{3})\cap F(B_{4})]\o}$, we have by, Corollary \ref{Cor3.22} 
applied to the pair $(B_{3},B_{4})$, $F\o(B_{3})\leq C\o([F(B_{3})\cap F(B_{4})]\o)
\leq C\o(K)$. Thus $d(B_{1})\geq d(B_{3})>d(B_{4})$ by Lemma \ref{corAsymMaxPair} 
applied again to the pair $(B_{3},B_{4})$. But as $F\o(B_{1})\leq B_{4}$ also, 
$d(B_{4})\geq d(B_{1})$, a contradiction. Hence when $H$ is not abelian $B_{1}$ is 
the unique Borel subgroup of $G$ containing $C\o(X\o)$. 
\qed

\begin{cor}\label{Cor3.29}
{\bf \cite[Corollary 3.29]{BurdgesJalLemma}}
Suppose $H$ not abelian. Then, for any infinite definable subgroup $X_{1}\leq X$, 
$B_{1}$ is the unique Borel subgroup of $G$ containing $C\o(X_{1})$. 
\end{cor}
\proof
Recall that $C\o(X)=C\o(X\o)$. 
$C\o(X)\leq C\o(X_{1})$, the latter being solvable by local$\o$ solvability$\o$ 
of $G$, so the preceding lemma gives the desired result. 
\qed

\begin{cor}
If $H$ is nonabelian, then $C\o(Y)\leq C\o(X\o)\leq B_{1}$. 
\end{cor}
\proof
$X\o\leq Y$. 
\qed

\begin{lem}\label{LemYUniqB2}
{\bf (Compare with \cite[Lemma 3.10]{Deloro05})}
If $H$ is nonabelian, then any Sylow $(\infty,r')$-subgroup of $G$ containing 
$Y$ is contained in $B_{2}$. 
\end{lem}
\proof
We want to show that $\Sigma \leq B_{2}$ for any Sylow $(\infty,r')$-subgroup 
$\Sigma$ of $G$ containing $Y$. One can assume $Y < \Sigma$, and then 
$Y<U_{(\infty,r')}(N\o_{\Sigma}(Y))$ by normalizer condition, 
Fact \ref{NormCond}. By Lemma \ref{Lem3.17}, $N\o(Y)\leq B_{2}$, and thus 
$U_{(\infty,r')}(N\o_{\Sigma}(Y)) \leq B_{2}$. 

If $U_{(\infty,r')}(N\o_{\Sigma}(Y))$ is abelian, then it centralizes $Y$. 
But $C\o(Y)\leq C\o(X\o)\leq B_{1}$ by Lemma \ref{lem3.28}. 
Hence $U_{(\infty,r')}(N\o_{\Sigma}(Y))\leq (B_{1}\cap B_{2})\o=H$ 
and then $U_{(\infty,r')}(N\o_{\Sigma}(Y))=Y$, a contradiction. 

Hence $U_{(\infty,r')}(N\o_{\Sigma}(Y))$ is nonabelian. Now it follows 
from the result obtained in Theorem \ref{ThmSylIntAb} that in a locally$\o$ 
solvable$\o$ group of finite Morley rank, a nonabelian definable connected 
nilpotent subgroup is contained in a unique Borel subgroup. 
As $U_{(\infty,r')}(N\o_{\Sigma}(Y))$ is in $B_{2}$ and in $\Sigma$, 
this gives $\Sigma \leq B_{2}$. 
\qed

\begin{lem}\label{lem3.30}
{\bf \cite[Lemma 3.30]{BurdgesJalLemma}}
Let $B$ be a Borel subgroup of $G$, distinct from $B_{1}$. Suppose that 
$(B,B_{1})$ is a maximal pair, that $H_{1}=(B\cap B_{1})\o$ is not abelian, and that 
$d(B_{1})\geq d(B)$. Then $B$ is $F\o(B_{1})$-conjugate to $B_{2}$. 
\end{lem}
\proof
We can apply the results of the above analysis to the maximal pair $(B_{1},B)$. 
We observe that $H_{1}'\leq F(B_{1})\cap F(B)$. By Corollary \ref{Cor3.26} 
and Theorem \ref{TheoXHomogene}, $r'=d(H_{1}')$, and 
both $H'$ and $H_{1}'$ are contained in $U_{(\infty,r')}(F\o(B_{1}))$. 
By Lemma \ref{Lem3.21}, $U_{(\infty,r')}(F\o(B_{1}))$ is contained in 
both $H$ and $H_{1}$. Let $Q$ and $Q_{1}$ be Carter subgroups of $H$ and $H_{1}$ 
respectively. By Theorem \ref{TheoQCarterB1}, $Q$ and $Q_{1}$ are Carter subgroups 
of $B_{1}$. By conjugacy of Carter subgroups in connected solvable groups, 
$Q_{1}=Q^{h}$ for some $h\in B_{1}$, and we may assume $h\in F\o(B_{1})$ 
by Facts \ref{GConSolvG/FGdivAb} and \ref{Carter=NProj}. 
By Facts \ref{GConSolvG/FGdivAb} and \ref{Carter=NProj}, $Q$ and $Q_{1}$ 
cover $H/H'$ and $H_{1}/H_{1}'$ respectively. By Lemma \ref{Lem3.21}, 
$$H^{h}=U_{(\infty,r')}(F\o(B_{1}))\cdot Q^{h}=U_{(\infty,r')}(F\o(B_{1}))\cdot Q_{1}=H_{1}.$$
Since $H_{1}$ is not abelian, $B_{2}^{h}=B$ by Proposition \ref{HnonabB1B2unique}. 
\qed

\subsubsection{Conclusions}

\begin{prop}\label{prop4.1OnInter}
{\bf \cite[Proposition 4.1]{BurdgesJalLemma}}
Let $G$ be a locally$\o$ solvable$\o$ group of finite Morley rank, 
$B_{1}$ and $B_{2}$ two distinct Borel subgroups of $G$, and 
$H$ a nontrivial definable connected subgroup of $B_{1}\cap B_{2}$. 
Then the following hold: 
\begin{itemize}
\item[$(1)$]
$H'$ is a homogeneous $(\infty,r')$-group for some $1\leq r'<\infty$ (or trivial). 
\item[$(2)$]
Every definable connected nilpotent subgroup of $H$ is abelian. 
\item[$(3)$]
$U_{(\infty,r')}(F\o(H))=U_{(\infty,r')}(H)$ is the unique 
Sylow $(\infty,r')$-subgroup of $H$. 
\item[$(4)$]
$U_{\tilde q}(F\o(H))\leq Z(H)$ for any $\tilde q \neq (\infty,r')$. 
\item[$(5)$]
$0\leq d_{\infty}(H)=d(H) \leq d(C(H')) \leq d(N(H'))\leq \infty$, all 
inequalities, except maybe the third one, being strict when $H$ is not abelian. 
\end{itemize}
\end{prop}
\proof
We may assume $H$ not abelian, as otherwise all statements are trivially true once 
one has noticed that $d_{\infty}(H)=d(H)$ 
by Corollary \ref{corUnicityBorelSuperMax}. 

Let $(B_{3},B_{4})$ be a maximal pair containing $H$, with $d(B_{3})\geq d(B_{4})$. 
The first two conclusions follow immediately from 
Theorems \ref{TheoXHomogene} and \ref{ThmSylIntAb}. The third conclusion 
follows from Lemma \ref{U0r'H<FB2}. For the fourth conclusion, if 
$\tilde q=(\infty,r)$ with $r\neq r'$, then 
$U_{\tilde q}(F\o(H))$ lies in a Carter subgroup $Q$ of $H$ by 
Fact \ref{StrucSylowUSolv}, and $H\leq QH'$ 
(Facts \ref{GConSolvG/FGdivAb} and \ref{Carter=NProj}) 
$\leq C\o(U_{\tilde q}(F\o(H)))$ 
(by the second point), which shows the fourth point. 

By Corollary \ref{corAsymMaxPair}, 
$\infty > d(B_{3})>d(B_{4})\geq d(H)=d_{\infty}(H)>0$ 
(be careful, this is not the same $H$, and one uses also the divisibility of $F\o(B_{3})$ 
and of $F\o(B_{4})$). 
By Fact \ref{StrucNilpGroups}, $U_{(\infty,d_{\infty}(B_{3}))}(B_{3})\leq C(H')$, 
thus $d_{\infty}(C(H'))>d_{\infty}(H)=d(H)$. Hence for the last point it suffices 
to show that $d(N(H'))<\infty$. Otherwise, 
$U_{p}(N(H'))$ is nontrivial for some prime $p$; 
now the nontrivial group $U_{(\infty,d_{\infty}(B_{3}))}(B_{3})$, 
which is also in $N(H')$, normalizes $U_{p}(N(H'))$, and centralizes it by 
Fact \ref{ActionOnLessUnip} $(4)$, so that 
$U_{p}(N(H')) \leq N\o(U_{(\infty,d_{\infty}(B_{3}))}(B_{3}))=B_{3}$ 
(by local$\o$ solvability$\o$ and Lemma \ref{CharactBorelXLocal}), 
a contradiction to the divisibility of $F\o(B_{3})$. Hence 
$d(N(H'))<\infty$ and this completes the proof of the fifth point. 
\qed

\begin{cor}\label{LemUniqSurNilpNonab}
{\bf \cite[Corollary 4.2]{BurdgesJalLemma}}
Let $G$ be a locally$\o$ solvable$\o$ group of finite Morley rank. 
Then a definable connected 
nonabelian nilpotent subgroup is contained in exactly one Borel subgroup of $G$. 
\end{cor}

\begin{cor}
{\bf \cite[Corollary 4.2']{BurdgesJalLemma}}
Let $G$ be a locally$\o$ solvable$\o$ group of finite Morley rank. 
If $Q$ is a Carter subgroup of a Borel subgroup $B$, and if $Q$ is not 
abelian, then $Q$ is a Carter subgroup of $G$. 
\end{cor}
\proof
$N\o(Q)$ is contained in a Borel subgroup $B_{1}$ of $G$ by local$\o$ solvability$\o$. 
As $Q\leq B\cap B_{1}$, $B=B_{1}$ by Corollary \ref{LemUniqSurNilpNonab}, 
and $N\o_{G}(Q)\leq N\o_{B_{1}}(Q)= N\o_{B}(Q)=Q$
\qed

\begin{lem}\label{lem4.4}
{\bf \cite[Lemma 4.4]{BurdgesJalLemma}}
Let $G$ be a locally$\o$ solvable$\o$ group of finite Morley rank, $B_{1}$ and $B_{2}$ 
two distinct Borel subgroups of $G$. Suppose that $H=(B_{1}\cap B_{2})\o$ is not abelian, 
and that $C\o(H')\leq B_{1}$. Then $B_{1}$ and $B_{2}$ are the only Borel subgroups 
containing $H$. 
\end{lem}
\proof
Suppose toward a contradiction that $G$ contains a Borel subgroup $B$ distinct from 
both $B_{1}$ and $B_{2}$ and which contains $H$. We may choose $B$ such that 
$H_{2}=(B\cap B_{2})\o$ is maximal subject to $B\neq B_{1}$ and $H\leq B$, $B_{2}$. 
Consider a maximal pair $(B_{3},B_{4})$ containing $H_{2}$ and such that 
$d_{\infty}(B_{3})\geq d_{\infty}(B_{4})$. Corollary \ref{Cor3.29} 
applied to $(B_{3},B_{4})$ implies that $B_{3}$ 
is the unique Borel subgroup of $G$ containing $C\o(H')$. So $B_{1}=B_{3}$. Thus 
$H=H_{2}$. By Proposition \ref{HnonabB1B2unique}, $B_{1}=B_{3}$ and $B_{4}$ 
are the unique Borel subgroups containing the connected component of their intersection. 
So we may assume $B_{4}\neq B_{2}$, as otherwise we are done. 
Therefore we may also assume that $B=B_{4}$. So $H_{1}=(B_{1}\cap B)\o=(B_{3}\cap B_{4})\o$ 
corresponds to the intersection$\o$ of maximal pairs, and we can apply the previous results 
to this intersection. We observe that 
$$r':=d_{\infty}(H')=d_{\infty}(F(B_{1})\cap F(B))$$
by Theorem \ref{TheoXHomogene}. 

Consider first the case $U_{(\infty,r')}(F\o(B_{2}))\leq B_{1}$. Since $H'$ is 
$(\infty,r')$-homoge\-neous, 
$U_{\tilde q}(F\o(B_{2}))\leq C\o(H')\leq B_{1}$ for every 
$\tilde q \neq (\infty,r')$. Hence $F\o(B_{2})\leq H$, 
and $H\trianglelefteq B_{2}$. 
By local$\o$ solvability$\o$ of $G$ and maximality of $B_{2}$, $N\o(H')=B_{2}$. But 
Corollary \ref{Cor3.22} applied to $(B,B_{1})$ yields 
$F\o(B_{1})\leq C\o([F(B)\cap F(B_{1})]\o)\leq C\o(H')\leq N\o(H')=B_{2}$. 
Then $U_{(\infty,d_{\infty}(B_{1}))}(B_{1})\leq H$ by Fact \ref{FactUnipHeaviest}, 
and $d_{\infty}(B_{1})\leq d_{\infty}(H_{1})$. This contradicts Lemma \ref{HypdB2leqdB1} 
applied with $(B,B_{1})$. 

Consider next the case $U_{(\infty,r')}(F\o(B_{2}))\nleq B_{1}$. Let 
$P=U_{(\infty,r')}(H)$, the unique Sylow $(\infty,r')$-subgroup of $H$ by 
Proposition \ref{prop4.1OnInter} (3), and $M=N\o(P)$, a solvable group 
by local$\o$ solvability$\o$. Since $P$ normalizes 
$U_{(\infty,r')}(F\o(B_{2}))$, 
$U_{(\infty,r')}(F\o(B_{2}))\cdot P$ is nilpotent 
by Fact \ref{ActionOnLessUnip}. 
By Normalizer Condition, Fact \ref{NormCond}, 
$P<{U_{(\infty,r')}(N_{U_{(\infty,r')}(F\o(B_{2}))\cdot P}(P))}$. 
Since $H=(B_{1}\cap B_{2})\o$ and $P=U_{(\infty,r')}(H)$, 
it follows that $H\leq {(M\cap B_{2})\o}\nleq B_{1}$ and that 
$H<{(M\cap B_{2})\o}$. Hence $M\leq B_{2}$ by 
maximality of $H_{2}$ ($=H$). 
By Lemma \ref{Fr'B2notinH}, $U_{(\infty,r')}(F\o(B))\nleq H_{1}$. Since $P$ 
normalizes $U_{(\infty,r')}(F\o(B))$, $U_{(\infty,r')}(F\o(B)) \cdot P$ is nilpotent 
by Fact \ref{ActionOnLessUnip}. By Normalizer Condition, 
Fact \ref{NormCond}, and using $M\leq B_{2}$, 
$P<{U_{(\infty,r')}(N\o_{U_{(\infty,r')}(F\o(B)) \cdot P}(P))}
\leq{(B\cap B_{2})\o}=H_{2}=H$, a contradiction to 
$P=U_{(\infty,r')}(H)$. 
\qed

\bigskip
We can now characterize maximal pairs with nonabelian intersections. 

\begin{theo}\label{TheoEquivMaxPair}
{\bf \cite[Theorem 4.3]{BurdgesJalLemma}}
Let $G$ be a locally$\o$ solvable$\o$ group of finite Morley rank, 
$B_{1}$ and $B_{2}$ two distinct Borel subgroups of $G$. 
Suppose $H=(B_{1}\cap B_{2})\o$ nonabelian. Then the following are equivalent: 
\begin{itemize}
\item[$(1)$]
$B_{1}$ and $B_{2}$ are the only Borel subgroups containing $H$. 
\item[$(2)$]
$(B_{1},B_{2})$ is a maximal pair. 
\item[$(3)$]
If $B_{3}\neq B_{1}$ is a Borel subgroup containing $H$, then $(B_{1}\cap B_{3})\o=H$. 
\item[$(4)$]
$C\o(H')$ is contained in $B_{1}$ or $B_{2}$. 
\item[$(5)$]
$B_{1}$ and $B_{2}$ are not conjugate under the action of $C\o(H')$. 
\item[$(6)$]
$d_{\infty}(B_{1})\neq d_{\infty}(B_{2})$. 
\end{itemize}
\end{theo}

\proof
Clearly $(1)$ implies $(2)$, $(2)$ implies $(3)$, and $(4)$ implies $(5)$. 
By Lemmas \ref{HypdB2leqdB1} and \ref{LemdH=dB2}, 
$(2)$ implies $(6)$. Clearly $(6)$ implies $(5)$. 
By local$\o$ solvability$\o$ of $G$, there exists a Borel subgroup $B_{c}$ of $G$ 
containing $N\o(H')$. 

We show now that $(3)$ implies $(4)$. Let $B_{x}$ denotes $B_{1}$, unless 
$B_{c}=B_{1}$, in which case we let $B_{x}$ denote $B_{2}$. By 
$(3)$, $H=(B_{c}\cap B_{x})\o$. By Lemma \ref{lem4.4} 
applied to the pair $(B_{c},B_{x})$, $B_{c}\geq C\o(H')$ must be one of 
$B_{1}$ or $B_{2}$, so $(4)$ holds. 

We show now that $(5)$ implies $(1)$. Assume $(1)$ fails. Then, for $i=1$ and $2$, 
$C\o(H')\nleq B_{i}$ by Lemma \ref{lem4.4}. But $(B_{c},B_{1})$ and 
$(B_{c},B_{2})$ are maximal pairs, by Lemma \ref{lem4.4} again. 
So $d(B_{c})\geq d(B_{1})$, $d(B_{2})$, by Lemma \ref{lem3.28}. 
By Lemma \ref{lem3.30}, $B_{1}$ is $F\o(B_{c})$-conjugate to $B_{2}$. 
By Corollary \ref{Cor3.22}, $F\o(B_{c})\leq C\o(H')$, so $(5)$ fails. 
\qed

\bigskip
We can now describe the maximal pairs having a nonabelian 
intersection$\o$, collecting the results from \cite{BurdgesJalLemma} 
with the additional results from \cite{Deloro05}. 
We slightly change the presentation in comparison to 
\cite[Theorem 4.5]{BurdgesJalLemma}, as we prefer to distinguish 
between a symmetric version and an asymmetric one. 
We start with the symmetric version. 

\begin{theo}\label{TheoPairMaxNAbAsym}
Let $G$ be a locally$\o$ solvable$\o$ group of finite Morley rank, and $(B_{1},B_{2})$ 
a maximal pair of Borel subgroups such that $H:=(B_{1}\cap B_{2})\o$ is 
nonabelian. Let $r'=d_{\infty}(H')$. 
\begin{itemize}
\item[$(1)$]
$0<d(B_{1})<\infty$ and $0<d(B_{2})<\infty$. 
\item[$(2)$]
$N\o(H)=H$. 
\item[$(3)$]
$[F(B_{1})\cap F(B_{2})]\o$ is $(\infty,r')$-homogeneous, and $r'>0$. 
\end{itemize}
Furthermore, if $Q$ denotes a Carter subgroup of $H$ and $Q_{r'}$ denotes 
$U_{(\infty,r')}(Q)$, then 
\begin{itemize}
\item[$(4)$]
$Q_{r'}\neq 1$, 
\end{itemize}
and exactly one of the following cases occur: 
\begin{itemize}
\item[$(4.a)$]
$N\o(Q_{r'})=H$. 
\item[$(4.b)$]
$H<N\o_{B_{1}}(Q_{r'})$; furthermore $N\o_{B_{2}}(Q_{r'})=H$ and $B_{1}$ 
is the unique Borel subgroup containing $N\o(Q_{r'})$. 
\item[$(4.c)$]
$H<N\o_{B_{2}}(Q_{r'})$; furthermore $N\o_{B_{1}}(Q_{r'})=H$ and $B_{2}$ 
is the unique Borel subgroup containing $N\o(Q_{r'})$. 
\end{itemize}
\end{theo}
\proof

$(1)$: \ref{FB1capFB2torsionfree}, \ref{LemFB1DivBender}, \ref{FactdeFB2dansH}. 

$(2)$: \ref{BenderHNonAbNH=H}. 

$(3)$: \ref{FB1capFB2torsionfree}, \ref{TheoXHomogene}. 

$(4)$: \ref{Lem3.23}, and proof of \cite[Lemme 3.9]{Deloro05} for the trichotomy. 
\qed

\bigskip
We finish with the description once the asymmetry is fixed. 

\begin{theo}\label{TheoPairMaxNAbSym}
Assume in addition to Theorem \ref{TheoPairMaxNAbAsym} 
that $d(B_{1})\geq d(B_{2})$. Then
\begin{itemize}
\item[$(1)$]
$0<d(B_{2})<d(H)=d(B_{1})<\infty$. 
\item[$(2)$]
$Q$ is a Carter subgroup of $B_{1}$. 
\item[$(3)$]
$U_{(\infty,r')}(F(B_{1}))=[F(B_{1})\cap F(B_{2})]\o$. 
\item[$(4)$]
$B_{1}$ is the unique Borel subgroup containing $C\o(U_{(\infty,r')}(F(B_{1})))$.
\item[$(5)$]
$N\o(Q)\leq B_{2}$.
\item[$(6)$]
$U_{(\infty,r')}(H)\leq F\o(B_{2})$, and $N\o(U_{(\infty,r')}(H))\leq B_{2}$.
\item[$(7)$]
$U_{\tilde q}(F(B_{2}))\leq Z(H)$ for any $\tilde q\neq (\infty,r')$, 
and $U_{(\infty,r')}(F(B_{2}))$ is nonabelian (in particular $U_{\tilde q}(F(B_{2}))$ 
is nonabelian iff $\tilde q=(\infty,r')$). 
\item[$(8)$]
Any Sylow $(\infty,r')$-subgroup of $G$ containing $U_{(\infty,r')}(H)$ is contained 
in $B_{2}$. 
\item[$(9)$]
If $Q$ is a Carter subgroup of $B_{2}$, then $U_{(\infty,r')}(F(B_{2}))$ 
is the unique Sylow $(\infty,r')$-subgroup of $B_{2}$, and 
in particular the unique Sylow $(\infty,r')$-subgroup of $G$ 
containing $U_{(\infty,r')}(H)$. 
\end{itemize}
\end{theo}
\proof

$(1)$: \ref{TheoPairMaxNAbAsym} $(1)$, \ref{corAsymMaxPair}. 

$(2)$: \ref{TheoQCarterB1}. 

$(3)$: \ref{Lem3.21}. 

$(4)$: \ref{Lem3.21}, \ref{lem3.28}. 

$(5)$: \ref{TheoQCarterB1}

$(6)$: \ref{U0r'H<FB2}, \ref{Lem3.17}. 

$(7)$: \ref{FactdeFB2dansH}, \ref{Fr'B2notinH}. 

$(8)$: \ref{LemYUniqB2}. 

$(9)$: \ref{CarterB1B2UB2Nilp}, \ref{LemYUniqB2}.
\qed

\bigskip
Finally, we record a point about exceptional elements, which applies in 
particular in 
Theorems \ref{TheoPairMaxNAbAsym} and \ref{TheoPairMaxNAbSym}. 

\begin{theo}
Let $G$ be a locally$\o$ solvable$\o$ group of finite Morley rank and 
$(B_{1},B_{2})$ 
a maximal pair of Borel subgroups such that $[F(B_{1})\cap F(B_{2})]\o$ is 
nontrivial. Then the finite subgroup $S$ of $F(B_{1})\cap F(B_{2})$ as 
in Lemma \ref{FB1capFB2torsionfree} is toral and central, 
both in $B_{1}$ and $B_{2}$. 
\end{theo}
\proof
$F\o(B_{1})$ and $F\o(B_{2})$ are divisible by 
Lemmas \ref{FactdeFB2dansH} and \ref{LemFB1DivBender}, and 
Remark \ref{RemarkTorXToral} applies. 
\qed

\subsection{An extra homogeneity result}

The following extra homogeneity result proved for the purpose of 
\cite{Deloro05} is essentially a corollary 
of Corollary \ref{LemUniqSurNilpNonab}. 

\begin{lem}\label{LemXHomFBNonAb}
{\bf (Compare with \cite[Lemme 3.8]{Deloro05})}
Let $G$ be a locally$\o$ solvable$\o$ group of finite Morley rank, $B$ and 
$B^{g}$ two distinct conjugates of a same Borel subgroup $B$. If 
$[F(B)\cap F(B^{g})]\o$ is not homogeneous, then $F\o(B)$ is abelian. 
\end{lem}
\proof
By assumption $[F(B)\cap F(B^{g})]\o$ contains two nontrivial Sylow subgroups 
$U_{1}$ and $U_{2}$ with two distinct unipotence parameters, 
say $\tilde p$ for $U_{1}$ and $\tilde q$ for $U_{2}$. 

By local$\o$ solvability$\o$ of $G$, $N\o(U_{1})$ is contained in a Borel subgroup 
$B_{1}$. If $B_{1}\neq B$, then by Fact \ref{StrucNilpGroups} $(2)$ 
Corollary \ref{LemUniqSurNilpNonab} 
implies that Sylow subgroups of $F\o(B)$ of unipotence parameters different from 
$\tilde p$ are abelian. If $B_{1}=B$, then $B_{1}\neq B^{g}$ and one sees 
similarly that Sylow subgroups of unipotence parameters different from $\tilde p$ 
of $F\o(B^{g})$, and thus also of $F\o(B)$, are also abelian. 

Considering a Borel subgroup $B_{2}$ containing $N\o(U_{2})$, one sees similarly 
that Sylow subgroups of $F\o(B)$ of unipotence parameters different from 
$\tilde q$ are abelian. 

Now $F\o(B)$ is abelian by Fact \ref{StrucNilpGroups} $(2)$. 
\qed

\subsection{Exceptional connected subgroups}

Section \ref{SectMaxPaire} concerned the analysis of intersections of 
maximal pairs of Borel subgroups. In the present section we continue 
a little bit in this vein when one of the two subgroups involved is not necessarily 
solvable, a possibility in the context of locally$\o$ solvable$\o$ 
groups of finite Morley rank in comparison to the context of minimal 
connected simple groups. 

\begin{defi}\label{DefiMaximalInterAsym}
Let $G$ be a group of finite Morley rank and $K$ a definable connected subgroup 
of $G$. We say that a Borel subgroup $B$ of $G$ has {\em maximal intersection} with 
$K$ if $B\nleq K$ and $(K\cap B)\o$ is maximal for inclusion among groups 
of the form $(K\cap B_{1})\o$ for some Borel subgroup $B_{1}$ of $G$ such that 
$B_{1}\nleq K$. 
\end{defi}

We note in Definition \ref{DefiMaximalInterAsym} that if $K$ is solvable and 
not a Borel subgroup, then it has a maximal intersection with any Borel subgroup 
containing it. If $K$ is a Borel subgroup of $G$, and a Borel subgroup $B$ has a maximal 
intersection with $K$, then $G\o$ is not solvable. 

\begin{lem}\label{LemKContHB1Prop}
Let $G$ be a group of finite Morley rank, $K$ a definable connected subgroup of $G$, 
and $B$ a Borel subgroup of $G$ having maximal intersection with $K$. Then any Borel 
subgroup $B_{1}$ of $G$ such that $(K\cap B)\o < (K\cap B_{1})\o$ is in $K$. 
\end{lem}
\proof 
This is immediate by definition. 
\qed

\bigskip
It follows that if $K$ is a Borel subgroup of a locally$\o$ solvable$\o$ group $G$ 
and $B$ is a Borel subgroup of $G$ having maximal intersection with $K$, 
then if $(K\cap B)\o$ is nonabelian any Borel 
subgroup $B_{3}$ of $G$ containing $(K\cap B)\o$ such that 
$(K\cap B)\o<(K\cap B_{3})\o$ must be $K$, and hence $(K,B)$ is a maximal pair of 
Borel subgroups of $G$ by the equivalence provided 
in Theorem \ref{TheoEquivMaxPair} $(3)$. 

In the general case of a locally$\o$ solvable$\o$ group $G$ a proper 
definable connected subgroup $K$ can be nonsolvable, and we slightly 
clarify the situation in this general case. 

\begin{lem}
Let $G$ be a locally$\o$ solvable$\o$ group of finite Morley rank, 
$K$ a nontrivial definable connected subgroup of $G$, $B$ a Borel subgroup 
of $G$ having maximal intersection with $K$, and let $H=(K\cap B)\o$. 
Then assuming $H$ nontrivial exactly one of the following cases occurs. 
\begin{itemize}
\item[$(1)$]
$H$ is an abelian Carter subgroup of $K$ and of $B$. 
\item[$(2)$]
$H$ is an abelian Carter subgroup of $K$ and $H<N\o_{B}(H)\leq B$. 
\item[$(3)$]
$H$ is an abelian Carter subgroup of $B$, and $H<N\o(H)\leq K$. In this 
case any Borel subgroup of $K$ containing $H$ is a Borel subgroup of $G$. 
\item[$(4)$]
$H$ is a nonabelian Borel subgroup of $K$. 
\item[$(5)$]
$H$ is nonabelian and not a Borel subgroup of $K$. In this 
case any Borel subgroup of $K$ containing $H$ is a Borel subgroup of $G$. 
\end{itemize}
\end{lem}
\proof
Notice that $N\o(H)$ is solvable by local$\o$ solvability$\o$ of $G$. 

Assume first $H$ abelian. If $H$ has finite index in its normalizer in $K$ and in $B$ 
then we are in case $(1)$. 

Assume now $H<N\o_{B}(H)$. Then the maximality 
of the intersection forces $N\o_{K}(H)=H$, and $H$ is an abelian Carter 
subgroup of $K$. Hence we are in case $(2)$. 

Assume now $H<N\o_{K}(H)$. Then the maximality of the intersection 
forces $N\o(H)\leq K$ with Lemma \ref{LemKContHB1Prop}. 
Now $N\o_{B}(H)\leq (K\cap B)\o=H$, and $H$ is an abelian Carter subgroup of $B$. 
Hence we are in case $(3)$ by Lemma \ref{LemKContHB1Prop}. 

This treats all cases corresponding to the case $H$ abelian, so we may now assume 
$H$ nonabelian. If $H$ is a Borel subgroup of $K$, then we are in case $(4)$. 

It remains only to consider the case in which $H$ is not abelian and not a Borel subgroup 
of $K$. By Lemma \ref{LemKContHB1Prop}, any Borel subgroup of $K$ 
containing $H$ is a Borel subgroup of $G$. We are in case $(5)$. 
\qed

\section{Homogeneous cases and torsion}\label{SectHomCasesTorsion}

In this final section we collect various additional results of specialized nature 
about locally$\o$ solvable$\o$ groups of finite Morley rank, generally pending on the 
uniqueness theorems of Section \ref{SectionJaligLemma}. 

The first type of results concerns the homogeneous cases. Recall 
from \cite{FreconJaligot07} or Section \ref{SectionAbstractUnipTh} 
that a group of finite Morley rank is homogeneous if is $\tilde p$-homo\-geneous 
for some unipotence parameter $\tilde p$, that is every definable connected 
nilpotent subgroup is a $\tilde p$-group. (This is weaker than the 
definition in \cite{FreconUnipotence}, which requires to consider all 
definable connected subgroups, not only the nilpotent ones.) In a 
$\tilde p$-homogeneous group one sees easily with 
Lemma \ref{LemGenSmallMaxUnipParam} and 
Fact \ref{FactUnipHeaviest} that any Borel subgroup is a 
(homogeneous) $\tilde p$-group, and in particular nilpotent. Hence we will 
more generally consider the case in which all Borel subgroups are nilpotent, and 
look at the homogeneous cases at various levels of generality. 

The torsion-free case will be fairly well understood in this context, and 
with torsion this connects to a bit of Sylow theory. 
As far as torsion is concerned, there is in general no Sylow theory as 
in Fact \ref{StructpSylSolvGps} available in an arbitrary group of finite Morley rank. 
The following fact shows however similarities with Fact \ref{StructpSylSolvGps} 
in the general case. 

\begin{fait}\label{FactToricvs.Unip}
{\bf \cite[Theorem 3 and Corollary 3.1]{BurdgesCherlinSemisimpleTorsion}} 
Let $G$ be a connected group of finite Morley rank, $t$ a 
$\pi$-element of $G$ for some set $\pi$ of primes $p$. 
If $U_{p}(C(t))=1$ for every $p$ in $\pi$, then $t$ belongs to a, and in fact to any, 
maximal $\pi$-torus of $G$ and of $C\o(t)$.  
\end{fait}

Notice that the second statement is a mere corollary of the first, together 
with the fact that toral elements belong to the connected components of their 
centralizers and Fact \ref{ConjdecentTori}. 

\subsection{Nilpotent Borel subgroups}

In this section we consider locally$\o$ solvable$\o$ groups in which each Borel 
subgroup is nilpotent. We start with a lemma concerning abelian Borel subgroups. 

\begin{lem}\label{LemAbelianBorel}
Let $G$ be a locally$\o$ solvable$\o$ group of finite Morley rank with an 
abelian Borel subgroup $B$. Let $B_{u}$ denote the maximal definable connected 
subgroup of $B$ of bounded exponent. Then 
$B$ has a subgroup $E$ finite modulo $B_{u}$ such that ${B\cap B^{g}} \leq E$ 
for any element $g$ of $G$ not in $N(B)$, 
and one of the following two cases occurs. 
\begin{itemize}
\item[$(1)$]
$B$ is a generous abelian Carter subgroup. 
\item[$(2)$]
$B$ is an abelian Carter subgroup of bounded exponent. 
\end{itemize}
\end{lem}
\proof
For any $g$ in ${G\setminus N(B)}$, $N\o(B\cap B^{g})$ contains 
$B$ and $B^{g}$, and if $B\cap B^{g}$ is infinite then $N\o(B\cap B^{g})$ 
is solvable by local$\o$ solvability$\o$ of $G$ and one gets 
$B$, $B^{g}\leq N\o(B\cap B^{g})$ and $B=B^{g}$ by maximality, 
a contradiction. 
Hence $B\cap B^{g}$ is finite for every $g\in{G\setminus N(B)}$. 

The uniformly definable family of finite subgroups $B\cap B^{g}$, 
for $g\in{B\setminus N(B)}$, consists of subgroups of uniformly bounded 
cardinals by elimination of infinite quantifiers. As Pr\"ufer $p$-ranks are finite 
for any prime $p$, all these subgroups must be contained modulo 
$B_{u}$ in a finite subgroup of the maximal definable decent 
torus of $B$. Calling $E$ the preimage in $B$ of this group, 
this proves our first statement. 

If $B_{u}<B$, then $E$ is not generic in $B$ and one can conclude that the 
Carter subgroup $B$ of $G$ is generous by the equivalence given in 
Fact \ref{CarterGeneiffgendis} $(3)$. This proves our alternative. 
\qed

\bigskip
We note that the two cases in Lemma \ref{LemAbelianBorel} are 
a priori not necessarily mutually exclusive. In the 
locally$\o$ solvable context $E$ is necessarily trivial, and $B$ is then necessarily 
generous in any case. 

We pass now to nilpotent Borel subgroups, replacing the commutativity assumption 
by a nilpotence assumption on all Borel subgroups of the ambient group. 
The first lemma is essentially the content of the first part of the proof of 
Lemma \ref{LemAbelianBorel} and typical of earlier work on bad groups 
\cite[Chapter 13]{BorovikNesin(Book)94}. 

\begin{lem}\label{LemAllBorelNilpFinitInter}
Let $G$ be a locally$\o$ solvable$\o$ group of finite Morley 
rank in which all Borel subgroups are nilpotent. Then any two distinct 
Borel subgroups have a finite intersection. 
\end{lem}
\proof
Assume the contrary. Let $B_{1}$ and $B_{2}$ contradict our claim, 
with $[B_{1}\cap B_{2}]\o$ of maximal rank. Call the latter group $H$, and 
notice that $H<B_{1}$ and $H<B_{2}$. 
By normalizer condition in nilpotent groups, 
\cite[Lemma 6.3]{BorovikNesin(Book)94}, 
$H<N\o_{B_{1}}(H)$ and $H<N\o_{B_{2}}(H)$. 
Now $N\o(H)$ is solvable by local$\o$ solvability$\o$ of $G$, and contained 
in a Borel subgroup $B_{3}$. As $H<(B_{1}\cap B_{3})\o$ and 
$H<(B_{2}\cap B_{3})\o$, our maximality assumption forces 
$B_{1}=B_{3}=B_{2}$, a contradiction. 
\qed

\bigskip
We get in any case conclusions similar to those of Lemma \ref{LemAbelianBorel}. 

\begin{lem}\label{LemAllBorelNilp}
Let $G$ be a locally$\o$ solvable$\o$ group of finite Morley 
rank in which all Borel subgroups are nilpotent. Let $B$ be a Borel subgroup of $G$. 
Then $B$ has a definable subgroup $E$, finite modulo the bounded 
exponent part of $B$, such 
that $B\cap B^{g} \leq E$ for any intersection $B\cap B^{g}$ 
with $g\in{G\setminus N(B)}$. Moreover one of the following two cases occurs. 
\begin{itemize}
\item[$(1)$]
$B$ is a generous Carter subgroup. 
\item[$(2)$]
$B$ is a Carter subgroup of bounded exponent. 
\end{itemize}
\end{lem}
\proof
With Lemma \ref{LemAllBorelNilpFinitInter} applied to distinct 
conjugates of $B$, the existence of $E$ follows as in the proof of 
Lemma \ref{LemAbelianBorel}. The alternative proposed follows similarly as well. 
\qed

\bigskip
As for Lemma \ref{LemAbelianBorel}, the two cases in 
Lemma \ref{LemAllBorelNilp} are a priori not mutually exclusive, and if 
the ambient group $G$ is locally$\o$ solvable then distinct conjugates 
of $B$ are necessarily pairwise disjoint by the same proof 
as in Lemma \ref{LemAllBorelNilpFinitInter}, and $B$ is always 
generous. 

\subsection{The torsion-free homogeneous case}\label{SectionTorsionFreeCase} 

We shall now evacuate, or rather collect in Pandora's box of bad groups, 
$\tilde p$-homogeneous locally$\o$ solvable$\o$ groups of finite Morley rank, 
with $\tilde p$ not of the form $(\infty,0)$ or $(p,\infty)$ for $p$ a prime. 
In this case Borel subgroups are nilpotent and torsion-free by 
Facts \ref{FactCharacttildepHom} and \ref{pneqq}. More generally, we have 
the following result for such groups. 

\begin{theo}\label{TheoTorsionFreeHomegeneous}
Let $G$ be a torsion-free locally$\o$ solvable$\o$ group of finite Morley rank in which 
Borel subgroups are all nilpotent. Then Borel subgroups are conjugate and either 
\begin{itemize}
\item[$(1)$]
$G$ is nilpotent, or 
\item[$(2)$]
$B<G$ is a full Frobenius group for some Borel subgroup $B$ of $G$. 
\end{itemize}
\end{theo}

As far a torsion is concerned there is a classical lifting result. 

\begin{fait}\label{FactLiftingTorsion}
{\bf \cite{BorovikNesin92}}
Let $G$ be a group of finite Morley rank, $H$ a definable normal subgroup, and 
$x$ a $p$-element modulo $H$, for some prime $p$. Then the definable 
hull $H(x)$ of $x$ contains a $p$-element. 
\end{fait}

\proof
Notice that $G$ is connected by absence of torsion and Fact \ref{FactLiftingTorsion}. 

By Lemmas \ref{LemAllBorelNilpFinitInter} and \ref{LemAllBorelNilp}, 
distinct Borel subgroups have trivial intersections, and each Borel subgroup is generous. 
As $G$ is connected it cannot have two disjoint generic subsets. If $B_{1}$ and 
$B_{2}$ are two Borel subgroups, then two conjugates of $B_{1}$ and $B_{2}$ 
must have a nontrivial intersection by generosity, and then are equal. 
This shows that Borel subgroups are conjugate. 

If $G$ is not nilpotent, then $B<G$ for some Borel subgroup $B$ of $G$. 
By Fact \ref{FactLiftingTorsion}, $N(B)=B$, 
and $B$ is malnormal in $G$ by disjointness 
of distinct Borel subgroups. As $B^{G}$ is generic, any element $g$ of 
$G$ has an infinite centralizer 
(this is also an easy consequence of the main result of \cite{BorovikBurdgesCherlin07} 
in arbitrary connected groups), 
and in particular normalizes a Borel subgroup by 
Lemma \ref{LemGenSmallMaxUnipParam} $(1)$ and the disjointness of 
Borel subgroups. Hence $G=B^{G}$, and $B<G$ is a full Frobenius group. 
\qed

\bigskip
We note that a connected $(\infty,r)$-homogeneous group of finite Morley rank, 
with $0<r<\infty$, is torsion-free by 
Facts \ref{pneqq} and \ref{FactToricvs.Unip}, and in particular 
Theorem \ref{TheoTorsionFreeHomegeneous} applies to such homogeneous 
connected locally$\o$ solvable$\o$ groups. 

Otherwise in the torsion free case all results of 
Section \ref{SectLocAnal} still apply, where all definable subgroups are 
connected. In this case Carter subgroups are conjugate by the same 
proof as in \cite{FreconConjCarter}. 

\subsection{The bounded exponent case}\label{SectBoundExpCase}

In presence of bounded exponent torsion the uniqueness theorems of 
Section \ref{SectionJaligLemma} can be applied in their most straightforward 
forms for dealing with generosity, 
as seen in Section \ref{SectionLocAnalandGen} already. 

\begin{lem}\label{LemmaNontrivialpUnipoExists}
Let $G$ be a locally$\o$ solvable$\o$ group of finite Morley rank such that 
$U_{p}(G)$ is nontrivial for some prime $p$. Then one of the following three 
cases occur. 
\begin{itemize}
\item[$(1)$]
Maximal $p$-unipotent subgroups are conjugate in $G\o$ and $N\o(U)$ is a 
generous Borel subgroup of unbounded exponent for any maximal $p$-unipotent 
subgroup $U$ (and in fact one may assume also $N\o(U)=UC\o(U)$). 
\item[$(2)$]
There is a maximal $p$-unipotent subgroup $U$ normalized but not centralized 
by a nontrivial $q$-torus $T$ for some (and in fact infinitely many) 
prime(s) $q\neq p$. Moreover $T$ is contained in a generous Carter subgroup of $G$. 
\item[$(3)$]
$N\o(U)$ is a Carter subgroup of bounded exponent for some 
maximal $p$-unipotent subgroup $U$. 
\end{itemize}
\end{lem}
\proof
First recall that $N\o(U)$ is a Borel subgroup of $G$ for any maximal 
$p$-unipotent subgroup $U$ of $G$ by Lemma \ref{NUnipHeavyBorel}. 

Assume case $(3)$ does not occur. This means that for any maximal 
$p$-unipotent subgroup $U$, $N\o(U)$ is not nilpotent of bounded 
exponent. By Fact \ref{FactUnipHeaviest}, this simply means that any 
such group $N\o(U)$ has unbounded exponent. 

If $UC\o(U)<N\o(U)$ for some maximal $p$-unipotent subgroup $U$, then 
Wagner's theorem \cite[Corollary 8]{MR1833472} gives a 
nontrivial $q$-torus in $N\o(U)$, for some prime $q\neq p$, acting 
nontrivially on $U$ (see for example \cite[Fact 2.5]{JaligotFrecon} and 
Zilber's field theorem \cite[\S9.1]{BorovikNesin(Book)94}). 
The fact that there are infinitely many primes $q$ occuring in the definable 
subgroup of the multiplicative group of the field of characteristic $p$ 
is due to \cite{MR1979004}. 
Then Fact \ref{DecentTorusInGenCarter} shows that we are in case $(2)$. 

This leaves us with the case in which $N\o(U)=UC\o(U)$ is a Borel 
subgroup of unbounded exponent for any maximal $p$-unipotent subgroup $U$. 

If $N\o(U)\cap N\o(U^{g})$ has a nontrivial connected component $X$ 
for some $g\in G$, then 
$N\o(X)$ is solvable by local$\o$ solvability$\o$ of $G$. As 
$N\o(U)=UC\o(U)$, $X$ centralizes a nontrivial $p$-unipotent 
subgroup of $U$ by Fact \ref{StructpSylSolvGps} $(2)$, and similarly 
a nontrivial $p$-unipotent subgroup of $N\o(U^{g})$. Now, as $N\o(X)$ 
is contained in a Borel subgroup,  
Lemma \ref{LemFusionBorel2} implies $N\o(U)=N\o(U^{g})$. Hence distinct 
conjugates of $N\o(U)$ have finite intersections. 

As $N\o(U)$ has unbounded exponent, these finite intersections cannot 
cover $N\o(U)$ generically by Lemma \ref{HUnifCovFinHBounded}. 
In particular they land in a (definable) non-generic subset of $N\o(U)$, and 
one concludes easily that $N\o(U)$ is generous in $G\o$ (see for instance 
\cite[Lemma 3.3]{CherlinJaligot2004}, bearing in mind that $N\o(U)$ is of 
finite index in its normalizer, as a Borel subgroup). 

We have thus $N\o(U)=UC\o(U)$ a generous Borel subgroup of unbounded 
exponent for any maximal $p$-unipotent subgroup $U$. 

Now let $U_{1}$ and $U_{2}$ be two maximal $p$-unipotent subgroups of $G$. 
By generosity of $N\o(U_{1})$ and \cite[Proposition 2.1]{Jaligot06}, 
a generic element $g$ of $G\o$ is in a conjugate of $N\o(U_{1})$, and in finitely many such. 
Similarly, $g$ is in a conjugate of $N\o(U_{2})$, say $N\o(U_{2})$ after conjugacy, 
and in finitely many such. Now $Z\o(U_{2})$ centralizes $g$ as 
$N\o(U_{2})=U_{2}C\o(U_{2})$. So it permutes naturally by conjugation 
the finitely many conjugates of $N\o(U_{1})$ containing $g$, and one can 
argue as in \cite[Fundamental Lemma 3.3]{Jaligot06}. 
By Fact \ref{ActionOnFiniteSets}, it fixes each of them, and in 
particular it normalizes a conjugate of $U_{1}$, say $U_{1}$ up to conjugacy. 
Hence $Z\o(U_{2})\leq N\o(U_{1})$, $Z\o(U_{2})\leq U_{p}(N\o(U_{1}))=U_{1}$, 
and $U_{1}=U_{2}$ by Theorem \ref{UniquenessLemma}. This shows that $U_{1}$ 
and $U_{2}$ are conjugate and completes our proof.  
\qed

\bigskip
First we note that this completes the proof of Theorem \ref{TheoProperCosetsNonGen}. 

The fact that $N\o(U)=UC\o(U)$ is stated between parentheses in case $(1)$ of 
Lemma \ref{LemmaNontrivialpUnipoExists} is to depreciate this aspect not 
true in the algebraic case. A conclusion closer to the algebraic case 
would be case $(2)$ combined with case $(1)$ without this aspect. But even in the well 
described context of \cite{CherlinJaligot2004} there are potentially Borel 
subgroups as in case $(2)$ but not as in case $(1)$ without this aspect 
(in sets of Borel subgroups usually denoted by 
$\mathfrak{B}$ in \cite{CherlinJaligot2004}). 

If the ambient group $G$ is locally$\o$ solvable in 
Lemma \ref{LemmaNontrivialpUnipoExists}, 
then one sees by the same argument, and using the results 
of Section \ref{SectionJaligLemma} adapted to the 
locally$\o$ solvable case, that Borel subgroups as in cases $(1)$ and $(3)$ 
have trivial intersections indeed, and are all generous. In particular 
a maximal $p$-unipotent subgroup $U$ as 
in case $(3)$ must satisfy $N\o(U)$ generous, and must be conjugate to 
one as in case $(1)$ if it exists. But in this case one also has 
$N\o(U)$ of unbounded exponent, and thus cases $(1)$ and $(3)$ are mutually exclusive. 
It follows also that cases $(2)$ and $(3)$ are mutually exclusive, and as 
cases $(1)$ and $(2)$ are obviously mutually exclusive all cases are pairwise 
mutually exclusive, and with a generous Carter subgroup in any case. 
One can summarize this as follows. 

\begin{lem}
Let $G$ be a locally$\o$ solvable group of finite Morley rank such that 
$U_{p}(G)$ is nontrivial for some prime $p$. Then exactly one of the following two 
cases occur. 
\begin{itemize}
\item[$(1)$]
Maximal $p$-unipotent subgroups are conjugate in $G\o$ and Borel subgroups 
of the form $N\o(U)$, for $U$ a maximal $p$-unipotent subgroup, are pairwise 
disjoint, generous, of the form $UC\o(U)$, and either of unbounded exponent or 
nilpotent of bounded exponent.  
\item[$(2)$]
There is a maximal $p$-unipotent subgroup $U$ normalized but not centralized 
by a nontrivial $q$-torus $T$ for some (and in fact infinitely many) 
prime(s) $q\neq p$. Moreover $T$ is contained in a generous Carter subgroup of $G$. 
\end{itemize}
\end{lem}

As in Section \ref{SectionTorsionFreeCase} one may wish to consider the 
$(p,\infty)$-homogeneous case for some prime $p$, or more 
generally the case in which all Borel subgroups are nilpotent but now 
of bounded exponent. In this case any Borel subgroup is a Carter subgroup 
of bounded exponent, and cases $(1)$ and $(2)$ of 
Lemma \ref{LemmaNontrivialpUnipoExists} cannot occur (recall that in 
case $(1)$ $N\o(U)$ has unbounded exponent). 
One can also see in this case that any two distinct Borel subgroups have a finite intersection, 
by using Corollary \ref{corUnicityBorelMax}. 

We continue with the mere presence of a nontrivial $p$-unipotent subgroup for some 
prime $p$. 

\begin{lem}\label{LemTUComBQ}
Let $G$ be a locally$\o$ solvable$\o$ group of finite Morley rank, 
$p$ and $q$ two primes (possibly the same). Assume a nontrivial 
$p$-unipotent subgroup $U$ of $G$ commutes with a nontrivial $q$-torus $T$ 
of $G$. Then there is a Borel subgroup $B$ of $G$ containing $T$, $U$, 
a Carter subgroup of $G$ (and $B$) containing $T$ and generous in $G$, and a maximal 
$p$-unipotent subgroup of $G$. 
\end{lem}
\proof
Let $Q$ be a Carter subgroup of $G$ containing $T$, which 
exists and is generous in $G$ by Fact \ref{DecentTorusInGenCarter}. 
We have $Q$ and $U$ in $N\o(T)$, and $N\o(T)\leq B$  for some 
Borel subgroup by local$\o$ solvability$\o$ of $G$. Now $B$ is the unique 
Borel subgroup of $G$ containing $U$ by the Uniqueness Theorem, 
here Corollary \ref{corUnicityBorelMax} or 
Corollary \ref{corUnicityBorelSuperMax}, 
and our claim follows. 
\qed

\begin{defi}
If $M$ is a proper definable subgroup of a group $G$ of finite Morley rank and $p$ a 
prime, we say that
\begin{itemize}
\item[$(1)$]
$M$ is {\em $p$-weakly embedded} in $G$ if $M$ has infinite $p$-subgroups and 
$M\cap M^{g}$ has no infinite $p$-subgroups for any $g$ in $G\setminus M$. 
\item[$(2)$]
$M$ is {\em $p$-strongly embedded} in $G$ if $M$ has nontrivial $p$-subgroups and 
$M\cap M^{g}$ has no nontrivial $p$-subgroups for any $g$ in $G\setminus M$. 
\end{itemize}
\end{defi}

Again the following remarks were obviously made around 
\cite{CherlinJaligot2004}, but they were not explicitely stated there to 
keep that paper not too long. 

\begin{lem}\label{LemBoundspEltsNU1NU2}
Let $G$ be a locally$\o$ solvable$\o$ group of finite Morley rank, $p$ a prime, 
$U_{1}$ and $U_{2}$ two distinct maximal $p$-unipotent subgroups of $G$ 
(which are then necessarily nontrivial). 
Then $p$-subgroups of $N(U_{1})\cap N(U_{2})$ are exceptional and 
have order at most $e(G)$. 
\end{lem}
\proof
By assumption, $N\o(U_{1})$ and $N\o(U_{2})$ are two distinct Borel subgroups 
of $G$. 

Assume toward a contradiction $N(U_{1})\cap N(U_{2})$ contains a $p$-subgroup 
$X$ with $C\o(X)$ solvable. We have then $C\o(X)\leq B$ for some Borel subgroup 
$B$. Notice that $C\o_{U_{1}}(X)$ and $C\o_{U_{2}}(X)$ are both nontrivial 
by Fact \ref{StructpSylSolvGps} $(2)$. Now Lemma \ref{LemFusionBorel2} 
implies $N\o(U_{1})=N\o(U_{2})$, a contradiction. 
\qed

\begin{cor}
Let $G$ be a locally$\o$ solvable$\o$ group of finite Morley rank with $G\o$ nonsolvable. 
Assume that for some prime $p$ maximal $p$-unipotent subgroups of $G$ are 
nontrivial. Then $N(U)$ is $p$-weakly embedded in $G$ for any such maximal 
$p$-unipotent subgroup $U$ of $G$, and $p$-strongly embedded 
whenever $G$ is locally$\o$ solvable. 
\end{cor}
\proof
Assume $N(U)\cap N(U^{g})$ has an infinite $p$-subgroup $S$ for some 
$g$ in $G$. We have $S\o\leq N\o(U)\cap N\o(U^{g})$, and as $S$ is infinite 
$S\o$ is infinite as well. By Fact \ref{StructpSylSolvGps} $(1)$, 
$S\o$ is a central product of a $p$-unipotent subgroup $V$ and a $p$-torus 
$T$, and one of the two factors is nontrivial by assumption. 
Now Lemma \ref{LemFusionBorel2} or 
Lemma \ref{LemBoundspEltsNU1NU2} gives in any case 
$N\o(U)=N\o(U^{g})$. Thus $g\in N(U)$. 

When $G$ is locally$\o$ solvable one proceeds similarly, but now the 
only exceptional $p$-element is the identity. 
\qed

\bigskip
We also observe that when a nontrivial $p$-unipotent subgroup commutes with a 
nontrivial $p$-torus, then a maximal $p$-torus commutes with a maximal 
$p$-unipotent subgroup by Lemma \ref{LemTUComBQ}. One can then 
build a $p$-weakly embedded subgroup as for the elimination of $2$-mixed type 
simple groups \cite{AltBorCher(Book)}. 
If $U$ is a definable $p$-unipotent subgroup of $G$, we denote by 
$U^{\perp}$ the definable connected subgroup $T_{p}(C(U))$, the subgroup of 
$C(U)$ generated by the definable hulls of its $p$-tori. 
By local$\o$ solvability$\o$, this group is solvable (for $U$ nontrivial). 
One sees easily that if $[U_{1},U_{2}]=1$, then $U_{1}^{\perp}=U_{2}^{\perp}$. 
Then one observes that the graph on the set of nontrivial $p$-unipotent subgroups, 
where adjacency is commutation, is not connected, as otherwise $U^{\perp}$ 
is independent of the choice of $U$, hence normal in 
$G$, and as it is nontrivial connected and solvable, $G\o$ is solvable by 
local$\o$ solvability$\o$, a contradiction to the assumption. 
Let $\cal C$ be a connected component of the graph. 
The group $G$ acts naturally on this graph. 
Let $M$ be the normalizer in $G$ of a connected component of $\cal C$. 
If $U\in \cal{C}$, then 
$M\leq N(U^{\perp})$. In particular $M\o$ is solvable, and $M$ has a 
unique maximal $p$-unipotent subgroup $U$. 
$M=N(U)=N(U^{\perp})$. Notice that $B=M\o$ is a Borel subgroup of $G$. 
But in any case one shows that $M$ is $p$-weakly embedded in $G$. 

With $p=2$ these notions will suffice to eliminate connected non-solvable mixed type 
locally$\o$ solvable$\o$ groups in \cite{DeloroJaligotII}, 
by methods and/or results from the simple case. 
For $p\neq 2$ Configuration \ref{BadConfig1} stands around. Pr\"ufer 
ranks will be controlled with the notion of strongly embedded subgroup. 

If one is not interested in conjugacy in Lemma \ref{LemmaNontrivialpUnipoExists} 
but merely in genericity, then one can notice that a connected 
locally$\o$ solvable$\o$ group with $U_{p}(G)$ nontrivial with no generous 
Borel subgroup must satisfy that $N\o(U)$ is a Carter subgroup of bounded exponent 
for each maximal $p$-unipotent subgroup $U$; otherwise $N\o(U)$ has unbounded 
exponent and one gets as in the proof of Lemma \ref{LemmaNontrivialpUnipoExists} 
either a nontrivial decent torus or $N\o(U)$ generous, a contradiction to the 
assumption. 

In particular, if the generic element of a connected 
locally$\o$ solvable$\o$ group $G$ of finite Morley rank is not in a connected 
nilpotent subgroup, then $G$ contains no decent tori 
(Fact \ref{DecentTorusInGenCarter}), contains nontrivial $p$-unipotent 
subgroups (Facts \ref{FactToricvs.Unip} and \ref{FactLiftingTorsion}), 
and $N\o(U)$ is a Carter subgroup of bounded exponent for each 
such maximal $p$-unipotent subgroup $U$, generically composed of 
exceptional elements by Lemma \ref{LemCarterWithHeavy}. 

\subsection{The toral homogeneous case}\label{SectionToralCase}

We shall now consider the case in which there is no bounded exponent subgroup, 
and more specifically the toral homogeneous case. Before studying 
this specific case precisely, we note that Carter subgroups are conjugate in any 
locally$\o$ solvable$\o$ group $G$ of finite Morley rank such that $d(G)<\infty$, 
by the same proof as in \cite{FreconConjCarter}. 

\begin{theo}\label{TheoToralHomogeneousNoException}
Let $G$ be a locally$\o$ solvable$\o$ group of finite Morley rank 
in which Borel subgroups are divisible abelian. Assume furthermore 
that nontrivial toral elements are not exceptional and that $G$ contains no involution. 
Then, either 
\begin{itemize}
\item[$(1)$]
$G\o$ is abelian, or 
\item[$(2)$]
$T<G\o$ is a full Frobenius group for some (any) Borel subgroup $T$. 
\end{itemize}
\end{theo}

We will use the following fact. 

\begin{fait}[Ali's Lemma]\label{NesinLemma}
Let $G$ be a group, $T_{1}$ and $T_{2}$ two disjoint subgroups, 
$x_{1}\in T_{1}\cap (N(T_{2})\setminus T_{2})$ and 
$x_{2}\in T_{2}\cap (N(T_{1})\setminus T_{1})$ satisfying 
$x_{1}T_{2}={x_{1}}^{T_{2}}$, 
$x_{2}T_{1}={x_{2}}^{T_{1}}$, and 
$(x_{1}^{2}){T_{2}}={({x_{1}^{2}})}^{T_{2}}$.
Then $x_{1}$ and $x_{1}^{2}$ are conjugate in $G$. 
In particular, if $x_{1}$ has prime order $p\neq 2$ and is central in $T_{1}$, 
$N(T_{1})$ controls fusion in $T_{1}$, and $N(T_{1})/T_{1}$ is finite, then 
some nontrivial prime divisor of $N(T_{1})/T_{1}$ divides $p-1$. 
\end{fait}
\proof
This is one of the essential contents of \cite{Nesin89-a}, already re-employed 
through the scope of \cite[Lemma 7.23]{CherlinJaligot2004}. 
By the fusion assumptions one can conjugate $x_{1}$ to $x_{1}x_{2}$ in 
$x_{1}T_{2}$, $x_{1}x_{2}$ to $x_{1}^{2}x_{2}$ in $x_{2}T_{1}$, 
and $x_{1}^{2}x_{2}$ to $x_{1}^{2}$ in $x_{1}^{2}T_{2}$, which yields 
the $G$-conjugacy of $x_{1}$ and $x_{1}^{2}$. 

For the second point we have now a nontrivial induced automorphism of $\<x_{1}\>$ in 
$N(T_{1})/T_{1}$, and the cyclic group 
$\<x_{1}\>$ of prime order $p$ has an automorphism group of order $p-1$. 
\qed

\bigskip
We now proceed to the proof of Theorem \ref{TheoToralHomogeneousNoException}. 

\bigskip
\proof
In order to prove Theorem \ref{TheoToralHomogeneousNoException} 
we consider now $G$ a connected 
locally$\o$ solvable$\o$ group, and fix a Borel subgroup $T$, which 
is divisible abelian by assumption. If $G$ is solvable, then $G=T$ and we 
are in case $(1)$. So we may assume $G$ not solvable. 

As in Lemma \ref{LemAbelianBorel}, any two Borel subgroups 
$T_{1}$ and $T_{2}$ must have a finite intersection $E$, and being 
a finite subgroup of a divisible abelian group $E$ must be toral if 
it is nontrivial. If $E\neq 1$, then $T_{1}$, $T_{2}\leq C\o(E)$, and one gets 
either $T_{1}=T_{2}$ when $C\o(E)$ is solvable, or a nontrivial exceptional 
toral subgroup otherwise, which is excluded by assumption. Hence 
distinct Borel subgroups are pairwise disjoint. As usual, each is generous, 
and they are all conjugate. 

As any element of $G$ also has an infinite centralizer, any such element 
must centralize an infinite abelian subgroup by 
Lemma \ref{LemGenSmallMaxUnipParam} $(1)$, 
and in particular normalizes the unique conjugate of $T$ containing it. 

This shows that $G=N(T)^{G}$. If $N(T)=T$, then $T$ is malnormal in $G$ by 
disjointness of pairwise distinct Borel subgroups, and $G=T^{G}$, and thus 
$T<G$ is a full Frobenius group as desired. 

Hence the analysis boils down to showing that $T$ is selfnormalizing. 
Assume on the contrary $T<N(T)$, and let $x$ be an element of order $p$ modulo $T$ 
for some prime $p$, which may be assumed to be itself a $p$-element of $G$ 
by Fact \ref{FactLiftingTorsion}, and in fact inside a $p$-torus. 

By conjugacy, one concludes that $T$ contains a maximal $p$-torus $T_{p}$ 
which is nontrivial. Now $x$ is in a conjugate $T_{p}^{g}$ of 
$T_{p}$ and $x^{p}\in T_{p}$. As $T^{g}\cap T=1$, as otherwise 
$T=T^{g}$ and $x\in T^{g}=T$, 
$x^{p}\in{T\cap T^{g}}=1$. 
For any element $y$ in $xT$, the definable hull $H(y)$ of $y$ contains 
also a $p$-element $y_{1}$ by Fact \ref{FactLiftingTorsion}, which similarly 
belongs to a maximal torus $T_{1}$ distinct from $T$. Now 
$C\o_{T}(y)\leq C\o_{T}(y_{1})\leq {T\cap T_{1}}$, and thus any 
element $y$ in $xT$ has a finite centralizer in $T$. Hence $y^{T}$ is generic in 
$xT$ for any $y$ in $xT$, and as the Morley degree is one one gets $xT=x^{T}$. 
Now $x$ normalizes $T$ and centralizes a nontrivial element $z$ in the elementary 
abelian $p$-subgroup of $T$. We have $z$ normalizing $T_{x}$, the torus containing $x$, 
without being inside, and similarly $zT_{x}=z^{T_{x}}$ 
(this is typical of \cite{Nesin89-a}. See also \cite[Lemma 7.19]{CherlinJaligot2004}). 
We are now in situation to apply 
Fact \ref{NesinLemma}. Noticing that $N(T)$ controls fusion in the torsion subgroup 
of $T$ by Corollary \ref{CorControlFusion}, this gives a contradiction by choosing 
for $p$ the smallest prime divisor $N(T)/T$. 
\qed

\bigskip
We note similarities between groups as in 
Theorem \ref{TheoToralHomogeneousNoException} $(2)$ with those of 
\cite{HoucineJaligot04}. These are far from 
being stable by \cite{JMN07}, but there are some hints for the existence of 
(at least partially) stable such groups, as envisionned in \cite[\S1]{JaligotPekin}. 

We also note that a reduction to Fact \ref{NesinLemma} yields involutions or triviality 
of Weyl groups in general in groups of finite Morley rank without non-trivial 
$p$-unipotent subgroups \cite{BurdgesCherlinSemisimpleTorsion}. In particular 
for the last paragraph of the proof of Theorem \ref{TheoToralHomogeneousNoException} 
we could have refered to this. 

With this triviality of Weyl groups in connected groups without involutions and 
without $p$-unipotent subgroups, one can give a general decomposition as 
in Theorem \ref{TheoToralHomogeneousNoException} without non-exceptionality 
assumption, with the Galois connection of Section \ref{SectionL0RvsL0R0}. 

\begin{theo}\label{TheoToralHomogeneous}
Let $G$ be a locally$\o$ solvable$\o$ group of finite Morley rank in which Borel 
subgroups are divisible abelian, and without involutions. Then $G\o$ has an abelian 
generous selfnormalizing Carter subgroup $T$ such that $G\o=T^{G\o}$. 
\end{theo}
\proof
This is similar to the proof of Theorem \ref{TheoToralHomogeneousNoException}, 
using Lemma \ref{LemAbelianBorel} for the generosity of (divisible abelian) 
Borel subgroups $T$. Notice also that such (conjugate) Borel subgroups $T$ are 
selfnormalizing by the above mentioned result of 
\cite{BurdgesCherlinSemisimpleTorsion}, or a more direct reduction to 
Fact \ref{NesinLemma} here, and cover $G\o$ by the same argument 
as in the proof of Theorem \ref{TheoToralHomogeneousNoException} again. 
\qed 

\bigskip
In general one cannot say much more in Theorem \ref{TheoToralHomogeneous}, 
except describing the full group $G\o$ by the graph of finite exceptional 
closed subsets of the divisible abelian Borel subgroup $T$ introduced at the end 
of Section \ref{SectionL0RvsL0R0} and delineated in 
Lemma \ref{extensionsminimales}. In fact, exceptional subsets of $T$ are in the 
divisible torsion subgroup of $T$, and in a finite subset of it by 
Lemma \ref{LemExcepEltsOfDecentTori}. One sees easily that closed exceptional 
subsets of $T$ correspond exactly to intersections of $T$ with distinct conjugates of $T$. 
If $(X_{0},\cdots , X_{k})$ is a maximal chain of exceptional closed subsets of $T$ in 
the graph of exceptional closed subsets of $T$ (i.e. with $(X_{i},X_{i+1})$ a minimal 
extension for each $i$), then $X_{0}=Z(G\o)$, $C\o(X_{i})=C_{G\o}(X_{i})$ for each $i$ 
(by a Frattini Argument following the conjugacy of generous Carter subgroups and 
triviality of the Weyl group $N(T)/T$ in $G\o$) and the center of this group is $X_{i}$, 
and each group 
$$C\o(X_{i})/X_{i}$$ 
also satisfies the assumptions of Theorem \ref{TheoToralHomogeneous}, 
with exception indices decreasing as $i$ increases. The last factor 
$C\o(X_{k})/X_{k}$ is as in Theorem \ref{TheoToralHomogeneousNoException} 
by Lemma \ref{LemKillExceptionsInTori}. 

Hence any group as in Theorem \ref{TheoToralHomogeneous} is entirely 
described as above by the finite graph of exceptional closed subsets of $T$. 
In particular the picture in Theorem \ref{TheoToralHomogeneous} 
looks like Configuration \ref{BadConfig1}, 
where all Borel subgroups involved are decent tori but potentially with more complexity 
involved in the finite graph of exceptional subsets of $T$. 
As for Theorem \ref{TheoToralHomogeneousNoException}, 
constructions of such abstract groups can be obtained as in \cite{HoucineJaligot04} 
with {\em any} finite graph for $T$ (compatible with the conditions of 
Lemma \ref{extensionsminimales}), similarly with a bad control on the complexity 
of their model theory by the general construction but perhaps with some stability 
if more care is taken. 

The case of groups as in Theorems \ref{TheoToralHomogeneousNoException} 
and \ref{TheoToralHomogeneous} but with involutions will be considered 
in \cite{DeloroJaligotII}, and eventually disappear by the analysis of this paper 
and the contents of \cite{Nesin89-a}. 

\subsection{Pr\"ufer ranks and strong embedding}\label{SectionPruferRanks}

In this final section we recast the dichotomy represented by 
Sections 6 and 7 of \cite{CherlinJaligot2004} in its actual content. 
If $S$ is an abelian $p$-group 
for some prime $p$ and $n$ a natural number, then we denote by $\Omega_{n}(S)$ 
the subgroup of $S$ generated by all elements of order $p^{n}$. 

\begin{theo}\label{TheoPruferRankspStrongEmbed}
Let $G$ be a connected nonsolvable locally$\o$ solvable$\o$ group of finite Morley rank 
of Pr\"ufer $p$-rank at least $2$ for some prime $p$, and fix a maximal $p$-torus 
$S$ of $G$. Assume that every proper definable connected 
subgroup containing $S$ is solvable, that elements of $S$ of order $p$ are 
not exceptional, and let 
$$B=\<C\o(s)~|~s\in {\Omega_{1}(S)\setminus \{1\}} \>.$$ 
Then either 
\begin{itemize}
\item[$(1)$]
$B<G$, in which case $B$ is a Borel subgroup of $G$, and moreover 
$N(B)$ is $p$-strongly embedded in $G$ assuming additionally that 
$U_{p}(C(s))=1$ for every element $s$ of order $p$ of $S$ and that $S$ is a 
Sylow $p$-subgroup of $N_{N(B)}(S)$, or  
\item[$(2)$]
$B=G$, in which case $S$ has Pr\"ufer $p$-rank $2$. 
\end{itemize}
\end{theo}

Let $M=N(B)$. As $B$ contains a Carter subgroup $Q$ of $G$ containing $S$, 
Facts \ref{FactConjGenCarter} and \ref{DecentTorusInGenCarter} 
and a Frattini Argument give 
$M=N(B)\subseteq BN(Q)$, and as $Q$ is almost selfnormalizing 
$B=M\o$. (With the same notation, this holds of course for an arbitrary $p$-torus $S$ 
in an arbitrary group $G$ of finite Morley rank.) 

Assume first 
$$B<G\leqno(1)$$

By assumption $B\leq B_{1}$ for some Borel subgroup $B_{1}$ of $G$. 
As $S\leq B \leq B_1$, Corollary \ref{CorBipgeneration} 
implies that $B=B_{1}$, and thus $B$ is a Borel subgroup of $G$. 

We adopt now the extra assumptions that $U_{p}(C(s))=1$ for every element 
$s$ of order $p$ of $S$ and that $S$ is a Sylow $p$-subgroup of 
$N_{N(B)}(S)$. 
We claim that $M=N(B)$ is $p$-strongly embedded in $G$ in this case by using 
a ``black hole" principle (a term going back to Harada) similar to the one 
used in \cite[\S2.2]{BurdgesCherlinJaligot07}, and already contained in 
\cite[Lemma 7.3]{CherlinJaligot2004}. 
We note that our additional assumptions imply that $S$ is a 
Sylow $p$-subgroup of $B$ 
(Facts \ref{HallGpConInConSolv} and \ref{StructpSylSolvGps}), and of 
$M$ as well, as $M=BN_{N(B)}(S)$ by a Frattini Argument. In particular 
$M/B$ has trivial Sylow $p$-subgroups by Fact \ref{FactLiftingTorsion}. 
 
Assume that $M\cap M^{g}$ contains an element $s$ of order 
$p$ for some $g$ in $G$. Notice that $s$ is actually in $B\cap B^g$, 
and $p$-toral. By connectedness and conjugacy of Sylow $p$-subgroups in 
connected solvable groups, the definition of $B$ implies that $C\o(s')\leq B$ 
for any element $s'$ of order $p$ of $B$. Similarly, $C\o(s')\leq B^{g}$ whenever 
$s'$ has order $p$ and is in $B^g$. By conjugacy in $B$ we may assume $s$ in $S$, 
and if $Q$ denotes a Carter subgroup of $G$ containing $S$ then 
$\Omega_1(S)\leq S\leq Q\leq C\o(s)\leq B\cap B^g$. 
By Lemma \ref{Bipgeneration} or Corollary \ref{CorBipgeneration} 
applied in $B$ and in $B^{g}$ we get 
$B^{g}={\<C\o(s)~|~s\in \Omega_{1}(S)\setminus \{1\} \>}=B$. 
Thus $g$ normalizes $B$, and is in $M$. 

Hence $M=N(B)$ is $p$-strongly embedded in $G$ under the two extra assumptions, 
and this proves clause $(1)$ of Theorem \ref{TheoPruferRankspStrongEmbed}. 

Now we pass to the second case 
$$B=G\leqno(2)$$
We will eventually show that clause $(2)$ of 
Theorem \ref{TheoPruferRankspStrongEmbed} holds by reworking the begining 
of Section 6 of \cite{CherlinJaligot2004}. We first put aside $p$-unipotent subgroups. 

\begin{lem}\label{LemAnyBorelContspPerp}
Any Borel subgroup containing a toral element of order $p$ has trivial 
$p$-unipotent subgroups. 
\end{lem}
\proof
Assuming the contrary, we may assume after conjugacy of decent tori that a 
Borel subgroup $L$ with $U_{p}(L)$ nontrivial contains an element $s$ of $S$ 
of order $p$. Then $U_{p}(C(s))$ is nontrivial by Fact \ref{StructpSylSolvGps}, 
contained in a unique Borel subgroup $B_{1}$ of $G$. (Actually $B_{1}=L$.) By 
Corollary \ref{corUnicityBorelSuperMax}, $B_{1}$ is the unique Borel subgroup 
containing any given nontrivial $p$-unipotent subgroup of $U_{p}(C(s))$. 
Now any element $s'$ of order $p$ of $S$ normalizes $U_{p}(C(s))$, and thus 
$U_{p}(C(s,s'))\neq 1$ by 
Fact \ref{StructpSylSolvGps}, and as $B_{1}$ is the unique Borel subgroup 
containing the latter group we get $C\o(s')\leq B_{1}$. This shows that $B\leq B_{1}$, 
a contradiction as $B=G$ is nonsolvable under the current assumption. 
\qed

\bigskip
In other words, nontrivial $p$-toral elements commute with no nontrivial $p$-unipotent 
subgroups. This can be stated more carefully as follows. 

\begin{cor}\label{CorLemAnyBorelContspPerp}
Any connected solvable subgroup $\<s\>$-invariant for some $p$-toral element 
$s$ of order $p$ has trivial $p$-unipotent subgroups. 
\end{cor}
\proof
Otherwise $s$ would normalize a nontrivial $p$-unipotent 
subgroup, and by Fact \ref{StructpSylSolvGps} it would 
centralize a nontrivial $p$-unipotent subgroup. 
\qed

Our assumption $(2)$ on $B$ yields similarly a property antisymmetric to 
the black hole principle implied by assumption $(1)$. 
Let $E$ denote the elementary abelian $p$-group $\Omega_1(S)$. 

\begin{lem}\label{LemDispatchCentInvBorels}
Let $E_1$ be a subgroup of $E$ of order at least $p^{2}$. Then for any proper 
definable connected 
subgroup $L$ there exists an element $s$ of order $p$ of $E_{1}$ such that 
$C\o(s)\nleq L$. In particular $G=\<C\o(s)~|~s\in E_{1}\setminus \{1\}\>$. 
\end{lem}
\proof
Assume on the contrary $C\o(s)\leq L$ for any element $s$ of 
order $p$ of $E_1$. 

We claim that $C\o(t)\leq L$ for any element $t$ of order $p$ of $E$. 
In fact, as $E_1\leq S\leq C\o(t)$, $C\o(t)$ is by Lemma \ref{Bipgeneration} 
generated by its subgroups of the form $C\o(t,s)$, with 
$s$ of order $p$ in $E_1$. As these groups are all contained in $L$ by assumption, 
our claim follows. 

Hence we have $B\leq L<G$. But under our current assumption 
$B=G$, and this is a contradiction. 

Our last claim follows, as proper definable connected subgroups of $G$ 
containing $S$ are solvable by assumption. 
\qed

\begin{cor}\label{CorCs<Cs}
There exists an element $s$ of order $p$ of $E$ such that $C\o(S)<C\o(s)$. 
\end{cor}
\proof
$C\o(S)$ is $S$-local$\o$, and thus solvable by local$\o$ solvability$\o$ of $G$ 
and Lemma \ref{LemGenCharLocSolvGps}. As $C\o(S)\leq C\o(s)$ for 
any element $s$ of order $p$ of $S$, it suffices to apply 
Lemma \ref{LemDispatchCentInvBorels}. 
\qed

\begin{lem}\label{LemExistsLdOp'L1}
There exists an element $s$ of order $p$ of $E$ such that 
$$d(O_{p'}(C\o(s)))\geq 1.$$
\end{lem}
\proof
Assume the contrary, and let $s$ be an arbitrary element of order $p$ of $E$. 
By our assumption that $d(O_{p'}(C\o(s)))\leq 0$, 
$O_{p'}(C\o(s))$ is trivial or a good torus by 
Lemma \ref{LemGenSmallMaxUnipParam}, and 
central in $C\o(s)$ by Fact \ref{ActionOnLessUnip} $(1)$. Notice that 
$U_{p}(C\o(s))=1$ by Lemma \ref{LemAnyBorelContspPerp}. As 
$C\o(s)/{O_{p'}(C\o(s))}$ is abelian by Fact \ref{FaitH/Op'HDivAbelian}, 
$C\o(s)$ is nilpotent. Now $S$ is central in $C\o(s)$ by Fact \ref{StrucNilpGroups} $(2)$. 
In particular $C\o(S)=C\o(s)$, and this holds for any element $s$ of order $p$ 
of $E$. We get a contradiction to Corollary \ref{CorCs<Cs}. 
\qed

\bigskip
It follows in particular from Lemma \ref{LemExistsLdOp'L1} 
that there exist definable connected subgroups $L$ 
containing $C\o(s)$ for some element $s$ of order $p$ of $E$ 
and such that $O_{p'}(L)$ is not a good torus. 
Choose then a unipotence parameter $\tilde q=(q,r)$ 
different from $(\infty,0)$ such that $r$ is maximal in the set of all 
$d_{q}(O_{p'}(L))$, where $L$ varies in the set of all definable connected 
solvable subgroups with the above property. 

Notice that there might exist several such maximal unipotence parameters 
$\tilde q$, maybe one for $q=\infty$ and several ones for $q$ prime, except 
for $q=p$ by Corollary \ref{CorLemAnyBorelContspPerp}. 

It will also be shortly and clearly visible below that the notion of 
maximality for $\tilde q$ is the same when $L$ varies in two smaller subsets 
of definable connected solvable subgroups 
containing $C\o(s)$ for some $s$ of order $p$ of $E$: the set of 
Borel subgroups with this property 
on the one hand, and exactly the finite set of subgroups of the form $C\o(s)$ on the other. 

\begin{lem}
Let $L$ be any definable connected solvable subgroup containing $C\o(s)$ 
for some element $s$ of order $p$ of $E$. Then $U_{\tilde q}(O_{p'}(L))$ 
is a normal definable connected nilpotent subgroup of $L$. 
\end{lem}
\proof
As $O_{p'}(L)$ is normal in $L$, it suffices to show that its definably 
characteristic subgroup $U_{\tilde q}(O_{p'}(L))$ is nilpotent. But the latter is in 
$F(O_{p'}(L))$ by Fact \ref{FactUnipHeaviest} and the maximality of $r$. 
\qed

\begin{cor}\label{CorEveryThingInFOp'}
Let $L$ be any definable connected solvable subgroup containing $C\o(s)$ 
for some element $s$ of order $p$ of $E$. Then any definable $\tilde q$-subgroup 
of $L$ without elements of order $p$ is in $U_{\tilde q}(F(O_{p'}(L)))$. 
\end{cor}
\proof
Let $U$ be such a subgroup. As $U_{p}(L)=1$ by Lemma \ref{LemAnyBorelContspPerp}, 
$L/O_{p'}(L)$ is (divisible) abelian by Fact \ref{FaitH/Op'HDivAbelian}, and 
thus $U\leq O_{p'}(L)$, and $U\leq U_{\tilde q}(O_{p'}(L))$. Now it 
suffices to apply the normality and the nilpotence of the latter. 
\qed

\bigskip
We now prove a version of the Uniqueness Theorem \ref{UniquenessLemma} 
with a combined action, more precisely where the assumption on unipotence 
degrees of centralizers is replaced by an assumption of invariance by 
a sufficiently ``large" $p$-toral 
subgroup. For this purpose we first note the following. 

\begin{lem}\label{LemActionUnipDegOnpPerpByE1}
Let $E_{1}$ be a subgroup of order at least $p^{2}$ of $E$, and $H$ a 
definable connected solvable $E_{1}$-invariant subgroup. Then 
$d_{q}(O_{p'}(H))\leq r$. 
\end{lem}
\proof
Assume toward a contradiction $r'>r$, where $r'$ denotes 
$d_{q}(O_{p'}(H))$. In this case $r$ is necessarily finite, and $q=\infty$. 
By Fact \ref{FactUnipHeaviest}, $U_{(\infty,r')}(O_{p'}(H))\leq F\o(O_{p'}(H))$, 
and this nontrivial definable $(\infty,r')$-subgroup is $E_{1}$-invariant. 
Fact \ref{BiGenerationWithpElts} gives an element $s$ of order $p$ 
in $E_{1}$ such that 
$$C_{U_{(\infty,r')}(O_{p'}(H))}(s)\neq 1.$$
But the latter is an $(\infty,r')$-group by 
Fact \ref{FaitDefGpAutpGp}. 
Now considering the definable connected solvable subgroup 
$C\o(s)$ gives a contradiction to the maximality of $r$, as $C\o(s)/O_{p'}(C\o(s))$ 
is (divisible) abelian as usual and the centralizer above is 
connected without elements of order $p$, and thus contained in $O_{p'}(C\o(s))$. 
\qed

\bigskip
As mentioned already around the definition of maximal parameters $\tilde q$, 
the same argument shows that $r$ is also exactly 
the maximum of the $d_{q}(O_{p'}(L))$ different from $0$, with $L$ varying 
in the set of {\em Borel} subgroups containing $C\o(s)$ for some element $s$ of 
order $p$ of $E$ 
(instead of all definable connected solvable subgroups $L$ with the same property), 
and similarly with $L$ varying in the set of subgroups $C\o(s)$ for some element $s$ of 
order $p$ of $E$. 

We now prove the specific version of the 
Uniqueness Theorem \ref{UniquenessLemma}. 

\begin{theo}\label{RelUniqLem}
Let $E_{1}$ be a subgroup of order at least $p^{2}$ of $E$. Then any 
$E_{1}$-invariant nontrivial definable $\tilde q$-subgroup without elements 
of order $p$ is contained in a unique maximal such. 
\end{theo}
\proof
Let $U_{1}$ be the $\tilde q$-subgroup under consideration. 
Fix $U$ a maximal $E_{1}$-invariant definable $\tilde q$-subgroup 
without elements of order $p$ containing $U_{1}$. 

Assume $V$ is another such subgroup, distinct from $U$, 
and chosen so as to maximize the rank of 
$U_{2}=U_{\tilde q}(U\cap V)$. As $1< U_{1}\leq U_{2}$, 
the subgroup $U_{2}$ is nontrivial. As $U_{2}$ is nilpotent, 
$N:=N\o(U_{2})$ is solvable by local$\o$ solvability$\o$ of $G$. 
Note that $U_{2}<U$, as otherwise $U=U_{2}\leq V$ and $U=V$ by maximality 
of $U$. Similarly $U_{2}<V$, as otherwise $V=U_{2}\leq U$ and $V=U$ 
by maximality of $V$. In particular by normalizer condition, 
Fact \ref{NormCond}, $U_{2}<U_{\tilde q}(N_{U}(U_{2}))$ and 
$U_{2}<U_{\tilde q}(N_{V}(U_{2}))$.

We claim that $d_{q}(O_{p'}(N))=r$. Actually 
$d_{q}(O_{p'}(N))\leq r$ by Lemma \ref{LemActionUnipDegOnpPerpByE1}, 
and as $O_{p'}(N)$ contains $U_{2}$ which is nontrivial and 
of unipotence degree $r$ in characteristic $q$ we get $d_{q}(O_{p'}(N))=r$. 

By Fact \ref{FactUnipHeaviest} and the fact that $r\geq 1$ 
we get $U_{\tilde q}(O_{p'}(N))\leq F\o(O_{p'}(N))$. 
In particular $U_{\tilde q}(O_{p'}(N))$ is 
nilpotent, and contained in a maximal definable $E_{1}$-invariant 
$\tilde q$-subgroup without elements of order $p$, say $\Gamma$. 
Notice that $N$, being $E_{1}$-invariant, satisfies $U_{p}(N)=1$, and 
$N/O_{p'}(N)$ is abelian as usual. 
Now $U_{1}\leq U_{2}<U_{\tilde q}(N_{U}(U_{2}))\leq \Gamma$, so our 
maximality assumption implies that $\Gamma=U$. In particular 
$U_{\tilde q}(N_{V}(U_{2}))\leq \Gamma =U$. But then 
$U_{2}<U_{\tilde q}(N_{V}(U_{2}))\leq U_{\tilde q}(U\cap V)=U_{2}$, 
a contradiction. 
\qed

\begin{cor}\label{CorUnipTheoAction}
Let $E_{1}$ be a subgroup of order at least $p^{2}$ of $E$. 
\begin{itemize}
\item[$(1)$]
If $U_{1}$ is a nontrivial $E_{1}$-invariant definable $\tilde q$-subgroup without 
elements of order $p$, then $U_{1}$ is contained in a unique maximal 
$E_{1}$-invariant definable connected solvable subgroup $B$. 
Furthermore $U_{\tilde q}(O_{p'}(B))$ is the unique 
maximal $E_{1}$-invariant definable $\tilde q$-subgroup without elements 
of order $p$ containing $U_{1}$, and, for any element $s$ of order $p$ 
of $E_{1}$ with a nontrivial centralizer in $U_{1}$, $C\o(s)\leq B$ and $B$ is 
a Borel subgroup of $G$. 
\item[$(2)$]
$U_{\tilde q}(O_{p'}(C\o(E_{1})))$ is trivial. 
\end{itemize}
\end{cor}
\proof
$(1)$. Assume $B_{1}$ and $B_{2}$ are two maximal $E_{1}$-invariant definable 
connected solvable subgroups 
containing $U_{1}$. We have $U_{p}(B_{1})=U_{p}(B_{2})=1$. Hence 
$B_{1}$ and $B_{2}$ are both abelian modulo their $O_{p'}$ subgroups. 

Let $U=U_{\tilde q}(O_{p'}(B_{1}\cap B_{2}))$. 
This group contains $U_{1}$ and is in particular nontrivial, and is $E_{1}$-invariant, 
as well as 
$U_{\tilde q}(O_{p'}(B_{1}))$ and $U_{\tilde q}(O_{p'}(B_{2}))$. 
Now all these three subgroups are contained in a (unique) common 
maximal $E_{1}$-invariant definable $\tilde q$-subgroup without elements 
of order $p$ by the Uniqueness Theorem \ref{RelUniqLem}, say $\tilde U$. 
Notice that $B_{1}=N\o(U_{\tilde q}(O_{p'}(B_{1})))$ 
and $B_{2}=N\o(U_{\tilde q}(O_{p'}(B_{2})))$ by maximality of 
$B_{1}$ and $B_{2}$. Now applying the normalizer condition, 
Fact \ref{NormCond}, in the subgroup $\tilde U$ without elements of order $p$ 
yields easily 
$U_{\tilde q}(O_{p'}(B_{1}))=\tilde U=U_{\tilde q}(O_{p'}(B_{2}))$. 
Now taking their common normalizers$\o$ yields $B_{1}=B_{2}$. 

Our next claim follows from the same argument. 

For the last claim, we note that there exists an element $s$ in $E_1$ of order $p$ 
such that $C_{U_{1}}(s)$ is nontrivial. By Fact \ref{FaitDefGpAutpGp} 
the latter is a $\tilde q$-group, and of course it is $E_{1}$-invariant. So 
the preceding uniqueness applies to $C_{U_{1}}(s)$, and as 
$C_{U_{1}}(s) \leq U_{1} \leq B$ we get that $B$ is the unique 
maximal $E_{1}$-invariant definable connected solvable subgroup 
containing $C_{U_{1}}(s)$. 
But $C_{U_{1}}(s) \leq C\o(s)\leq B_{s}$ for some Borel subgroup 
$B_{s}$ and $E_{1}\leq B_{s}$, so $B_{s}$ satisfies the same conditions as $B$, 
so $B_{s}\leq B$ and $B=B_{s}$ is a Borel subgroup of $G$. 

$(2)$. Suppose toward a contradiction 
$U:=U_{\tilde q}(O_{p'}(C\o(E_{1})))$ nontrivial. It is of course 
$E_{1}$-invariant. Recall that $Q$ is a fixed Carter 
subgroup of $G$ containing the maximal $p$-torus $S$. As $Q\leq C\o(E_{1})$, 
$Q$ normalizes the subgroup $U$. Now for any element $s$ of order $p$ or $E_{1}$ 
we have $UQ\leq C\o(s)$. 

As $E_{1}\leq Q$, any Borel subgroup containing $UQ$ is $E_{1}$-invariant, and by 
the first point there is a {\em unique} Borel subgroup containing $UQ$. 
Now $C\o(s)$ is necessarily contained in this unique 
Borel subgroup containing $UQ$, and this holds for any element $s$ of order 
$p$ of $E_{1}$. We get a contradiction 
to Lemma \ref{LemDispatchCentInvBorels}. 
\qed

\bigskip
We note that the proof of the second point in Corollary \ref{CorUnipTheoAction} 
actually shows that any definable connected subgroup containing $E_{1}$ and $U_{1}$ 
for some nontrivial $E_{1}$-invariant definable $\tilde q$-subgroup $U_{1}$ 
without elements of order $p$ is contained in a {\em unique} Borel subgroup 
of $G$. Furthermore with the notation of Corollary \ref{CorUnipTheoAction} $(1)$ 
we have in any case 
$N(U_{1})\cap N(E_{1})\leq N(U_{\tilde q}(O_{p'}(B)))=N(B)$. 

There are two possible ways to prove that the Pr\"ufer $p$-rank is $2$. 
One may use the Uniqueness Theorem \ref{RelUniqLem} provided by 
the local$\o$ solvability$\o$ of the ambient group, or use the general 
signalizer functor theory, which gives similar consequences in more general 
contexts. We explain now how to use the signalizer functor theory to get 
the bound on the Pr\"ufer $p$-rank, but we will rather continue the analysis 
with the Uniqueness Theorem \ref{RelUniqLem} which is closer in spirit 
to \cite[Lemma 6.1]{CherlinJaligot2004}, and our original proof anyway. 
It also gives much more information in the specific context under consideration, 
including when the Pr\"ufer $p$-rank is $2$, while the general 
signalizer functor theory just provides the bound. 

For $s$ a nontrivial element of $E$ we let 
$$\theta(s)=U_{\tilde q}(O_{p'}(C(s))).$$ 
If $t$ is another nontrivial element 
of $E$, then it normalizes the connected nilpotent $\tilde q$-group without 
elements of order $p$ $\theta(s)$, and by 
Facts \ref{FaitDefGpAutpGp} and \ref{FaitH/Op'HDivAbelian} 
$C_{\theta(s)}(t)\leq U_{\tilde q}(O_{p'}(C(t)))=\theta(t)$. Hence 
one has the two following properties:
\begin{itemize}
\item[$(1)$]
$\theta(s)^{g}=\theta(s^{g})$ for any $s$ in $E\setminus \{1\}$ and any $g$ in $G$. 
\item[$(2)$]
$\theta(s) \cap C_{G}(t) \leq \theta(t)$ for any $s$ and $t$ in $E\setminus \{1\}$. 
\end{itemize}
In the parlance of finite group theory one says that $\theta$ is an 
{\em $E$-signalizer functor} on $G$. In groups of finite Morley rank 
one says that $\theta$ is a {\em connected nilpotent} 
$E$-signalizer functor, as any $\theta(s)$ is connected (by definition) and 
nilpotent, which follows from Corollary \ref{CorEveryThingInFOp'}. 
When $E_{1}$ is a subgroup of $E$ one defines 
$$\theta(E_{1})=\<\theta(s)~|~s\in E_{1}\setminus \{1\}\>$$

In groups of finite Morley rank there is no Solvable Signalizer Functor Theorem 
available as in the finite case \cite[Chapter 15]{MR1264416} 
(see \cite{MR0323904, MR0297861, MR0417284, Bender75} for the story in the 
finite case). However Borovik imported from finite group theory a 
Nilpotent Signalizer Functor Theorem for groups of finite Morley rank 
\cite{Borovik90, Borovik95} \cite[Theorem B.30]{BorovikNesin(Book)94}, 
stated as follows in \cite[Theorem A.2]{Burdges03} (and which suffices by the 
unipotence theory of \cite{Burdges03} for which it has been designed 
originally). 

\begin{fait}
{\bf (Nilpotent Signalizer Functor Theorem)}
Let $G$ be a group of finite Morley rank, $p$ a prime, and $E\leq G$ a finite elementary 
abelian $p$-group of order at least $p^{3}$. Let $\theta$ be a connected nilpotent 
$E$-signalizer functor. Then $\theta(E)$ is nilpotent. Furthermore 
$\theta(E)=O_{p'}(\theta(E))$ and $\theta(s)=C_{\theta(E)}(s)$ for any 
$s$ in $E\setminus \{1\}$. 
\end{fait}

(From the finite group theory terminology one says that $\theta$ is {\em complete} 
when it satisfies the two properties of the last statement.) 

In our situation one has thus, assuming toward a contradiction the Pr\"ufer $p$-rank 
is at least $3$, that $\theta(E)$ is nilpotent. Notice that the 
definable connected subgroup $\theta(E)$ is nontrivial, 
as $\theta(s)$ is nontrivial at least 
for some $s$ by Facts \ref{FaitDefGpAutpGp} and \ref{BiGenerationWithpElts}. 
In particular $N\o(\theta(E))$ is solvable by 
local$\o$ solvability$\o$ of $G$. 

From this point on one can use arguments formally identical to those of 
\cite[\S6.2-6.3]{Borovik95} used there for dealing with 
``proper $2$-generated cores". 

If $E_{1}$ and $E_{2}$ are two subgroups of $E$ of order at least $p^{2}$, then 
for any $s$ in $E_{1}\setminus \{1\}$ one has 
$\theta(s)\leq \<C_{\theta(s)}(t)~|~t\in E_{2}\setminus \{1\}\>\leq 
\theta(E_{2})$ and thus $\theta(E_{1})=\theta(E_{2})$. 

In particular $\theta(E)=\theta(E_{1})$ for any subgroup $E_{1}$ of $E$ of order at 
least $p^{2}$. 

Now if $g$ in $G$ normalizes such a subgroup $E_{1}$, then 
$\theta(E)^{g}=\theta(E_{1})^{g}=\theta(E_{1}^{g})=
\theta(E_{1})=\theta(E)$ and thus $g\in N(\theta(E))$. 

Take now as in Lemma \ref{LemDispatchCentInvBorels} 
an element $s$ of order $p$ in $E$ such that $C\o(s)\nleq N\o(\theta(E))$. 

Then, still assuming $E$ of order at least $p^{3}$, there exists a subgroup $E_{2}$ of $E$ 
of order at least $p^{2}$ and disjoint from $\<s\>$. 
By Lemma \ref{Bipgeneration}, 
$$C\o(s)=\<C_{C\o(s)}(t)~|~t\in E_{2}\setminus \{1\}\>.$$
But now if $t$ is in $E_{2}$ as in the above equality, then $E_1:=\<s,t\>$ has order 
$p^{2}$ as $E_{2}$ is disjoint from $\<s\>$, hence 
$C_{C\o(s)}(t) \leq C(s,t)\leq N(\<s,t\>)=N(E_1)\leq N(\theta(E))$, 
and this shows that 
$C\o(s)\leq N\o(\theta(E))$. This is a contradiction, and as our only extra 
assumption was that the Pr\"ufer $p$-rank was at least $3$, it must be $2$. 

Anyway, we can get the bound similarly, by using more directly the 
Uniqueness Theorem \ref{RelUniqLem} here instead of the axiomatized signalizer 
functor machinery. Actually the proof below is the core of the proof of the 
Nilpotent Signalizer Functor Theorem, and the Uniqueness Theorem here gives a 
shortcut to the passage to a quotient for the induction in the general case. 

\begin{theo}
$S$ has Pr\"ufer $p$-rank $2$. 
\end{theo}
\proof
Assume towards a contradiction $E$ has order at least $p^{3}$. 

We then claim that there exists a {\em unique} maximal nontrivial $E$-invariant definable 
$\tilde q$-subgroup without elements of order $p$. 
Let $U_{1}$ and $U_{2}$ be two such subgroups. Then by 
Facts \ref{FaitDefGpAutpGp} and \ref{BiGenerationWithpElts} 
$C_{U_1}(E_1)$ and $C_{U_2}(E_2)$ are nontrivial $\tilde q$-subgroups for 
some subgroups $E_1$ and $E_2$ of $E$, each of index $p$ in $E$. Assuming 
$|E|\geq p^3$ gives then an element $s$ of order $p$ in $E_1 \cap E_2$. 
Now $C_{U_1}(s)$ and $C_{U_2}(s)$ are both nontrivial, and these are 
both $\tilde q$-subgroups by Fact \ref{FaitDefGpAutpGp}. 
Clearly they are both $E$-invariant, as $E$ 
centralizes $s$, and in $U_{\tilde q}(O_{p'}(C\o(s)))$ as usual, which is also 
$E$-invariant. Now the Uniqueness Theorem \ref{RelUniqLem} gives $U_{1}=U_{2}$, 
as desired. 

Hence there is a unique maximal $E$-invariant definable 
$\tilde q$-subgroup without elements of order $p$, say ``$\theta(E)$" in the 
notation of the signalizer functor theory. For the same reasons as mentioned 
above, Facts \ref{FaitDefGpAutpGp} and \ref{BiGenerationWithpElts}, it is nontrivial. 

Now by Facts \ref{FaitDefGpAutpGp} and \ref{BiGenerationWithpElts} again, 
$C_{\theta(E)}(E_{1})$ is a nontrivial definable $\tilde q$-subgroup of 
$\theta(E)$ for some subgroup $E_{1}$ of $E$ of index $p$. 
As $U_{p}(C\o(E_{1}))=1$, $C\o(E_{1})/O_{p'}(C\o(E_{1}))$ is abelian 
as usual, and the definable connected subgroup $C_{\theta(E)}(E_{1})$ is in 
$O_{p'}(C\o(E_{1}))$, and in $U_{\tilde q}(O_{p'}(C\o(E_{1})))$. 

But as $|E|\geq p^{3}$, $|E_{1}|\geq p^{2}$, and we get a contradiction to 
Corollary \ref{CorUnipTheoAction} $(2)$. 
\qed

\medskip
This proves clause $(2)$ of 
Theorem \ref{TheoPruferRankspStrongEmbed} and completes the proof 
of Theorem \ref{TheoPruferRankspStrongEmbed}.
\qed

\bigskip
We can also record informally some information gained along the proof of case 
$(2)$ of Theorem \ref{TheoPruferRankspStrongEmbed}, which can be 
compared to \cite[6.1-6.6]{CherlinJaligot2004}. We let 
$G$ and $S$ be as in case $(2)$ of Theorem \ref{TheoPruferRankspStrongEmbed}, 
and $Q$ be a Carter subgroup of $G$ containing $S$. Then $Q$ is contained in 
at least two distinct Borel subgroups of $G$ by 
Lemma \ref{LemDispatchCentInvBorels}, and in particular 
$Q$ is divisible abelian by Corollary \ref{corUnicityBorelSuperMax} 
and Proposition \ref{prop4.1OnInter}. Now there are unipotence parameters 
$\tilde q\neq (\infty,0)$ as in the proof of case $(2)$ of 
Theorem \ref{TheoPruferRankspStrongEmbed} (maybe one for $q=\infty$, 
several for $q$ prime, but none for $q=p$ by Lemma \ref{LemAnyBorelContspPerp}). 
All the results of the above analysis apply, now with $|\Omega_{1}(S)|=p^{2}$ 
necessarily. 

By Corollary \ref{CorUnipTheoAction}, 
$$U_{\tilde q}(O_{p'}(C\o(\Omega_{1}(S))))=1.$$ 
As $\Omega_{1}(S)$ has order $p^{2}$, it contains in particular 
$$\frac{p^{2}-1}{p-1}=p+1$$
pairwise noncolinear elements. It follows that there are at most 
$p+1$ nontrivial subgroups of the form $U_{\tilde q}(O_{p'}(C\o(s)))$ 
for some nontrivial element $s$ of order $p$ of $S$, and at most $p+1$ 
Borel subgroups $B$ containing $Q$ (actually $\Omega_{1}(S)$-invariant suffices 
as noticed after Corollary \ref{CorUnipTheoAction}) 
and such that $U_{\tilde q}(O_{p'}(B))\neq 1$. 
By Corollary \ref{CorUnipTheoAction}, any such Borel subgroup 
would contain $C\o(s)$ for any element 
$s$ of order $p$ of $S$ having a nontrivial centralizer in $U_{\tilde q}(O_{p'}(B))$, 
and $\Omega_{1}(S)$ has a trivial centralizer in 
$U_{\tilde q}(O_{p'}(B))$. 

The following corollary of Theorem \ref{TheoPruferRankspStrongEmbed} 
will be of crucial use in \cite{DeloroJaligotII} to get a bound on Pr\"ufer ranks. 

\begin{cor}\label{CorDichotomyHighPruferRk}
Let $G$ be a connected nonsolvable locally$\o$ solvable$\o$ group of finite Morley 
rank and of Pr\"ufer $p$-rank at least $2$ for some prime $p$, and fix a maximal $p$-torus 
$S$ of $G$. Let $X$ be a maximal exceptional (finite) subgroup of $S$ 
(as in Lemma \ref{LemKillExceptionsInTori}), 
$\overline{H}=C\o(X)/X$, $\overline{K}$ a minimal definable 
connected nonsolvable subgroup of $\overline{H}$ containing $\overline{S}$, 
and let 
$$\overline{B}=\<C\o_{\overline{K}}(\overline{s})~|~
\overline{s}\in \Omega_{1}(\overline{S}) \setminus \{ \overline{1} \}\>.$$ 
Then either  
\begin{itemize}
\item[$(1)$]
$\overline{B}<\overline{K}$, in which case $\overline{B}$ is a Borel subgroup 
of $\overline{K}$, and moreover $N_{\overline{K}}(\overline{B})$ 
is $p$-strongly embedded in $\overline{K}$ assuming additionally that 
$U_p(C_{\overline{K}}(\overline{s}))=1$ for every element 
$\overline{s}$ of order $p$ of $\overline{S}$ and that $\overline{S}$ is a 
Sylow $p$-subgroup of $N_{N_{\overline{K}}(\overline{B})}(\overline{S})$, or  
\item[$(2)$]
$\overline{B}=\overline{K}$, in which case 
$\overline{S}$, as well as $S$, has Pr\"{u}fer $p$-rank $2$. 
\end{itemize}
\end{cor}
\proof
It suffices to apply Theorem \ref{TheoPruferRankspStrongEmbed} in 
$\overline{K}$. We note that $\overline{S}$ and $S$ have the same 
Pr\"ufer $p$-rank, as $X$ is finite by Lemma \ref{LemCentNonSolvFinite}. 
\qed

\bigskip
Cases $(1)$ and $(2)$ of Theorem \ref{TheoPruferRankspStrongEmbed} and 
Corollary \ref{CorDichotomyHighPruferRk} correspond respectively to 
Sections 7 and 6 of \cite{CherlinJaligot2004} in presence of divisible torsion. 
The remaining analysis of both of these sections, as well as the treatment 
without the extra assumption for $p$-strong embedding in case $(1)$, 
will be considered in our separate paper on Weyl groups, mentioned 
already in Section \ref{SectionLocAnalandGen}. 

For $p=2$ case $(1)$ will entirely disappear in \cite{DeloroJaligotII} by an 
argument similar to the one used in \cite[Case I]{BurdgesCherlinJaligot07}.

\bibliographystyle{alpha}
\bibliography{biblio}

\end{document}